\documentclass{amsart}
\usepackage{amssymb}
\usepackage{amsmath}
\usepackage{pstricks}
\usepackage[all]{xy}
\usepackage{latexsym}

\newtheorem{thm}{Theorem}[subsection]
\newtheorem*{conv*}{Conventions}
\newtheorem{cor}[thm]{Corollary}
\newtheorem{lem}[thm]{Lemma}
\newtheorem{prop}[thm]{Proposition}
\theoremstyle{definition}
\newtheorem{defn}[thm]{Definition}
\newtheorem{notation}[thm]{Notation}
\newtheorem{remk}[thm]{Remark}
\newtheorem*{ack*}{Acknowledgments}
\newtheorem{ex}[thm]{Example}

\numberwithin{equation}{section}
 
\def \Ig {(I^2_g)_\bullet}
\def \bpt {\{*\}}
\def \scup {\uplus}
\def \dualized {dualized iterated integral}

\newcommand{\Triangle}{T}
\newcommand{\R}{\mathbb R}

\newcommand{\Z}{\mathbb Z}

\newcommand{\T}{\mathbb T}
\newcommand{\itint}{\int_\mathcal{C}}
\newcommand{\Co}{\mathcal{C}}
\newcommand{\Do}{\mathcal{D}}
\newcommand{\No}{\mathcal{N}}
\newcommand{\So}{\mathcal{S}}
\newcommand{\Fu}{\mathcal{F}}
\newcommand{\Ch}{\mathcal It}
\newcommand{\It}{\Ch}
\newcommand{\Chen}{\mathcal{C}\mathit{hen}}
\newcommand{\Sets}{\mathcal Sets}
\newcommand{\piwedge} {\pi^{\bigvee_{i=1}^{2g}S^1}}
\newcommand \pinch[1] {{\rm Pinch}_{#1}}
\newcommand \HHS[1] {HH^\bullet_{\Sigma^{#1}_\bullet}}
\newcommand \HH[2] {HH^{#2}_{\Sigma^{#1}_{\bullet}}}
\newcommand \CHS[1] {CH^\bullet_{\Sigma^{#1}_\bullet}}
\newcommand \CH[2] {CH^\bullet_{\Sigma^{#1}_{#2}}}
\newcommand \CHW[3] {CH^\bullet_{(\Sigma^{#1}\vee \Sigma^{#2})_{#3}}}
\newcommand \cop[1] {\delta^{#1}}
\newcommand{\eps}[1]{\varepsilon^{\Sigma^{#1}}}
\newcommand{\Map}{\operatorname{Map}}

\begin{document}

\title[A Chen model for mapping spaces and the surface product]{A Chen model for mapping spaces\\ and the surface product}

\author[G.~Ginot]{Gr\'egory Ginot}
\address{Gr\'egory Ginot, Institut  Math\'ematiques de Jussieu, Analyse Alg\'ebrique, UPMC - Universit\'e Pierre et Marie Curie, Case 82, 4, place Jussieu, F-75252 Paris Cedex 05}
\email{ginot@math.jussieu.fr}

\author[T.~Tradler]{Thomas~Tradler}
\address{Thomas Tradler, Department of Mathematics, College of Technology of the City University
of New York, 300 Jay Street, Brooklyn, NY 11201, USA}
\email{ttradler@citytech.cuny.edu}

\author[M.~Zeinalian]{Mahmoud~Zeinalian}
\address{Mahmoud Zeinalian, Department of Mathematics, C.W. Post Campus
of Long Island University, 720 Northern Boulevard, Brookville, NY
11548, USA} 
\email{mzeinalian@liu.edu}



\maketitle
\begin{abstract}
We develop a machinery of Chen iterated integrals for higher Hochschild complexes. These are complexes whose differentials are modeled on an arbitrary simplicial set much in the same way the ordinary Hochschild differential is modeled on the circle.  We use these to give algebraic models for general mapping spaces and define and study the surface product operation on the homology of mapping spaces of surfaces of all genera into a manifold. This is an analogue of the loop product in string topology. As an application, we show this product is homotopy invariant.  We prove Hochschild-Kostant-Rosenberg type theorems and use them to give explicit formulae for the surface product of odd spheres and Lie groups.

\bigskip

\noindent {\sc R\'esum\'e. }
 Dans cet article, on \'etend le formalisme des int\'egrales it\'er\'ees de Chen aux  complexes de Hochschild sup\'erieurs. Ces derniers sont des complexes de (co)cha\^{\i}nes model\'es sur un espace (simplicial) de la m\^eme mani\`ere que le complexe de Hochschild classique est model\'e sur le cercle. On en d\'eduit des mod\`eles alg\'ebriques pour les espaces fonctionnels que  l'on utilise pour  \'etudier le produit surfacique. Ce produit, d\'efini sur l'homologie des espaces de fonctions continues de surfaces (de genre quelconque) dans une vari\'et\'e, est  un analogue du produit de Chas-Sullivan sur les espaces de lacets en topologie des cordes. En particulier, on en d\'eduit que le produit surfacique est un invariant homotopique. On d\'emontre \'egalement un th\'eor\`eme du type Hochschild-Kostant-Rosenberg pour les complexes de Hochschild model\'es sur les surfaces qui permet d'obtenir des formules explicites pour le produit surfacique des sph\`eres de dimension impaires ainsi que pour les groupes de Lie.
\end{abstract}


\section{Introduction}
An element of the Hochschild chain complex $CH_\bullet(A, A)$ of an associative algebra $A$ is by definition an element in the multiple tensor product $A\otimes  \dots\otimes A$.  When defining the differential $D:CH_\bullet(A, A) \to CH_{\bullet-1}(A, A)$ however, it is instructive to picture this linear sequence of tensor products in a circular configuration, because the differential multiplies any two adjacent tensor factors starting from the beginning until the end and at the very end multiplies the last factor of the sequence with the first factor, as shown below.
$$ \begin{pspicture}(-.5,0)(7,3.5)
\psline[arrows=|->](3.2,2)(4.1,2)\rput(3.6, 2.3){$D$} \rput(5.5,2){$\sum $}
\rput(5.5,1.5){$_{\text{multiply }a_i\text{ to }a_{i+1}}$}
\rput(5.5,1.1){$_{\text{for all }i= 0,\dots, n,}$}
\rput(5.5,.7){$_{\text{ with }a_{n+1}\equiv a_0}$}
\pscircle[linestyle=dotted](1.5,2){.9}
\rput(1.5, .8){$a_0$} \rput(2, .9){$\otimes$} 
\rput(2.5, 1.2){$a_1$} \rput(2.7, 1.65){$\otimes$} 
\rput(1, 3.1){$\otimes$} 
\rput(1.5, 3.2){$a_i$} \rput(2, 3.1){$\otimes$}
\rput(.4, 1.5){$\otimes$} \rput(1.1,0.9){$\otimes$}
\rput(.75, 1.15){$a_n$}
\end{pspicture}
\begin{pspicture}(5,0)(8.5,3.5)
\pscircle[linestyle=dotted](6.5,2){.9}
\rput(5,1.8){$\pm$}
\rput(6.5, .8){$a_0$} \rput(7, .9){$\otimes$} 
\rput(7.5, 1.2){$a_1$} \rput(7.7, 1.65){$\otimes$} 
\rput(5.7, 3){$\otimes$} 
\rput(6.5, 3.15){$a_{i+1} \cdot a_i$} \rput(7.3, 3){$\otimes$}
\rput(6.1,.9){$\otimes$}
\rput(5.4, 1.5){$\otimes$} 
\rput(5.75, 1.15){$a_n$}
\end{pspicture}
$$
As it turns out this is not just a mnemonic device but rather an explanation of the fundamental connection between the Hochschild chain complex and the circle, which, for instance gives rise to the cyclic structure of the Hochschild chain complex and thus to cyclic homology, see~\cite{L}. This connection is also at the heart of the relationship between the Hochschild complex of the differential forms $\Omega^\bullet M$ on a manifold $M$, and the differential forms $\Omega^\bullet (LM)$ on the free loop space $LM$ of $M$, which is the space of smooth maps from the circle $S^1$ to the manifold $M$; see \cite{GJP}. At the core of this connection is the fact that the Hochschild complex is the underlying complex of a simplicial module whose simplicial structure is modelled on a particular simplicial model $S^1_\bullet$ of the circle. The principle behind this can be fruitfully used to construct new complexes whose module structure and differential are combinatorially governed by a given  simplicial set $X_\bullet$, much in the same way that the ordinary Hochschild complex is governed by $S^1_\bullet$; see \cite{P}. However carrying the construction to higher dimensional simplicial sets turns out to require associative and commutative algebras. The result of these constructions define for any (differential graded) commutative algebra $A$, any $A$-module $N$, and any pointed simplicial set $X_\bullet$, the \emph{(higher) Hochschild chain complex  $CH^{X_\bullet}_\bullet(A,N)$ of $A$ and $N$ over $X_\bullet$} as well as the Hochschild cochain complex $CH_{X_\bullet}^\bullet(A,N)$ over $X_\bullet$; see~\cite{P, G}. These Hochschild (co)chain complexes are functorial in all three of their variables $A$, $N$ and $ X_\bullet$.

\smallskip

The analogy with the usual Hochschild complex and its connection to the free loop space is in fact complete, because the  Hochschild complex $CH_\bullet^{X_\bullet}(\Omega^\bullet M, \Omega^\bullet M)$ over a simplicial set $X_\bullet$ provides an algebraic model of the differential forms on the mapping space $M^X=\{f:X \to M\}$, where $X=|X_\bullet|$ is the geometric realisation of $X_\bullet$. This is one of the main result of Section $2$ of this paper; see Section~\ref{section-qi}. 
The main tool to prove this result is a machinery of iterated integrals that we develop and use to obtain a quasi-isomorphicm $\Ch^{X_\bullet}:CH_\bullet^X(\Omega^\bullet M, \Omega^\bullet M) \to \Omega^\bullet({M^X})$, for any $k$-dimensional simplicial set $X_\bullet$ and  $k$-connected manifold $M$; see Proposition~\ref{qi}, and Corollary~\ref{qi*} for the dual statement. Further, for any simplicial set $X_\bullet$ and (differential graded) commutative algebra $A$, the Hochschild chain complex $CH^{X_\bullet}_\bullet(A,A)$ has a natural structure of a differential graded commutative algebra given by the shuffle product  $sh_{X_\bullet}$ (Proposition~\ref{P:shuffleinvariance}). We show that the iterated integral $\Ch^{X_\bullet}:(CH_\bullet^X(\Omega^\bullet M, \Omega^\bullet M), sh_{X_\bullet)} \to (\Omega^\bullet({M^X}), \wedge)$  is an algebra map sending  the shuffle product to the wedge product of differential forms on the mapping space; see Proposition~\ref{P:it-algebra}. 

\smallskip

Two important features of Hochschild (co)chain complexes over simplicial sets are their naturality in the simplicial set $X_\bullet$, and that two simplicial models of quasi-isomorphic spaces have naturally quasi-isomorphic Hochschild (co)chain complexes, see \cite{P}. In particular one (usually) obtains many different models to study $CH^{X_\bullet}_\bullet(A,N)$ for a given space $X=|X_\bullet |$. These facts are used, in Sections~\ref{S:surfaceproduct} and~\ref{S:surfacesymmetric}, to carry certain geometric and topological constructions over to the Hochschild complexes modeled on compact surfaces $\Sigma^g$ of genus $g$. 

The collection of  compact surfaces of any genus is naturally equipped with a product similar to the loop product of string topology~\cite{CS}, also see~\cite{S}.  
The idea behind this product, that we call the surface product, is shown in the following picture.
\begin{equation} \label{eq:rhoinrhoout}
\begin{pspicture}(1,0)(12.7,5)
\psline[arrows=->](4.1,2.2)(4.9,2.9) \rput(3.9,2.65){$wedge$} \psline[arrows=->](10.3,2.2)(9.5,2.9) \rput(10.5,2.65){$pinch$}
\psccurve(4.7,4)(5.4,5)(6.2,4.5)(7,5)(7.7,4)(7,3)(6.2,3.5)(5.4,3)
\pscurve(5.1,4.2)(5.4,3.8)(5.7,4.2) \pscurve(5.2,4)(5.4,4.1)(5.6,4)
\pscurve(6.7,4.2)(7,3.8)(7.3,4.2) \pscurve(6.8,4)(7,4.1)(7.2,4)
\pscircle(8.7,4){1}
\pscurve(8.4,4.2)(8.7,3.8)(9,4.2) \pscurve(8.5,4)(8.7,4.1)(8.9,4)
\psccurve(1,1)(1.7,2)(2.5,1.5)(3.3,2)(4,1)(3.3,0)(2.5,.5)(1.7,0)
\pscurve(1.4,1.2)(1.7,0.8)(2,1.2) \pscurve(1.5,1)(1.7,1.1)(1.9,1)
\pscurve(3,1.2)(3.3,0.8)(3.6,1.2) \pscurve(3.1,1)(3.3,1.1)(3.5,1)
\pscircle(5.3,1){1}
 \pscurve(5,1.2)(5.3,0.8)(5.6,1.2) \pscurve(5.1,1)(5.3,1.1)(5.5,1)
 \psccurve(8,1)(8.7,2)(9.5,1.5)(10.3,2)(11.1,1.5)(11.9,2)(12.7,1)(11.9,0)(11.1,.5)(10.3,0)(9.5,.5)(8.7,0)
 \pscurve(8.4,1.2)(8.7,0.8)(9,1.2) \pscurve(8.5,1)(8.7,1.1)(8.9,1)
 \pscurve(10,1.2)(10.3,0.8)(10.6,1.2) \pscurve(10.1,1)(10.3,1.1)(10.5,1)
 \pscurve(11.6,1.2)(11.9,0.8)(12.2,1.2) \pscurve(11.7,1)(11.9,1.1)(12.1,1)
 \pscurve(11.1,.5)(11.2,1)(11.1,1.5) \pscurve[linestyle=dashed](11.1,.5)(11,1)(11.1,1.5)
\end{pspicture}
\end{equation}
In Section~\ref{S:surfacemodel}, we describe an explicit simplicial model for the string topology type operation induced by the map $$ \Map(\Sigma^g,M)\times \Map(\Sigma^{h},M) \stackrel {\rho_{in}} \longleftarrow \Map(\Sigma^g \vee \Sigma^h,M) \stackrel {\rho_{out}}\longrightarrow \Map(\Sigma^{g+h},M) $$ coming from the above picture~\eqref{eq:rhoinrhoout}. More precisely, we obtain a surface product $\scup: H_\bullet(\Map(\Sigma^g,M))\otimes H_\bullet(\Map(\Sigma^{h},M))\to H_{\bullet+dim(M)}(\Map(\Sigma^{g+h},M))$, which is given by the composition of the Umkehr map $(\rho_{in})_!$ and the map induced by $\rho_{out}$,
\begin{multline*}
 \scup: H_\bullet(\Map(\Sigma^g,M))\otimes H_\bullet(\Map(\Sigma^{h},M)) \\ \stackrel{(\rho_{in})_!}{\longrightarrow} H_\bullet(\Map(\Sigma^{g}\vee \Sigma^{h},M))\stackrel{(\rho_{out})_*}{\longrightarrow} H_\bullet(\Map(\Sigma^{g+h},M)).
\end{multline*}
We prove that the surface product makes  $$\big(\bigoplus_g\mathbb{H}_\bullet(\Map(\Sigma^g,M)), \scup\big)= \big(\bigoplus_g H_{\bullet+dim(M)} (\Map (\Sigma^g, M )), \scup)$$ into an associative bigraded\footnote{we always use a cohomological grading convention in this paper hence the plus sign in our degree shifting} algebra  with $\mathbb{H}_\bullet(\Map(\Sigma^0,M))$ in its center; see Theorem~\ref{T:surfaceproduct} and Proposition~\ref{P:bimodule}.  The restriction of the surface product to genus zero (\emph{i.e.} spheres), $\mathbb{H}_{\bullet}(\Map(\Sigma^0,M))$, coincides with the Brane topology product defined by Sullivan and Voronov $\mathbb{H}_{\bullet}(\Map(S^2,M))^{\otimes 2} \to \mathbb{H}_{\bullet}(\Map(S^2,M))$ see~\cite{H, CV}. In fact, in these papers it is shown, that $H_{\bullet+\dim(M)} (M^{S^n})$ is an algebra over $H_\bullet(fD_{n+1})$, the homology of the framed $n$-disc operad; see also \cite{S} for related algebraic structures.

\smallskip

In Section~\ref{Surface Hochschild cup product}, we apply the machinery of (higher) Hochschild cochain complexes over simplicial sets to  give a fully algebraic description of the surface product. In fact, for positive genera, we define an associative cup product $\cup: \CHS{g}(A,B)\otimes \CHS{h}(A,B) \to \CHS{g+h}(A,B)$ for the Hochschild cochains over surfaces, where $B$ is a differential graded commutative and unital $A$-algebra, viewed as a symmetric bimodule; see Definition~\ref{D:cupg}. The construction of the cup-product is based on the fact that for any pointed simplicial sets $X_\bullet$ and $Y_\bullet$, the multiplication $B\otimes B\stackrel{\cdot}\to B$ induces a cochain map $\CH{g}{\bullet}(A,B)\otimes \CH{h}{\bullet}(A,B) \stackrel{\vee} \to\CHW{g}{h}{\bullet}(A,B)$ which can be composed with the pullback $\CHW{g}{h}{\bullet}(A,B)  \stackrel{\pinch{g,h}^*}\to\CH{\bullet}{i+j}(A,B)$ induced by a simplicial model $\pinch{g,h}:\Sigma^{g+h}_\bullet\to \Sigma^g_\bullet \vee \Sigma^h_\bullet$ for the map $pinch$ in Figure~\eqref{eq:rhoinrhoout}.

However, this product is initially only defined for surfaces of positive genera, and more work is required to include the genus zero case of the $2$-sphere in this framework. To this end, we first recall the cup product defined in \cite{G} for genus zero, and then define a left and right action, $\tilde{\cup}$, of $CH^{\Sigma^0_\bullet}_\bullet(A,B)$ on $CH^{\Sigma^g_\bullet}_\bullet(A,B)$. Taking advantage of the fact that one can choose different simplicial models for a given space, we show that  $\tilde{\cup}$ is after passing to homology, equivalent to operations similar to the cup-product but defined on different simplicial models for the sphere and surfaces of genus $g$; see Definition~\ref{D:edgewisecup} and Proposition~\ref{P:edgewisecup}.  Putting everything together, we  obtain a well-defined cup product $\cup:\big(\bigoplus_{g\geq 0}\HHS{g}(A,A) \big)^{\otimes 2} \to \bigoplus_{g\geq 0}\HHS{g}(A,A)$ for all genera on cohomology; see Theorem~\ref{T:cup}. More precisely, we prove that
 for a differential graded commutative algebra $(A,d_A)$,
\begin{itemize} \item[i)] the cup product makes $\bigoplus_{g\geq 0}\HHS{g}(A,A)$ into an associative algebra that is bigraded with respect to the total degree grading and  the genus of the surfaces and has a unit induced by the unit $1_A$ of $A$.  
\item[ii)]   $\HHS{0}(A,A)$ lies in the center of $\bigoplus_{g\geq 0}\HHS{g}(A,B)$.\end{itemize}
The cup product is functorial with respect to both arguments $A$ and $B$ (see Proposition~\ref{P:cupinvariance}). 

\smallskip

 The connection to topology is precise and the cup-product models the surface product. We prove in Theorem~\ref{T:surface=cup}, using rational homotopy techniques as in~\cite{FT}, that for a 2-connected compact manifold $M$, the (dualized) iterated integral $\Ch^{\Sigma^\bullet}:(\bigoplus_{g\geq 0 } \mathbb{H}_{-\bullet}(\Map(\Sigma^g,M)),\scup)  \to (\bigoplus_{g\geq 0 }\HH{g}{-\bullet}(\Omega,\Omega),\cup) $ is an  isomorphism of algebras. As a  corollary of this, it follows that the surface product $\scup$  is  homotopy invariant, meaning that if $M$ and $N$ are  2-connected compact manifolds with equal dimensions, and $i:M\to N$ is a homotopy equivalence, then $i_*: (\bigoplus_{g\geq 0}\mathbb{H}_\bullet(\Map(\Sigma^g,M)),\scup) \to  (\bigoplus_{g\geq 0}\mathbb{H}_\bullet(\Map(\Sigma^g,N)),\scup)$ is an isomorphism of algebras.  
 
\smallskip

Section~\ref{S:surfacesymmetric} is devoted to the Hochschild homology of (differential graded) symmetric algebras $(S(V),d)$ and a Hochschild-Kostant-Rosenberg type theorem. Recall that classically, when $d=0$, this theorem states that the Hochschild homology $HH_\bullet(S(V), S(V))$, thought of as the Hochschild homology of functions on the dual space $V^\ast$, can be identified with the space of K\"ahler differential $\Omega_{S(V)}$, which is the space of polynomial differential forms $\Omega_{S(V)}\cong \Omega^\bullet V^\ast = S(V) \otimes \Lambda(V) = S (V \oplus V[1])$.   This result, and its extension to the case $d\neq 0$, are main tools for computing Hochschild homology in algebra and topology; see~\cite{L, BV, FTV1, FTV} for instance.  Note that there is an obvious identification $ \Omega_{S(V)}= S (H_\bullet(S^1) \otimes V)$. Similarly, it is shown in \cite{P} that $HH^{S^2}_\bullet(S(V), S(V))= S(V \oplus V[2])= S(H_\bullet(S^2) \otimes V)$. 

We prove a similar kind of theorem for surfaces by showing that for a surface $\Sigma^g$ of genus $g$, there is an algebra isomorphism $$H_\bullet\big( S(H_\bullet(\Sigma^g)\otimes V), d^{\Sigma^g}\big) \stackrel{\sim}\longrightarrow HH^{\Sigma^g_\bullet}_\bullet(S(V), S(V))$$ where $d^{\Sigma^g}$ is a differential build out of the differential $d$ and the coalgebra structure of $H_\bullet(\Sigma^g)$, see Theorem~\ref{T:HKR}.  Indeed, the left-hand side in the above quasi-isomorphism  coincides with the Haefliger model~\cite{Hae, BS}. It is worth noticing that $d^{\Sigma^g}=0$ iff $d=0$. Further, there is a commutative diagram
$$\xymatrix{H_\bullet \big( S(H_\bullet(\mathop{\bigvee}\limits_{i=1}^{2g}S^1)\otimes V)\big) \ar[d]_{\cong } \ar[r]^{p} & H_\bullet \big( S(H_\bullet(\Sigma^g)\otimes V)\big) \ar[d]_{\cong}\ar[r]^{q}&  H_\bullet \big( S(H_\bullet(S^2)\otimes V) \big) \ar[d]_{\cong}\\ HH^{\bigvee_{i=1}^{2g}S^1}_\bullet(S(V),S(V))  \ar[r]^{(\bigvee_{i=1}^{2d} S^1 \hookrightarrow \Sigma^g)_\bullet}  & HH^{\Sigma^g}_\bullet(S(V),S(V)) \ar[r]^{(\Sigma^g \twoheadrightarrow S^2)_\bullet}  & HH^{S^2}_\bullet(S(V),S(V))} $$
where the horizontal maps $p,q$ are the algebra homomorphisms respectively induced by the homology maps $H_\bullet(\mathop{\bigvee}\limits_{i=1}^{2g}S^1)\otimes V)\to H_\bullet(\Sigma^g)$ and $H_\bullet(\Sigma^g)\to H_\bullet(S^2)$   obtained by the obvious inclusion and surjection of spaces and the lower maps are obtained by functoriality of Hochschild homology.

The main idea that is used to prove this result is, that if $X$ is a space obtained by attaching various spaces along attaching maps, then the Hochschild homology $HH^X_\bullet(A,A)$ can be computed by the Hochschild homology of the various pieces and the attaching maps via derived tensor products. For instance, using the fact that a genus $g$ surface can be obtained by suitably gluing a square along its boundary to a bouquet of circles, we show that there is an isomorphism 
$CH^{\Sigma^g}_\bullet(A,N) \cong CH_\bullet^{\bigvee_{i=1}^{2g}S^1}(A,N)\mathop{\otimes}\limits^{\mathbb{L}}_{CH_\bullet^{S^1}(A,A)} CH_\bullet^{I^2}(A,A)$ 
for any $A$-module $N$; see Section~\ref{S:derivedtensor}.
 With this tool in hand, we reduce the proof of the main theorem to appropriate statements about $CH_\bullet^{\bigvee_{i=1}^{2g}S^1}(A,M)$ and $CH_\bullet^{I^2}(A,A)$ for which the usual Hochschild-Kostant-Rosenberg Theorem and contractibility of the square $I^2$ can be used. . 

\smallskip

Furthermore, the Hochschild-Kostant-Rosenberg type theorem for surfaces allows us to explicitly compute the surface product  for odd spheres $S^{2n+1}$ and a Lie group $G$. The idea is that the differential graded algebras $\Omega^\bullet S^{2n+1}$ and $\Omega^\bullet G$, are both quasi-isomorphic to symmetric algebras with zero differentials; see Examples \ref{E:oddspheres} and \ref{E:LieGroups} in Section~\ref{S:Liegroups}.

\begin{ack*}
We would like to thank the referee for pointing out to us the Haefliger model.
We also would like to thank the Max-Planck Institut f\"ur Mathematik in Bonn, where this project initiated. We thank the Einstein Chair at the City University of New York, and the University Pierre et Marie Curie of Paris 6 for their invitations. The last two authors are partially supported by the NSF FRG grant DMS-0757245. The second author is partially supported by the PSC-CUNY grant PSCREG-40-228.
\end{ack*}

\begin{conv*}\quad
\begin{itemize}
\item[-] In this paper we use a cohomological grading for all our (co)homology groups and graded spaces, even when we use subscripts to denote the grading. In particular, all differentials are of degree $+1$, of the form $d:A^i\to A^{i+1}$ and the homology groups $H_i(X)$ of a space $X$ are concentrated in non-positive degree (unless otherwise stated).
\item[-] We follow the Koszul-Quillen sign convention: ``Whenever something of degree $p$ is moved past something of degree $q$ the sign $(-1)^{pq}$ accrues", see \cite{Q}. In particular, we will often write ``$\pm$" in the case that the sign is obtained by a permutation of elements of various degrees, following the Koszul-Quillen sign rule.
\end{itemize}
\end{conv*}
\section{Chen iterated integrals for Hochschild complexes over simplicial sets}

\subsection{Higher Hochschild chain complex over simplicial sets}

In this section, we define Chen iterated integral map for any simplicial set $Y$, and give explicit versions of it in the three examples of the circle $S^1$, the sphere $S^2$ and the torus $\mathbb T$, and combinations of those. In the next section, we will prove that (under suitable connectivity conditions) this gives in fact a quasi-isomorphism to the cochains of the corresponding mapping space $M^Y=\Map(Y,M)$. We start by recalling the Hochschild chain  complex of a differential graded, associative, commutative algebra $A$ over a simplicial set $Y_\bullet$ from T. Pirashvili, see \cite{P}.

\begin{defn}\label{D:Delta}
Denote by $\Delta$ the category whose objects are the ordered sets $[k]=\{0,1,\dots,k\}$, and morphisms $f:[k]\to [l]$ are non-decreasing maps $f(i)\geq f(j)$ for $i>j$. In particular, we have the morphisms $\delta_i:[k-1]\to[k], i=0,\dots, k$, which are injections, that miss $i$, and we have surjections $\sigma_j:[k+1]\to [k], i=0,\dots,k$, which send $j$ and $j+1$ to $j$.

A finite simplicial set $Y_\bullet$ is, by definition a contravariant functor from $\Delta$ to the category of finite sets $\Sets$, or written as a formula, $Y_\bullet:\Delta^{op}\to\Sets$. Denote by $Y_k=Y_\bullet([k])$, and call its elements simplicies. The image of $\delta_i$ under $Y_\bullet$ is denoted by $d_i:=Y_\bullet(\delta_i):Y_{k}\to Y_{k-1}$, for $i=0,\dots,k$, and is called the $i$th face. Similarly, $s_i:=Y_\bullet(\sigma_i):Y_{k}\to Y_{k+1}$, for $i=0,\dots,k$, is called the $i$th degeneracy. An element in $Y_k$ is called a degenerate simplex, if it is in the image of some $s_i$, otherwise it is called non-degenerate. 

We will mainly be interested in pointed finite simplicial sets. These are defined to be contravariant functors into the category $\Sets_*$ of pointed finite sets, $Y_\bullet:\Delta^{op}\to \Sets_*$. In this case, each $Y_k=Y_\bullet([k])$ has a preferred element called the basepoint, and all differentials $d_i$ and degeneracies $s_i$ preserve this basepoint.

Furthermore, we may extend these definitions to define simplicial sets, respectively pointed simplicial sets, by allowing the target of $Y_\bullet$ to be any (not necessarily finite) set.

Now, a morphism of (finite or not, pointed or not) simplicial sets is a natural transformation of functors $f_\bullet:X_\bullet\to Y_\bullet$. Note that $f_\bullet$ is given by a sequence of maps $f_k:X_k\to Y_k$ (preserving the basepoint in the pointed case), which commute with the faces $f_k d_i = d_i f_{k+1}$, and degeneracies $f_{k+1} s_i=s_i f_{k}$ for all $k\geq 0$ and $i=0,\dots, k$.
\end{defn}
\begin{defn}\label{D:Hoch}
Let $Y_\bullet:\Delta^{op}\to \Sets_*$ be a finite pointed simplicial set, and for $k\geq 0$, we set $y_k:=\# Y_k -1$, where $\#Y_k$ denotes the cardinal of the set $Y_k$. Furthermore, let $(A=\bigoplus_{i\in \Z} A^i, d, \cdot)$ be a differential graded, associative, commutative algebra, and $(M=\bigoplus_{i\in \Z}M^i,d_M)$ a differential graded module over $A$ (viewed as a symmetric bimodule). Then, the Hochschild chain complex of $A$ with values in $M$ over $Y_\bullet$ is defined as $CH_\bullet^{Y_\bullet}(A,M):=\bigoplus_{n\in \Z} CH_n^{Y_\bullet}(A,M)$, where $$ CH_n^{Y_\bullet}(A,M):=\bigoplus_{k\geq 0} (M\otimes A^{\otimes y_k})_{n+k} $$ is given by a sum of elements of total degree $n+k$. In order to define a differential $D$ on $CH^{Y_\bullet}_\bullet(A,M)$, we define morphisms $d_i:Y_k\to Y_{k-1}$, for  $i=0,\dots,k$ as follows. First note that for any map $f:Y_k\to Y_l$ of pointed sets, and for $m\otimes a_1\otimes \dots\otimes a_{y_k}\in M\otimes A^{\otimes y_k}$, we denote by $f_*:M\otimes A^{\otimes y_k}\to M\otimes A^{\otimes y_l}$, 
\begin{equation}\label{f_*}
 f_*(m\otimes a_1\otimes \dots\otimes a_{y_k})=(-1)^{\epsilon}\cdot n\otimes b_1\otimes \dots\otimes b_{y_l},
 \end{equation}
  where $b_{j}=\prod_{i\in f^{-1}(j)} a_i$ (or $b_j=1$ if $f^{-1}(j)=\emptyset$) for $j=0,\dots,y_{l}$, and $n=m\cdot \prod_{i\in f^{-1}(\text{basepoint}), i\neq \text{basepoint}}a_i$. The sign $\epsilon$ in equation \eqref{f_*} is determined by the usual Koszul sign rule of $(-1)^{|x|\cdot |y|}$ whenever $x$ moves across $y$. In particular, there are induced boundaries $(d_i)_*:CH_k^{Y_\bullet}(A,M)\to CH_{k-1}^{Y_\bullet}(A,M)$ and degeneracies $(s_j)_*:CH_k^{Y_\bullet}(A,M)\to CH_{k+1}^{Y_\bullet}(A,M)$, which we denote by abuse of notation again by $d_i$ and $s_j$. Using these, the differential $D:CH_\bullet^{Y_\bullet}(A,A)\to CH_\bullet^{Y_\bullet}(A,A)$ is defined by letting $D(a_0\otimes a_1\otimes \dots\otimes a_{y_k})$ be equal to 
\begin{equation*}
\sum_{i=0}^{y_k} (-1)^{k+\epsilon_i} a_0\otimes \dots\otimes d(a_i)\otimes \dots\otimes a_{y_k}+\sum_{i=0}^k (-1)^i d_i (a_0\otimes \dots\otimes a_{y_k}),
\end{equation*}
where $\epsilon_i$ is again given by the Koszul sign rule, \emph{i.e.}, $(-1)^\epsilon_i=(-1)^{|a_0|+\cdots+|a_{i-1}|}$. The simplicial conditions on $d_i$ imply that $D^2=0$.

More generally, if $Y_\bullet:\Delta^{op}\to \Sets$ is a finite (not necessarily pointed) simplicial set, we may still define $CH^{Y_\bullet}_\bullet(A,A):=\bigoplus_{n\in \Z} CH_n^{Y_\bullet}(A,A)$ via the same formula as above, $ CH_n^{Y_\bullet}(A,A):=\bigoplus_{k\geq 0} (A\otimes A^{\otimes y_k})_{n+k}$. Formula \eqref{f_*} again induces boundaries $d_i$ and degeneracies $s_i$, which produce a differential $D$ of square zero on $CH_\bullet^{Y_\bullet}(A,A)$ as above.
\end{defn}
\begin{remk}
Note that due to our grading convention, if $A$ is non graded (that is, concentrated in degree 0), then $HH_\bullet^{Y_{\bullet}}(A,A)$ is concentrated in non positive degrees. In particular our grading is  opposite of the one in~\cite{L}.
\end{remk}

Note that the equation~\eqref{f_*} also makes sense for any map of simplicial pointed sets $f:X_k\to Y_k$. 
Since $A$ is graded commutative and $M$ symmetric, $(f\circ g)_* =f_*\circ g_*$, hence $Y_\bullet \mapsto CH^{Y_\bullet}(A,M)$ is a functor from the category of finite pointed simplicial sets to the category of simplicial $k$-vector spaces, see~\cite{P}. If $M=A$, $CH^{Y_\bullet}(A,A)$ is a functor from the category of finite simplicial sets to the category of simplicial $k$-vector spaces.

\begin{defn}\label{normal}
Denote by $D_{k+1}$ the subspace of $CH^{Y_{k+1}}_\bullet(A,M)$ spanned by all degenerate objects, $D_{k+1}=Im((s_0)_*)+\dots+Im((s_k)_*)$. It is well-known (\cite[1.6.4, 1.6.5]{L}), that the $D_{k+1}$ form an acyclic subcomplex $D_\bullet$ of $CH^{Y_\bullet}_{\bullet}(A,M)$, which therefore implies that the projection $CH^{Y_\bullet}_{\bullet}(A,M)\to CH^{Y_\bullet}_{\bullet}(A,M)/D_\bullet$ is  a quasi-isomorphism. We denote this quotient by $NH^{Y_\bullet}_\bullet(A,M)$, and call it the normalized Hochschild complex of $A$ and $M$ with respect to $Y_\bullet$.
\end{defn}

The tensor product of  differential graded commutative algebras is naturally a differential graded commutative algebra. Thus, for any finite simplicial set $Y_\bullet$, $CH^{Y_\bullet}_\bullet(A,A)$ is a simplicial  differential graded commutative algebra, and further $CH^{Y_\bullet}_\bullet(A,M)$ is a (simplicial) module over $CH^{Y_\bullet}_\bullet(A,A)$. 
\begin{defn}\label{D:wedge}
Let $X_\bullet, Y_\bullet$, and $Z_\bullet$ be simplicial sets, and let $f_\bullet:Z_\bullet\to X_\bullet$ and $g_\bullet:Z_\bullet\to Y_\bullet$ be maps of simplicial sets. We define the wedge $W_\bullet=X_\bullet\cup_{Z_\bullet} Y_\bullet$ of $X_\bullet$ and $Y_\bullet$ along $Z_\bullet$ as the simplicial space given by $W_k= (X_k \cup Y_k)/\sim$, where $\sim$ identifies $f_k(z)=g_k(z)$ for all $z\in Z_k$. The face maps are defined as $d^{W_\bullet}_i(x)=d^{X_\bullet}_i(x), d^{W_\bullet}_i(y)=d^{Y_\bullet}_i(y)$ and the degeneracies are $s^{W_\bullet}_i(x)=s^{X_\bullet}_i(x), s^{W_\bullet}_i(y)=s^{Y_\bullet}_i(y)$ for any $x\in X_k \hookrightarrow W_k$ and $y\in Y_k \hookrightarrow W_k$. It is clear that $W_\bullet$ is well-defined and there are simplicial maps $X_\bullet\stackrel{i_\bullet}\to W_\bullet$ and $Y_\bullet\stackrel{j_\bullet}\to W_\bullet$.

If $X_\bullet$ is a pointed simplicial set, then we can make $W_\bullet$ into a pointed simplicial set by declaring the basepoint to be the one induced from the inclusion $X_\bullet\to W_\bullet$. (Note that this is in particular the case, when $X_\bullet, Y_\bullet, Z_\bullet, f_\bullet$ and $g_\bullet$ are in the pointed setting.)
\end{defn}

\begin{lem}\label{L:pushout}
There is a map of Hochschild chain complexes 
$$  CH^{X_\bullet}_\bullet(A,M)\otimes_{CH^{Z_\bullet}_\bullet(A,A)} CH^{Y_\bullet}_\bullet(A,A) \to CH^{W_\bullet}_\bullet(A,M). $$
If $Z_\bullet$ injects into either $Z_\bullet\stackrel {f_\bullet}\to X_\bullet$ or $Z_\bullet\stackrel{g_\bullet}\to Y_\bullet$, then this map is in fact an isomorphism of simplicial vector spaces.
\end{lem}
\noindent The tensor product in Lemma~\ref{L:pushout} is the tensor product of (simplicial) modules over the simplicial differential graded commutative algebra $CH^{Z_\bullet}_\bullet(A,A)$.
\begin{proof}
Using the functoriality of the Hochschild chain complex \cite{P}, there is a commutative diagram,
\begin{equation}\label{eq:pushout}
\xymatrix{
 CH_\bullet^{Z_\bullet} (A,A) \ar[r]^{f_*} \ar[d]_{g_*}  & CH_\bullet^{X_\bullet}(A,A) \ar[d]^{i^*} \\
 CH_\bullet^{Y_\bullet} (A,A)  \ar[r]^{j_*} &  CH_\bullet^{W_\bullet} (A,A) }
\end{equation}
which induces the claimed map. If $Z_\bullet$ injects for example into $X_\bullet$, then the tensor product $A^{\otimes x_k+1}\otimes_{A^{\otimes z_k+1}} A^{\otimes y_k+1}$ is isomorphic to $A^{\otimes (x_k-z_k)+y_k+1}$, which gives the isomorphism of the Hochschild complexes.
\end{proof}

\subsection{Definition of Chen iterated integrals}\label{S:Chen-map}
Assume now, that $M$ is a compact, oriented manifold, and denote by $\Omega=\Omega^\bullet (M)$ the space of differential forms on $M$.  For any simplicial set $Y_\bullet$, we now define the space of Chen iterated integrals $\Chen (M^Y)$ of the mapping space $M^Y=\Map(Y,M)$, and relate it to the  Hochschild complex over $Y_\bullet$ from the previous section.
\begin{defn}\label{Chen-map}
Denote by $Y:=|Y_\bullet|= \coprod \Delta^\bullet \times Y_\bullet /\sim$ the geometric realization of $Y_\bullet$. Furthermore, let $S_\bullet(Y)$  be the simplicial set associated to $Y$, i.e. $S_k(Y):=\Map(\Delta^k,Y)$ is the mapping space from the standard $k$-simplex $\Delta^k=\{0\leq t_1\leq \dots\leq t_k\leq 1\}$ to $Y$. By adjunction, there is the canonical simplicial map $\eta: Y_\bullet \to S_\bullet (Y) $, which is given for $i\in Y_k=\{0,\dots,y_k\}$ by maps $\eta(i):\Delta^k \to Y$ in the following way,
\begin{equation*}
\eta (i)(t_1\leq\cdots\leq t_k):= [(t_1\leq\cdots\leq t_k)\times \{i\}]\in \left(\coprod \Delta^\bullet \times Y_\bullet /\sim\right) =Y.
\end{equation*}
Here we identify $0$ with the base point of $Y_k$.

Denote by $M^Y=\Map(Y,M)$ the space of continuous maps from $Y$ to $M$, which are smooth on the interior of each simplex $Image(\eta(i))\subset Y$.
Recall from Chen \cite[Definition 1.2.1]{C}, that to give a differentiable structure on $M^Y$, we need to specify a set of plots $\phi:U\to M^Y$, where $U\subset \R^n$ for some $n$, which include the constant maps, and are closed under smooth transformations and open coverings. To this end, we denote the adjoint of a map $\phi:U\to M^Y$ by $\phi_\sharp:U\times Y\to M$. Then, define the plots of $M^Y$ to consist those maps $\phi:U\to M^Y$, for which $\phi_\sharp:U\times Y\to M$ is continuous on $U\times Y$, and smooth on the restriction to the interior of each simplex of $Y$, i.e. $\phi_\sharp|_{U\times (\text{simplex of }Y)^\circ}$ is smooth.

Following \cite[Definition 1.2.2]{C}, a $p$-form $\omega\in \Omega^p(M^Y)$ on $M^Y$ is given by a $p$-form $\omega_\phi\in \Omega^p(U)$ for each plot $\phi:U\to M^Y$, which is invariant with respect to smooth transformations of the domain. Let us recall from~\cite{C} that a smooth transformation between two plots $\phi:U\to M^Y$ and $\psi:V\to M^Y$ is a smooth map $\theta: U\to V$ such that the following diagram $$\xymatrix{U\times Y \ar[rr]^{\phi_\sharp} \ar[d]_{\theta} && M \\
V\times Y \ar[urr]_{\psi_\sharp}&& }$$ commutes. The invariance of a $p$-form means that $\theta^{*}(\omega_\psi) =\omega_\phi$. The differential, wedge product and pullback of forms are all defined plotwise. We will now define certain forms on $M^Y$, which we call (generalized) iterated integrals.

To this end, assume that $\phi:U\to M^Y$ is a plot of $M^Y$, and we are given $y_k+1=\# Y_k$ many forms on $M$, $a_0,\dots,a_{y_k}\in \Omega=\Omega^\bullet(M)$. We define $\rho_\phi:=ev\circ(\phi\times id)$,
\begin{equation}\label{rho_phi}
\rho_\phi: U\times \Delta^k \stackrel{\phi\times id} \longrightarrow M^Y \times \Delta^k\stackrel {ev} \longrightarrow M^{y_k+1},
\end{equation}
where $ev$ is defined as the evaluation map,
\begin{multline}\label{eval}
ev(\gamma:  Y\to M, t_1\leq\cdots\leq t_k)=\left( \big(\gamma\circ \eta(0)\big) (t_1\leq\cdots\leq t_k) ,\right.\\
 \left.\big(\gamma\circ \eta(1)\big) (t_1\leq\cdots\leq t_k),\cdots, \big(\gamma\circ \eta(y_k)\big) (t_1\leq\cdots\leq t_k)\right).
\end{multline}
Now, the pullback $(\rho_\phi)^*(a_0\otimes\dots\otimes a_{y_k})\in \Omega^\bullet (U\times \Delta^k)$, may be integrated along the fiber $\Delta^k$, and is denoted by
\begin{equation*}
\left(\itint a_0\dots a_{y_k} \right)_{\phi}:=\int_{\Delta^k} (\rho_\phi)^*(a_0\otimes\dots\otimes a_{y_k})\quad \in \Omega^\bullet (U).
\end{equation*}
The resulting $p=(\sum_i deg(a_i) -k)$-form $\itint a_0\dots a_{y_k}\in \Omega^p(M^Y)$ is called the (generalized) iterated integral of $a_0,\dots, a_{y_k}$. Here, we used the symbol $\itint$ instead of $\int$, since our notation differs slightly from the usual one, where iterated integrals refer to the integration over the interior of a path, see also example \ref{S^1} below.

The subspace of the space of De Rham forms $\Omega^\bullet(M^Y)$ generated by all iterated integrals is denoted by $\Chen(M^Y)$. In short, we may picture an iterated integral as the pullback composed with the integration along the fiber $\Delta^k$ of a form in $M^{y_k+1}$,
\begin{equation*}
\xymatrix{
  M^Y\times \Delta^k \ar[r]^{ev} \ar[d]_{\int_{\Delta^k}} & M^{y_k+1} \\ M^Y& }
\end{equation*}
\end{defn}

We now use the above definition to obtain a chain map $\Ch^{Y_\bullet}:CH_\bullet^{Y_\bullet}(\Omega,\Omega)\to \Chen(M^Y)$. In sections \ref{section-qi} and \ref{section-product}, we will see that $\Ch^{Y_\bullet}$ is in fact a quasi-isomorphism and an algebra map. In detail, for homogeneous elements $a_0\otimes\dots\otimes a_{y_k}\in CH_\bullet^{Y_\bullet}(\Omega,\Omega)$, we define \begin{equation} \label{eq:chen-map}\It^{Y_\bullet}(a_0\otimes\dots\otimes a_{y_k}):=\itint a_0\dots a_{y_k}.\end{equation}
\begin{lem}
$\Ch^{Y_\bullet}:CH_\bullet^{Y_\bullet}(\Omega,\Omega) \to \Chen (M^Y)$ is a chain map.
\end{lem}
\begin{proof} Since, $\itint a_0\dots a_{y_k}$ is a $(\sum_i deg(a_i) -k)$-form,  $\Ch^{Y_\bullet}$ is a degree $0$ linear map.
We need to show that $\Ch^{Y_\bullet}$ respects the differentials. Using the definitions together with Stokes theorem for integration along a fiber, $$ d\left(\int_{\Delta^k} \omega\right)=\int_{\Delta^k} d\omega \pm \int_{\partial \Delta^k}\omega, $$ we see that for some plot $\phi:U\to M^Y$, $d\left(\Ch^{Y_\bullet}(a_0\otimes a_1\otimes\cdots \otimes a_{y_k})\right)_{\phi}$ is equal to
\begin{multline*}
 d\int_{\Delta^k} (\rho_\phi)^*(a_0\otimes a_1\otimes\cdots \otimes a_{y_k})\\
= \int_{\Delta^k} d \Big((\rho_\phi)^*(a_0\otimes a_1\otimes\cdots \otimes a_{y_k}) \Big)\pm \int_{ \partial(\Delta^{k})} (\rho_\phi)^*(a_0\otimes a_1\otimes\cdots \otimes a_{y_k})\\
= \int_{\Delta^k}  (\rho_\phi)^*\left(\sum_{j=0}^{y_k} a_0\otimes \cdots\otimes d(a_j)\otimes\cdots \otimes a_{y_k}\right) \quad\quad\quad \quad\quad\quad \quad\quad\quad \quad\quad \\ \pm \int_{\Delta^{k-1}} (\rho_\phi)^* \left( \sum_{j=0}^{k} b^j_0\otimes\cdots \otimes b^j_{y_{k-1}}\right),
\end{multline*}
where, for the boundary component, with $t_1\leq \cdots\leq t_j=t_{j+1}\leq \cdots\leq t_k$, $b^j_i$ is the product of all the $a_l$'s for which $l$ satisfies $d_i(l)=j$. Thus, we recover precisely $\Ch^{Y_\bullet}(D(a_0\otimes a_1\otimes\cdots \otimes a_{y_k}))$. 
\end{proof}

Recall also from \cite[Definition 1.3.2]{C} that the space of smooth $p$-chains $C_p(M^Y)$ on $M^Y$ is the space generated by plots of the form $\sigma:\Delta^p\to M^Y$. We denote by $C_\bullet(M^Y)=\bigoplus_p C_p(M^Y)$, and give it a differential in the usual way. It follows from \cite[Lemma 1.3.1]{C}, that the induced homology is the usual one, $H(C_\bullet(M^Y))=H_\bullet(M^Y)$ (that is, the singular homology of $M^Y$). The canonical chain map $\Omega^\bullet(M^Y)\otimes C_\bullet(M^Y)\to \R$, given by $<\omega,\sigma>:=\int_{\Delta^p} \omega_\sigma$, induces a similar chain map on the space of iterated integrals, $\Chen(M^Y)\otimes C_\bullet(M^Y)\to \R$.

Finally, we remark that the construction is clearly functorial in $Y_\bullet$. Thus it extends by limits to the case where $Y_\bullet$ is non necessarily finite. In particular, this allows us to define a Chen iterated integral map for a topological spaces that are not simplicial. 
\begin{defn}\label{Singular}
Let $X$ be a (pointed) topological space. Then we have the (pointed) simplicial set $S_\bullet(X)$. By definition the Hochschild chain complex of $A$ over $S_\bullet(X)$ with value in an $A$-module $N$ is, in (external) degree $k$ the limit 
$$CH^{S_k(X)}_\bullet(A,N):= \lim_{i_+\to S_k(X)} N\otimes A^{\otimes i} $$ over all $i_+$ where $i_+=\{0,\cdots,i\}$, with $0$ for base point.
If $X=|Y_\bullet|$ is the realization of a simplicial set $Y_\bullet$, then the natural map $\eta:Y_\bullet \to S_\bullet(|Y_\bullet|)\cong S_\bullet(X)$ (see Definition~\ref{Chen-map}) induces a natural quasi-isomorphism $CH_\bullet^{Y_\bullet}(A,N) \stackrel{\eta_*} \longrightarrow CH_\bullet^{S_\bullet(X)}(A,N)$, see~\cite{P}.
Let us define  $\Ch^X: CH^{S_\bullet(X)}_k(\Omega,\Omega)\to \Omega^\bullet(\Map(X,M)) $. It is enough to define, for all $k\geq 0$, natural maps $\Ch_\beta:\Omega \otimes \Omega^{\otimes i} \to \Omega^\bullet(\Map(X,M))$ for all $\beta: i_+\to S_k(X)$. We define, for each plot $\phi:U\to \Map(X,M)$, and $a_0,\dots,a_i\in \Omega=\Omega^\bullet(M)$,
$$\Ch_{\beta}(a_0\otimes \cdots \otimes a_i)_\phi=\int_{ \Delta^k} (\rho_{\phi,\beta})^*(a_0\otimes \cdots \otimes a_i), $$
where $\rho_{\phi,\beta}= ev_\beta \circ (\phi\times id)$ and $ev_\beta: \Map(X,M)\times \Delta^k \to M^{i+1}$ is given by
\begin{multline*}
ev_\beta (\gamma:  X\to M, t_1\leq\cdots\leq t_k):=\left( \big(\gamma\circ \beta(0)\big) (t_1\leq\cdots\leq t_k) ,\right.\\
  \left.\big(\gamma\circ \beta(1)\big) (t_1\leq\cdots\leq t_k),\cdots, \big(\gamma\circ \beta(i)\big) (t_1\leq\cdots\leq t_k)\right).
\end{multline*}
The naturality of $ev_\beta$ is obvious and induces the one of $\Ch_{\beta}$. Hence we get a well defined map $\Ch^X: CH_\bullet^{S_\bullet(X)}(\Omega,\Omega) \to \Omega^\bullet(\Map(X,M)) $. It is a chain map for the same reason as above. The image of $\Ch^X$ in $\Omega^\bullet(\Map(X,M))$ is denoted by $\Chen(\Map(X,M))$.
\end{defn}
\begin{remk}Note that if $X =|Y_\bullet|$, then $\Ch^{Y_\bullet}=\Ch^{X}\circ \eta_*$ where $\eta_*: CH_\bullet^{Y_\bullet}(\Omega,\Omega)\to CH_\bullet^{S_\bullet(X)}(\Omega,\Omega)$ is the chain map induced by  $\eta: Y_\bullet \to S_\bullet(|Y_\bullet|)=S_\bullet(X)$. 
\end{remk}

\begin{remk}
Definition~\ref{Singular} easily extends to any (infinite) pointed simplicial set $Y_\bullet$, see \cite{P}.
\end{remk}

\subsection{Examples: the circle, the sphere and the torus}
We will demonstrate the above definitions in three examples provided by the circle $S^1$, the torus $\T$, and the $2$-sphere $S^2$. We then demonstrate how to wedge two squares along an edge, and how to collapse an edge to a point. We start with the circle $S^1$.
\begin{ex}[The circle $S^1$]\label{S^1}
The pointed simplicial set $S^1_\bullet$ is defined in degree $k$ by $S^1_k=\{0,\dots,k\}$, i.e. it has exactly $k+1$ elements. We define the face maps $d_i:S^1_{k}\to S^1_{k-1}$, for $0\leq i \leq k$, and degeneracies $s_i:S^1_k\to S^1_{k+1}$, for $0\leq i\leq k$, as follows, cf. \cite[6.4.2]{L}. For $i=0,\dots,k-1$, let $d_i(j)$ be equal to $j$ or $j-1$ depending on $j=0,\dots,i$ or $j=i+1,\dots,k$. Let $d_k(j)$ be equal to $j$ or $0$ depending on $j=0,\dots, k-1$ or $j=k$. For $i=0,\dots,k$, let $s_i(j)$ be equal to $j$ or $j+1$ depending on $j=0,\dots,i$ or $j=i+1,\dots, k$.

In this case, we have $(\# S^1_k-1)=k$, so that we obtain for the Hochschild chain complex over $S^1_\bullet$, $ CH^{S^1_\bullet}_\bullet(A,A)= \bigoplus_{k\geq 0} A\otimes A^{\otimes k}$. The differential is given by
\begin{multline*}
D(a_0\otimes\dots\otimes a_{k}) = \sum_{i=0}^k \pm a_0\otimes\dots\otimes d(a_i)\otimes \dots\otimes a_k\\
+\sum_{i=0}^{k-1} \pm a_0\otimes\dots\otimes(a_i\cdot a_{i+1})\otimes\dots\otimes a_k \pm (a_k\cdot a_0)\otimes a_1\otimes\dots\otimes a_{k-1},
\end{multline*} (see Definition~\ref{D:Hoch} for the signs)
which is just the usual Hochschild chain complex $CH_\bullet(A,A)$ of a differential graded algebra.

Note that $|S^1_\bullet |=S^1$, cf. \cite[7.1.2]{L}, whose only non-degenerate simplicies are $0\in S^1_0$ and $1\in S^1_1$. Now, if we view $S^1$ as the interval $I=[0,1]$ where the endpoints $0$ and $1$ are identified, then the map $\eta(i):S^1_k\to S_k(S^1)=\Map(\Delta^k,S^1)$ from definition \ref{Chen-map} is given by $\eta(i)(0\leq t_1\leq \dots\leq t_k\leq 1)=t_i$, where we have set $t_0=0$. Thus, the evaluation map \eqref{eval} becomes $ev(\gamma:S^1\to M, t_1\leq\dots\leq t_k)=(\gamma(0),\gamma(t_1),\dots,\gamma(t_k))\in M^{k+1}$. Furthermore, we can recover the classical Chen iterated integrals $\It^{S^1_\bullet}:CH_\bullet(A,A)\to \Omega^\bullet(M^{S^1})$ as follows. For a plot $\phi:U\to M^{S^1}$ we have,
\begin{multline*}
\It^{S^1_\bullet}(a_0\otimes\dots\otimes a_k)_\phi=\left(\itint a_0\dots a_k\right)_{\phi} = \int_{\Delta^k} (\rho_\phi)^*(a_0\otimes\dots\otimes a_k) \\
=(\pi_0)^*(a_0)\wedge \int_{\Delta^k} (\widetilde {\rho_\phi})^*(a_1\otimes\dots\otimes a_k)= (\pi_0)^*(a_0)\wedge \int a_1\dots a_k,
\end{multline*}
where $\widetilde {\rho_\phi}: U\times \Delta^k\stackrel {\phi\times id}\longrightarrow M^{S^1}\times \Delta^k\stackrel {\widetilde{ev}} \to M^k$ induces the classical Chen iterated integral $\int a_1\dots a_k$ from \cite{C} and $\pi_0:M^{S^1} \to M$ is the evaluation at the base point $\pi_0:\gamma \mapsto \gamma(0)$.
\end{ex}

\begin{ex}\label{E:torus}[The torus $\mathbb T$]
In this case, we can take $\T_\bullet$ to be the diagonal simplicial set associated to the bisimplicial set $S^1_\bullet\times S^1_\bullet$, i.e. $\T_k=S^1_k\times S^1_k$, see \cite[Appendix B.15]{L}. Thus, $\T_k$ has $(k+1)^2$ elements, so that we may write $\T_k=\{(p,q)  \,|\,\, p,q=0,\dots,k\}$ which we equipped with the lexicographical ordering. The face maps $d_i:\T_k\to \T_{k-1}$ and degeneracies $s_i:\T_k\to \T_{k+1}$, for $i=0,\dots,k$, are given as the products of the differentials and degeneracies of $S^1_\bullet$, i.e. $d_i(p,q)=(d_i(p), d_i(q))$ and $s_i(p,q)=(s_i(p), s_i(q))$.

With this description, we obtain $CH_\bullet^{\T_\bullet}(A,A)=\bigoplus_{k\geq 0} A\otimes A^{\otimes (k^2+2k)}$. If we index forms in $M$ by tuples $(p,q)$ as above, then we obtain homogenous elements of $CH_\bullet^{\T_\bullet}(A,A)$ as linear combinations of tensor products $a_{(0,0)}\otimes\dots \otimes a_{(k,k)}\in CH_\bullet^{\T_\bullet}(A,A)$. The differential $D(a_{(0,0)}\otimes\dots \otimes a_{(k,k)})$ on $CH_\bullet^{\T_\bullet}(A,A)$ consists of a sum
\begin{equation*}
\sum_{(p,q)=(0,0)}^{(k,k)} \pm a_{(0,0)}\otimes \dots\otimes d a_{(p,q)}\otimes \dots\otimes a_{(k,k)}+\sum_{i=0}^k \pm d_i(a_{(0,0)}\otimes \dots\otimes a_{(k,k)}).
\end{equation*}
The face maps $d_i$  can be described more explicitly, when placing $a_{(0,0)}\otimes \dots\otimes a_{(k,k)}$ in a $(k+1)\times(k+1)$ matrix. For $i=0,\dots,k-1$, we obtain $d_i(a_{(0,0)}\otimes\dots \otimes a_{(k,k)})$ by multiplying the $i$th and $(i+1)$th rows and the $i$th and $(i+1)$th columns simultaneously, \emph{i.e.}, $d_i(a_{(0,0)}\otimes\dots \otimes a_{(k,k)})$ is equal to
\begin{equation*}
\begin{matrix}
 a_{(0,0)}\otimes& \dots &\otimes a_{(0,i)}a_{(0,i+1)}\otimes& \dots &\otimes a_{(0,k)} \\  
 \vdots && \vdots &&\vdots \\
 a_{(i-1,0)}\otimes &\dots & \otimes a_{(i-1,i)}a_{(i-1,i+1)}\otimes & \dots &\otimes a_{(i-1,k)} \\ 
 a_{(i,0)} a_{(i+1,0)}\otimes& \dots &\otimes a_{(i,i)}a_{(i,i+1)}a_{(i+1,i)}a_{(i+1,i+1)}\otimes& \dots &\otimes a_{(i,k)}a_{(i+1,k)} \\  
 a_{(i+2,0)}\otimes& \dots &\otimes a_{(i+2,i)}a_{(i+2,i+1)}\otimes& \dots &\otimes a_{(i+2,k)} \\  
 \vdots && \vdots &&\vdots \\
 a_{(k,0)}\otimes& \dots &\otimes a_{(k,i)}a_{(k,i+1)}\otimes& \dots &\otimes a_{(k,k)} \\  
\end{matrix} 
\end{equation*}
The differential $d_k$ is obtained by multiplying the $k$th and $0$th rows and the $k$th and $0$th columns simultaneously, \emph{i.e.}, $d_k(a_{(0,0)}\otimes\dots \otimes a_{(k,k)})$  equals
\begin{equation*}
\pm \quad\begin{matrix}
a_{(0,0)}a_{(0,k)}a_{(k,0)} a_{(k,k)}\otimes &a_{(0,1)}a_{(k,1)}\otimes & \dots& \otimes a_{(0,k-1)}a_{(k,k-1)}\\
a_{(1,0)}a_{(1,k)}\otimes & a_{(1,1)}\otimes & \dots & \otimes a_{(1,k-1)}\\
 \vdots & \vdots & & \vdots\\
a_{(k-1,0)}a_{(k-1,k)}\otimes & a_{(k-1,1)}\otimes & \dots & \otimes a_{(k-1,k-1)}
\end{matrix}
\end{equation*}
where the sign $\pm$ is the Koszul sign (with respect to the lexicographical order) given by moving the $k$th row and lines across the matrix.

Note that $|\T_\bullet|=|S^1_\bullet|\times |S^1_\bullet|=\T$, which has non-degenerate simplicies $(0,0)\in \T_0, (0,1),(1,0),(1,1) \in \T_1$ and $(1,2),(2,1)\in \T_2$. Now, if we view the torus $\T$ as the square $[0,1]\times [0,1]$ where horizontal and vertical boundaries are identified, respectively, then the map $\eta(p,q):\T_k\to \Map(\Delta^k,\T)$ is given by $\eta(p,q)(0\leq t_1\leq\dots\leq t_k\leq 1)=(t_p,t_q)\in \T$, for $p,q=0,\dots,k$ and $t_0=0$. Thus, the evaluation map in definition \ref{Chen-map} becomes
\begin{equation*}
ev(\gamma: \mathbb T\to M, t_1\leq\cdots\leq t_k)=\quad
\begin{matrix}
 \big(\gamma(0,0),\, \gamma(0,t_1),\, \cdots,\, \gamma( 0,t_k),\\ \quad \gamma (t_1,0),\gamma(t_1,t_1),\cdots, \gamma(t_1,t_k),\\ \quad ...\\ \quad
\gamma (t_k,0),\gamma(t_k,t_1),\cdots, \gamma(t_k,t_k)\big)
\end{matrix}
\end{equation*} 
According to definition \ref{Chen-map}, the iterated integral $\Ch^\T(a_{(0,0)}\otimes \dots\otimes a_{(k,k)})$ is given by a pullback under the above map $M^\T\times \Delta^k\stackrel {ev}\longrightarrow M^{(k+1)^2}$, and integration along the fiber $\Delta^k$.

Note that a similar description works for any higher dimensional torus $\T^d=S^1\times \cdots \times S^1$ ($d$ factors) by taking $(\T^d)_k=S^1_k\times \dots\times S^1_k$. Its underlying Hochschild chain complex $CH^{\T^d_\bullet}_\bullet(A,A)=\bigoplus_{k\geq 0} A^{\otimes (k+1)^d}$ has the $i$th face map $d_i$, for $i=0,\dots,k-1$, given by simultaneously multiplying each $i$th with $(i+1)$th hyperplane in each dimension, and a similar description for $d_k$.
\end{ex}

\begin{ex}\label{E:2sphere}[The $2$-sphere $S^2$]
In this case, we define $S^2_\bullet$ to be the simplicial set with $\# S^2_k=k^2+1$ elements in simplicial degree $k$. In order to describe the faces and degeneracies, we write $S^2_k=\{(0,0)\}\cup \{(p,q) \,|\,\, p,q=1,\dots,k\}$, and set the degeneracy to be the same as for the torus in the previous example \ref{E:torus}, \emph{i.e.} $s^{S^2_\bullet}_i(p,q)=s^{\T_\bullet}_i(p,q)$ for $(p,q)\in S^2_k$ and $i=0,\dots,k$. The $i$th differential is also obtained from the previous examples by setting $d^{S^2_\bullet}_i(p,q)=(0,0)$ in the case that $d^{S^1_\bullet}_i(p)=0$ or $d^{S^1_\bullet}_i(q)=0$, or setting otherwise $d^{S^2_\bullet}_i(p,q)=d^{\T_\bullet}_i(p,q)=(d^{S^1_\bullet}_i(p), d^{S^1_\bullet}_i(q))$.

We thus obtain $CH^{S^2_\bullet}_\bullet(A,A)=\bigoplus_{k\geq 0} A\otimes A^{k^2}$ with a differential similar to the one in example \ref{E:torus}. For example, we have for $D| A\otimes A^{2^2}\to (A\otimes A^{2^2})\oplus (A\otimes A^{1^2})$,
\begin{multline*}
D\left(
\begin{matrix}
a_{(0,0)}&&\\ 
& \otimes a_{(1,1)}& \otimes a_{(1,2)} \\
& \otimes a_{(2,1)}& \otimes a_{(2,2)} 
\end{matrix}
\right)=\sum_{(p,q)} \pm \,\text{apply $d$ to each $a_{(p,q)}$}\\
+ \;\;\left(\begin{matrix}
a_{(0,0)}a_{(1,1)}a_{(1,2)}a_{(2,1)}&\\ & \otimes a_{(2,2)} 
\end{matrix}\right)\quad 
-\;\; \left(\begin{matrix}
a_{(0,0)} &\\ &\otimes a_{(1,1)}a_{(1,2)}a_{(2,1)}a_{(2,2)} 
\end{matrix}\right) \\
\pm\;\; \left(\begin{matrix}
a_{(0,0)}a_{(1,2)}a_{(2,1)}a_{(2,2)}&\\ & \otimes a_{(1,1)} 
\end{matrix}\right).\quad\quad
\end{multline*}
where the last $\pm$ sign is the Koszul sign (in the lexicographical order) given by moving $a_{(1,2)},a_{(2,1)},a_{(2,2)}$ across $a_{(1,1)}$.

It can be seen that $|S^2_\bullet|=S^2$. If we view the $2$-sphere as a square $[0,1]\times [0,1]$ where the boundary is identified to a point, then we obtain the evaluation map,
\begin{equation*}
ev(\gamma:S^2\to M, t_1\leq\cdots\leq t_k)=
\left(\begin{matrix}
\gamma(0,0),& \\ &\gamma(t_1,t_1),\cdots, \gamma(t_1,t_k),\\  & \vdots\quad \quad \quad \quad \vdots \\ &\gamma(t_k,t_1),\cdots, \gamma(t_k,t_k)
\end{matrix}\right) \quad\quad\in M^{k^2+1}.
\end{equation*}
\end{ex}
This completes our three examples $S^1, \T$, and $S^2$. Later, in section \ref{S:surfacemodel}, we will describe a simplicial model for a surface $\Sigma^g$ of genus $g$, which is build out of collapsing an edge to a point and wedging squares along vertices or edges. The essential ideas in these constructions will be demonstrated in the next example.
\begin{ex}\label{E:wedge}[Wedge along an edge or a vertex]
The simplicial model for the point $pt_\bullet$  is given by $pt_k=\{0\}$ for all $k\geq 0$, with trivial faces and degeneracies. Next, we can give a simplicial model for the interval $I_\bullet$ by taking $I_k=\{0,\dots,k+1\}$ with differential for $i=0,\dots,k$, $d_i(j)$ equal to $j$ or $j-1$ depending on $j\leq i$ or $j>i$. The associated Hochschild chain complex is just the two sided bar complex, $CH_\bullet^{I_\bullet}(A,A)=\bigoplus_{k\geq 0} A\otimes A^{k}\otimes A= B(A,A,A)$. Similarly, we have the simplicial square $I^2_\bullet:=I_\bullet\times I_\bullet$, i.e. $I^2_k=I_k\times I_k$ has $(k+2)^2$ elements with differential $d^{I^2_\bullet}_i=d^{I_\bullet}_i\times d^{I_\bullet}_i$.

We can use the pushout construction from Lemma \ref{L:pushout} to glue two squares along an edge. In fact, this can be easily done using the inclusion $inc:I_\bullet\to I^2_\bullet$ twice to obtain $I^2_\bullet\cup_{I_\bullet} I^2_\bullet$. Similarly, we can wedge two squares at a vertex, where we use the inclusion $inc':pt_\bullet \to I^2_\bullet$ to obtain $I^2_\bullet\cup_{pt_\bullet} I^2_{\bullet}$. Note that in this case, we do not need to assume that the inclusion $inc':pt_\bullet \to I^2_\bullet$ preserves the basepoints.

A more interesting operation may be obtained via the collapse map $col:I_\bullet\to pt_\bullet$, which together with the inclusion $inc:I_\bullet\to I^2_\bullet$ induces the square with one collapsed edge $I^2_\bullet\cup_{I_\bullet} pt_\bullet$. We will use this type of  construction in subsection \ref{S:surfacemodel} to obtain our model for the surface $\Sigma^g$ of genus $g$.

There are two maps $s,t:pt_\bullet \to I_\bullet$ given by $s(0)=0$ and $t(0)=k+1$ in simplicial degree $k$ and a (unique) projection $p:  I_\bullet \to pt_\bullet$. These 3 maps induces 4 inclusions 
$s_{ij}: pt_\bullet \to I^2_\bullet=I_\bullet\times I_\bullet$  mapping the point to one of the corner of the square, 2 projections $p_j:I^2_\bullet\to I_\bullet$ ($i,j\in \{1,2\}$) and a collapse map $p\times p:I^2_\bullet\to pt_\bullet\times pt_\bullet\cong pt_\bullet$. The following Lemma is trivial but useful.
\begin{lem}\label{L:collapse}
 The maps $s$, $t$, $p$, $s_{ij}$ and $p_j$ are maps of simplicial sets.
\end{lem}
 
\end{ex}

\subsection{The cup product for the mapping space $\Map(Y,M)$}\label{section-product}
We now show that there is a differential graded algebra structure on $CH_\bullet^{Y_\bullet}(\Omega,\Omega)$ and $\Chen(M^Y)$, and that the iterated integral $\It^{Y_\bullet}$ preserves this algebra structure. The algebra structure on $\Chen(M^Y)$ is the one induced by the wedge product in $\Omega^\bullet(M^Y)$, cf. definition \ref{Chen-map}. To define the product on $CH_\bullet^{Y_\bullet}(A,A)$, we first recall the shuffle product for simplicial vector spaces $V_\bullet$ and $W_\bullet$, see e.g. \cite[Lemma 1.6.11]{L}. 

\begin{defn}\label{cross}
For two simplicial vector spaces $V_\bullet$ and $W_\bullet$, one defines a simplicial structure on the simplicial space $(V\times W)_k:=V_k\otimes W_k$ using the boundaries $d^V_i\otimes d^W_i$ and degeneracies $s^V_i\otimes s^W_i$. There is a shuffle product $sh:V_p\otimes W_q\to (V\times W)_{p+q}$, $$ sh(v\otimes w)=\sum_{(\mu,\nu)} sgn(\mu,\nu) (s_{\nu_q}\dots s_{\nu_1}(v)\otimes s_{\mu_p}\dots s_{\mu_1}(w)), $$
where $(\mu,\nu)$ denotes a $(p,q)$-shuffle, i.e. a permutation of $\{0,\dots,p+q-1\}$ mapping $0\leq j\leq p-1$ to $\mu_{j+1}$ and $p\leq j\leq p+q-1$ to $\nu_{j-p+1}$, such that $\mu_1<\dots<\mu_p$ and $\nu_1<\dots<\nu_q$. In particular, for $V_\bullet=W_\bullet$, this becomes $sh:V_p\otimes V_q\to V_{p+q}\otimes V_{p+q}$.

Since $CH_\bullet^{Y_\bullet}(A,A)$ is a simplicial vector space, we obtain an induced shuffle map $sh$ on $CH_\bullet^{Y_\bullet}(A,A)$. Composing this with $CH^{Y_\bullet}_\bullet(\mu,\mu)$, where $\mu:A\otimes A\to A$ denotes the product of $A$, which is an algebra map since $A$ is graded commutative, we obtain the desired shuffle product $sh_{Y_\bullet}$ of $CH_\bullet^{Y_\bullet}(A,A)$,
$$ sh_{Y_\bullet}:CH_\bullet^{Y_\bullet}(A,A)\otimes CH_\bullet^{Y_\bullet}(A,A)\stackrel {sh} \to CH_\bullet^{Y_\bullet}(A\otimes A,A\otimes A)\stackrel {CH^{Y_\bullet}_\bullet(\mu,\mu)} \longrightarrow  CH_\bullet^{Y_\bullet}(A,A).  $$
\end{defn}
\begin{prop}\label{P:shuffleinvariance}  The shuffle product $sh_{Y_\bullet}:CH_\bullet^{Y_\bullet}(A,A)^{\otimes 2} \to CH_\bullet^{Y_\bullet}(A,A)$ makes $CH_\bullet^{Y_\bullet}(A,A)$ a differential graded commutative algebra, which is natural in $A$. If $f_\bullet: X_\bullet\to Y_\bullet$ is a map of simplicial sets, then the induced map $f_*:CH_\bullet^{X_\bullet}(A,A)\to CH_\bullet^{Y_\bullet}(A,A)$ is a map of algebras and a quasi-isomorphism if the map $H_\bullet(f_\bullet):H_\bullet(X_\bullet) \to H_\bullet(Y_\bullet)$  induced by $f_\bullet$ is an isomorphism. 
\end{prop}
\begin{proof}
  The shuffle product $V_\bullet \otimes V_\bullet \to (V\times V)_\bullet$ for simplicial vector spaces is associative and graded commutative (see \cite[Section 1.6]{L}). Further  $\mu: A\otimes A\to A$ is map of differential graded algebras since $A$ is (differential graded) commutative. Hence $sh_{Y_\bullet}$ is an associative and graded commutative multiplication and a map of chain complexes. 

 That $f_*:CH_\bullet^{X_\bullet}(A,A)\to CH_\bullet^{Y_\bullet}(A,A)$ is a map of algebras follows by naturality of the shuffle product and the last claim is proved in~\cite[Proposition 2.4]{P}. 
\end{proof}

  In view of the above Proposition~\ref{P:shuffleinvariance}, Lemma~\ref{L:pushout} has the following counterpart.
\begin{cor} \label{C:pushoutalgebra} Let $X_\bullet$, $Y_\bullet$, $W_\bullet$ and $Z_\bullet$, $f_\bullet: Z_\bullet\to X_\bullet$, $g_\bullet:Z_\bullet\to Y_\bullet$ be as in Lemma~\ref{L:pushout} and Definition~\ref{D:wedge}. There is a natural morphism of differential graded algebras $$CH_\bullet^{X_\bullet}(A,A)\otimes_{CH_\bullet^{Z_\bullet}(A,A)} CH_\bullet^{Y_\bullet}(A,A) \to CH_\bullet^{W_\bullet}(A,A)$$  If $Z_\bullet$ injects into either $Z_\bullet\stackrel {f_\bullet}\to X_\bullet$ or $Z_\bullet\stackrel{g_\bullet}\to Y_\bullet$, this natural map  is a quasi-isomorphism.
\end{cor}
\begin{proof} Proposition~\ref{P:shuffleinvariance} and the commutative diagram~\eqref{eq:pushout}  implies that $ CH_\bullet^{X_\bullet}(A,A)$ and $CH_\bullet^{Y_\bullet}(A,A)$ are $CH^{Z_\bullet}_\bullet(A,A)$-algebras and that the maps $i_*:CH_\bullet^{X_\bullet}(A,A)\to CH_\bullet^{W_\bullet}(A,A)$ $j_*:CH_\bullet^{Y_\bullet}(A,A)\to CH_\bullet^{W_\bullet}(A,A)$ are maps of (differential graded) commutative $CH^{Z_\bullet}_\bullet(A,A)$-algebras. Since $i_*\circ f_*=j_*\circ g_*$, the composition
\begin{equation*}
CH_\bullet^{X_\bullet}(A,A)\otimes CH_\bullet^{Y_\bullet}(A,A) \stackrel{i_*\otimes j_*}\longrightarrow CH_\bullet^{W_\bullet}(A,A) \otimes CH_\bullet^{W_\bullet}(A,A) \stackrel{sh_{W_\bullet}} \longrightarrow CH_\bullet^{W_\bullet}(A,A)
\end{equation*}
induces a natural morphism $CH_\bullet^{X_\bullet}(A,A)\otimes_{CH_\bullet^{Z_\bullet}(A,A)} CH_\bullet^{Y_\bullet}(A,A) \to CH_\bullet^{W_\bullet}(A,A)$ of $CH_\bullet^{Z_\bullet}(A,A)$-algebras. The last statement follows from Lemma~\ref{L:pushout} and the fact that the shuffle product $V_\bullet\otimes W_\bullet \to (V\times W)_\bullet$ is a natural quasi-isomorphism of simplicial vector spaces.
\end{proof}

By Proposition~\ref{P:shuffleinvariance}, for any simplicial set $X_\bullet$ and commutative (differential graded) algebra $A$, $(CH^{X_\bullet}_\bullet(A,A), D, sh_{X_\bullet})$ is again a commutative differential graded algebra. Thus its Hochschild chain complex $CH^{Y_\bullet}_\bullet\big(CH^{X_\bullet}_\bullet(A,A),CH^{X_\bullet}_\bullet(A,A)\big)$ 
is defined for any simplicial set $Y_\bullet$ and is a commutative differential graded algebra.
\begin{cor}\label{C:product}
There is a natural (with respect to $X_\bullet$, $Y_\bullet$ and $A$) quasi-isomorphism of algebras
$$CH^{X_\bullet}_\bullet\big(CH^{Y_\bullet}_\bullet(A,A),CH^{Y_\bullet}_\bullet(A,A)\big) \stackrel{\sim}\to CH^{(X\times Y)_\bullet}_\bullet(A,A)$$ where $(X\times Y)_\bullet$ is the diagonal of the bisimplicial set $X_\bullet \times Y_\bullet$, that is $(X\times Y)_n=X_n\times Y_n$.
\end{cor}
\begin{proof}
The quasi-isomorphism is induced by the Eilenberg-Zilber quasi-isomorphism (\emph{i.e.} the shuffle product). In order to define it more explicitly, we need the following definition of Hochschild chain complex $CH^{X_\bullet}(R_\bullet)$ over a simplicial set $X_\bullet$ for a \emph{simplicial algebra} $R_\bullet$. This is a bisimplicial vector space which is given in simplicial bidegree $(p,q)$ by $CH^{X_p}(R_q)=R_q\otimes R_q^{\otimes x_p}$ where $x_p+1=\# X_p$. Its face maps $CH^{X_p}(R_q)\stackrel{d_i}\to CH^{X_{p-1}}(R_q)$ are given by the usual ones for the Hochschild complex over $X_\bullet$ as in Definition~\ref{D:Hoch} and similarly for its degeneracies along the $X_\bullet$ direction. Thus, $CH^{X_\bullet}(R_q)=CH^{X_\bullet}_\bullet(R_q,R_q)$  is the Hochschild chain complex of the algebra $R_q$ over $X_\bullet$. The simplicial face maps $CH^{X_p}(R_q)\stackrel{d_i}\to CH^{X_p}(R_{q-1})$ along the $R_\bullet$ direction are induced by the face maps of $R_\bullet$: 
$$CH^{X_p}(R_q)= R_q\otimes R_q^{\otimes x_p} \stackrel{d_i^{\otimes 1+x_p}}\longrightarrow R_{q-1}\otimes R_{q-1}^{\otimes x_{p}}=CH^{X_p}(R_{q-1})$$ and similarly for the degeneracies. Thus,
$CH^{X_p}(R_\bullet)=(R \times \cdots\times R)_{\bullet}$ is the $(1+x_p)$-times iterated cross-product of the simplicial algebra $R_\bullet$ (see Definition~\ref{cross} above).
  If $R_\bullet$ is a simplicial commutative differential graded algebra, then $CH^{X_\bullet}(R_\bullet)$ is a bisimplicial commutative differential graded algebra. Note that for the standard Hochschild chain complex over $S^1_\bullet$, this definition  was first introduced by Goodwillie~\cite{Go}.

By the (generalized) Eilenberg-Zilber theorem~\cite{MacL, GJ}, there is a natural quasi-isomorphism
$EZ:CH^{X_\bullet}(R_\bullet) \to diag(CH^{X_\bullet}(R_\bullet))_\bullet$, where $diag (CH^{X_\bullet}(R_\bullet))_\bullet$ is the diagonal simplicial set associated to the bisimplicial set $CH^{X_\bullet}(R_\bullet)$. Hence, in simplicial degree $n$, $diag(CH^{X_\bullet}(R_\bullet))_n=CH^{X_n}(R_n)=R_n\otimes  R_n^{\otimes x_n}$. Note that the map $EZ:CH^{X_p}(R_q)=R_q^{\otimes 1+x_p}\to CH^{X_{p+q}}(R_{p+q})=R_{p+q}^{\otimes 1+x_{p+q}}$ is given by a formula similar to the one in Definition~\ref{cross}, that is 
\begin{equation}\label{E:EZproduct} EZ=\sum_{(\mu,\nu)} sgn(\mu,\nu) \big(s^{X_\bullet}_{\nu_q}\dots s^{X_\bullet}_{\nu_1}\big) \circ \left( \big( s_{\mu_p}^{R_\bullet}\dots s_{\mu_1}^{R_\bullet}\big)^{\otimes 1+x_p}\right) \end{equation} where $(\mu,\nu)$ denotes a $(p,q)$-shuffle and $s_j^{X_\bullet}$, $s_k^{R_\bullet}$ are the degeneracies along the $X_\bullet$ and $R_\bullet$ simplicial directions respectively.

From Definition~\ref{D:Hoch}, it is clear, that $CH^{Y_\bullet}_\bullet(A,A)$ is a simplicial differential graded commutative algebra, and, that $CH^{(X\times Y)_\bullet}_\bullet(A,A)\cong diag(CH^{X_\bullet}(CH^{Y_\bullet}_\bullet(A,A)))_\bullet$, since $CH^{(X\times Y)_n}_\bullet(A,A)=A^{\otimes (1+x_n)(1+y_n)}\cong \left(A^{\otimes (1+y_n)} \right)^{\otimes 1+x_n}$. Thus, the Eilenberg-Zilber map~\eqref{E:EZproduct} gives a quasi-isomorphism of underlying chain complexes
$$ CH^{X_\bullet}(CH^{Y_\bullet}_\bullet(A,A)) \stackrel{EZ}\longrightarrow CH^{(X\times Y)_\bullet}_\bullet(A,A).$$
Note that on the left hand side, $CH^{Y_\bullet}_\bullet(A,A)$ is considered as a simplicial differential graded algebra. Now we need to define a quasi-isomorphism $$ CH^{X_\bullet}_\bullet(CH^{Y_\bullet}_\bullet(A,A),CH^{Y_\bullet}_\bullet(A,A))  \to CH^{X_\bullet}(CH^{Y_\bullet}_\bullet(A,A)) $$ where on the left, $CH^{Y_\bullet}_\bullet(A,A)$ is equipped with its structure of commutative differential graded algebra given by the shuffle product $sh_{Y_\bullet}$ and the Hochschild differential $D$. 
Iterating the Eilenberg-Zilber map $x_p$-times, we get a quasi-isomorphism
\begin{multline*} 
EZ^{(x_p)}:CH^{X_p}_\bullet(CH^{Y_\bullet}_\bullet(A,A),CH^{Y_\bullet}_\bullet(A,A))= \left(CH^{Y_\bullet}_\bullet(A,A)\right)^{\otimes (1+x_p)}\\
 \stackrel{sh^{\otimes x_p}}\longrightarrow \left(CH^{Y_\bullet}_\bullet(A,A)\right)^{\times (1+x_p)}
 \cong CH^{X_p}(CH^{Y_\bullet}_\bullet(A,A)),
\end{multline*} 
where $sh$ is the shuffle product as in Definition~\ref{cross}. Thus the composition
$$ CH^{X_\bullet}_\bullet(CH^{Y_\bullet}_\bullet(A,A),CH^{Y_\bullet}_\bullet(A,A))\stackrel{EZ\circ \big( EZ^{(x_\bullet)}\big)}\longrightarrow  CH^{(X\times Y)_\bullet}_\bullet(A,A)$$ is a natural quasi-isomorphism. 

Since the algebra structure on $CH^{X_\bullet}_\bullet(CH^{Y_\bullet}_\bullet(A,A),CH^{Y_\bullet}_\bullet(A,A))$ is the composition of the shuffle product 
\begin{multline*}
sh:CH^{X_{p_1}}_\bullet(CH^{Y_\bullet}_\bullet(A,A),CH^{Y_\bullet}_\bullet(A,A))\otimes  CH^{X_{p_2}}_\bullet(CH^{Y_\bullet}_\bullet(A,A),CH^{Y_\bullet}_\bullet(A,A))\\ \to CH^{X_{p_1+p_2}}_\bullet\big(CH^{Y_\bullet}_\bullet(A,A)^{\otimes 2}, CH^{Y_\bullet}_\bullet(A,A)^{\otimes 2} \big)
\end{multline*}
with the map $CH_\bullet^{X_\bullet}(sh_{Y_\bullet})$ (also induced by the shuffle product see Definition~\ref{cross}), the fact that the natural map $CH^{X_\bullet}_\bullet(CH^{Y_\bullet}_\bullet(A,A),CH^{Y_\bullet}_\bullet(A,A))\to  CH^{(X\times Y)_\bullet}_\bullet(A,A)$ is a map of algebras follows from the associativity and commutativity of the shuffle product.
\end{proof}

\begin{ex}\label{E:torusproduct}
If $\T_\bullet$ is the simplicial model for the torus given in Example~\ref{E:torus}, then, by Corollary~\ref{C:product} above, $CH^{T_\bullet}_\bullet(A,A)$ is quasi-isomorphic, as an algebra, to $CH^{S^1_\bullet}_\bullet\big(CH^{S^1_\bullet}_\bullet(A,A),CH^{S^1_\bullet}_\bullet(A,A))$, that is to the standard Hochschild complex of the standard Hochschild complex of  the algebra $A$. 
\end{ex}

Using a decomposition of the product $\Delta^k\times \Delta^l$ into a union of $(k+l)$-simplices $\Delta^{k+l}$, which is indexed by the set of all shuffles, we obtain the following proposition.
\begin{prop}\label{P:it-algebra}
For any compact, oriented manifold $M$, the iterated integral map $\Ch^{Y_\bullet}:(CH^{Y_\bullet}_\bullet(\Omega,\Omega),sh_{Y_\bullet})\to (\Omega^\bullet(M^Y),\wedge)$ is a map of algebras.
\end{prop}
\begin{proof}
We need to show, that for $a_0\otimes\dots\otimes a_{y_k}, b_0\otimes \dots\otimes b_{y_l}\in CH_\bullet^{Y_\bullet}(\Omega,\Omega)$, 
\begin{equation}\label{Eq:It-product}
\It^{Y_\bullet}(a_0\otimes\dots\otimes a_{y_k})\wedge \It^{Y_\bullet}(b_0\otimes\dots\otimes b_{y_l})=\It^{Y_\bullet}(sh_{Y_\bullet}(a_0\otimes\dots\otimes a_{y_k}, b_0\otimes\dots\otimes b_{y_l})).
\end{equation}
Let $\phi:U\to M^Y$ be a plot, and $\rho_\phi^{k+l}:U\times \Delta^{k+l}\to M^{y_{k+l}+1}$ the map from \eqref{rho_phi}. Note that each degeneracy $s_i:Y_{r}\to Y_{r+1}$ induces a map $M^{s_i}:M^{y_{r+1}+1}\to M^{y_{r}+1}$, which in turn induces the degeneracy $s_i:\Omega^{\otimes y_{r}+1}\to \Omega^{\otimes y_{r+1}+1}$ on $CH_\bullet^{Y_\bullet}(\Omega,\Omega)$. Since the multiplication $\mu:\Omega^{\otimes 2}\to \Omega$ is obtained by pullback along the diagonal $\delta:M\to M\times M$,  the term on the right side of \eqref{Eq:It-product} becomes, 
\begin{multline*}
 \sum_{(\mu,\nu)} \pm \int_{\Delta^{k+l}} (\rho_\phi^{k+l})^*\circ (\delta^{y_{k+l}+1})^*\Big(s_{\nu_l}\dots s_{\nu_1}(a_0\otimes\dots\otimes a_{y_k})\otimes s_{\mu_k}\dots s_{\mu_1}(b_0\otimes\dots\otimes b_{y_l}) \Big) \\
 = \sum_{(\mu,\nu)} \pm \int_{\Delta^{k+l}} (\rho_\phi^{k+l})^*\circ \alpha^*\Big(a_0\otimes\dots\otimes a_{y_k}\otimes b_0\otimes\dots \otimes b_{y_l}\Big),
 \end{multline*}
 where $\alpha:M^{y_{k+l}+1}\to M^{y_{k+l}+1}\times M^{y_{k+l}+1}\to M^{y_k+1}\times M^{y_l+1}$ is the composition of the diagonal of $M^{y_{k+l}+1}$ with the map $(M^{s_{\nu_l}}\dots M^{s_{\nu_1}})\times (M^{s_{\mu_k}}\dots M^{s_{\mu_1}})$.
 
Recall the degeneracies $\sigma_i:\Delta^{r+1}\to \Delta^{r}, (0\leq t_1\leq \dots\leq t_{r+1}\leq 1)\mapsto (0\leq t_1\leq \dots\leq \hat t_i\leq \dots\leq t_{r+1} \leq 1)$, which removes the $i$th coordinate from the standard simplex, cf. \cite[Appendix B.6]{L}. Then for any $\gamma\in M^Y$, the map $ev$ from \eqref{eval} makes the following diagram commutative,
$$
\xymatrix{\Delta^{r+1} \ar[rr]^{ev(\gamma,-)} \ar[d]_{\sigma_i} && M^{y_{r+1}+1} \ar[d]^{M^{s_i}} \\
 \Delta^r \ar[rr]^{ev(\gamma,-)} && M^{y_r+1} }
$$
Thus, for any shuffle $(\mu,\nu)$, we obtain the commutative diagram,
$$
\xymatrix{
U\times \Delta^{k+l}\ar[dd]_{id\times \beta^{(\mu,\nu)}}\ar[rr]^{\phi\times id} && M^Y\times \Delta^{k+l}\ar[r]^{ev}\ar[d]_{=} & M^{y_{k+l}+1}\ar[d] \\
   && M^Y\times \Delta^{k+l} \ar[r]^{(ev,ev)\quad }\ar[d]_{id\times \beta^{(\mu,\nu)}} & M^{y_{k+l}+1} \times M^{y_{k+l}+1} \ar[d] \\
U\times \Delta^k\times \Delta^l \ar[rr]^{\phi\times id\times id} && M^Y\times \Delta^k\times \Delta^l \ar[r]^{(ev,ev)} & M^{y_{k}+1} \times M^{y_{l}+1} 
}
$$
where the right vertical map is $\alpha$, and $\beta^{(\mu,\nu)}=(\sigma_{\nu_l}\dots\sigma_{\nu_1},\sigma_{\mu_k}\dots\sigma_{\mu_1}):\Delta^{k+l}\to \Delta^k\times \Delta^l$. Thus, $\alpha\circ \rho_\phi^{k+l}=\rho_\Phi \circ(id\times \beta^{(\mu,\nu)})$, where $\rho_\Phi$ denotes the bottom map. Using the decomposition of $\Delta^{k}\times\Delta^{l}=\amalg_{(\mu,\nu)} \beta^{(\mu,\nu)}(\Delta^{k+l})$, we can simplify the right hand side of \eqref{Eq:It-product} to
\begin{multline*}
\int_{\amalg_{(\mu,\nu)} \beta^{(\mu,\nu)}(\Delta^{k+l})} (\rho_\Phi)^*(a_0\otimes\dots\otimes a_{y_k}\otimes b_0\otimes\dots \otimes b_{y_l})  \\
=\int_{\Delta^k} (\rho_\phi)^*(a_0\otimes\dots\otimes a_{y_k})\wedge \int_{\Delta^l} (\rho_\phi)^*(b_0\otimes\dots\otimes b_{y_l})\\
=\It^{Y_\bullet}(a_0\otimes\dots\otimes a_{y_k})\wedge \It^{Y_\bullet}(b_0\otimes\dots\otimes b_{y_l}),
\end{multline*}
which is the claim.
\end{proof}
The previous proposition shows, that the wedge product of two iterated integrals is again an iterated integral, so that ``$\wedge$'' preserves $\Chen(M^Y)\subset \Omega^\bullet(M^Y)$. We thus have the following Corollary.
\begin{cor}\label{P:it-algebra2}
$\Ch^{Y_\bullet}:(CH^{Y_\bullet}_\bullet(\Omega,\Omega),sh_{Y_\bullet})\to (\Chen(M^Y),\wedge)$ is a map of algebras.
\end{cor}
\begin{remk}
The proof of Proposition \ref{P:it-algebra} is essentially the same as the proof given by Patras and Thomas in \cite[Proposition 2]{PT}, and could have been deduced from \cite{PT}. We will use the relationship with \cite{PT} in the next subsection.
\end{remk}

\subsection{Chen iterated integrals as a quasi-isomorphism}\label{section-qi}

In this subsection, we show that the iterated integral map $\Ch^{Y_\bullet}: CH_\bullet^{Y_\bullet}(\Omega,\Omega)\longrightarrow \Omega^\bullet(M^Y)$ is a quasi-isomorphism under suitable connectivity conditions on $M$, where we set as usual $\Omega=\Omega^\bullet(M)$. For the proof we will apply a related result by Patras and Thomas \cite{PT}, which uses a simplicial description of cochains of $M^Y$. We start with a slight generalization of the simplicial cochain model used in \cite{PT}.
\begin{defn}
Let $Y_\bullet$ be a simplicial space, and $M$ a compact manifold. Denote by $\Co^\bullet$ a cochain functor, such as simplicial cochains, singular cochains, or De Rham forms. We define the simplicial chain complex $\Co^\bullet_\bullet(M^Y)$ by letting $\Co^\bullet_k(M^Y)=(\Co^\bullet(M^{Y_k})=\bigoplus_{p\geq 0} \Co^p(M^{Y_k}),\partial_{k})$, where $\partial_k$ is the differential on $M^{Y_k}$ induced by $\Co^\bullet$. The face maps $d_i$ and degeneracies $s_i$ of $Y_\bullet$ induce face maps $D_i:=\Co^\bullet(M^{d_i})$ and degeneracies $S_i:=\Co^\bullet(M^{s_i})$ on $\Co^\bullet_\bullet(M^Y)$.

The total complex $\Co(M^Y)^\bullet$ is defined by $\Co(M^Y)^p=\bigoplus_{k\geq 0} \Co^{p+k}(M^{Y_k})$, and has the differential $D:\Co(M^Y)^p\to \Co(M^Y)^{p+1}$, which on $\Co^{p+k}(M^{Y_k})$ is a sum of the differentials $(-1)^k\partial_k:\Co^{p+k}(M^{Y_k})\to \Co^{p+k+1}(M^{Y_{k}})$ and $\sum_{i=0}^k (-1)^i D_i:\Co^{p+k}(M^{Y_k})\to \Co^{p+k}(M^{Y_{k-1}})$.

The normalized complex $\No\Co(M^Y)^\bullet$ is defined as the quotient of $\Co(M^Y)^\bullet$ by the subspace generated by the images of the degeneracies $S_i$. It is well-known, that the projection $ \Co(M^Y)^\bullet\to \No\Co(M^Y)^\bullet$ is a quasi-isomorphism of chain complexes, see e.g. \cite{MacL}.
\end{defn}
\begin{lem}\label{L:CD}
Assume that $Y=|Y_\bullet|$ is $n$-dimensional, i.e. the highest degree of any non-degenerate simplex is $n$, and assume that $M$ is $n$-connected. Then any two cochain functors $\Co^\bullet$ and $\Do^\bullet$ induce quasi-isomorpic complexes $\Co(M^Y)^\bullet$ and $\Do(M^Y)^\bullet$.
\end{lem}
\begin{proof}
It is enough to show that $\No\Co(M^Y)^\bullet$ and $\No\Do(M^Y)^\bullet$ are quasi-isomorphic. For $\No\Co(M^Y)^\bullet$, we define the filtration by simplicial degree $F^k_\Co=\bigoplus_{0\leq l \leq k} \No\Co^\bullet(M^{Y_l})$. The $E^1$ term for this filtration is computed as the reduced homology $\bigoplus_{k\geq 0} \overline H^*(M^{Y_k})$. Using the assumptions on the connectivity of $M$, it is easy to see that the $E^1$ page is first quadrant, and thus the filtration converges to the homology $H(\No\Co(M^Y)^\bullet)$. Similar arguments give a spectral sequence converging to the homology $H(\No\Do(M^Y)^\bullet)$. Now any natural equivalence $\Fu:\Co^\bullet\to \Do^\bullet$ induces a map of spectral sequences, which is an isomorphism on the $E^1$ level. Since any two cochain models $\Co^\bullet$ and $\Do^\bullet$ can be connected by a sequence of natural equivalences, the claim follows.
\end{proof}
\begin{prop}\label{qi}
Under the assumption from Lemma \ref{L:CD}, the iterated integral map $\It^{Y_\bullet}:CH_\bullet^{Y_\bullet}(\Omega,\Omega)\to \Omega^\bullet(M^Y)\cong C^\bullet(M^Y)$ is a quasi-isomorphism.
\end{prop}
\begin{proof}
First, notice, that $\It^{Y_\bullet}:CH_\bullet^{Y_\bullet}(\Omega,\Omega)\to \Omega^\bullet(M^Y)$  factors through $\Omega(M^Y)^\bullet=\bigoplus_{k\geq 0}\Omega^\bullet(M^{Y_k})$ via
$$ \Omega\otimes \Omega^{\otimes y_k}\stackrel {Z} \to \Omega^\bullet (M^{Y_k})\stackrel {ev^*} \to \Omega^\bullet(M^Y\times \Delta^k)\stackrel {\int_{\Delta^k}} \to \Omega^\bullet(M^Y), $$
where $Z$ is the natural quasi-isomorphism obtained as the wedge of the pullbacks  of the $y_{k}+1=\#Y_k$  projections $M^{Y_k}\to M$, $ev:M^Y\times \Delta^k\to M^{Y_k}$ is the map in \eqref{eval}, and $\int_{\Delta^k}$ denotes integration along the fiber. An argument similar to Lemma \ref{L:CD} shows, that $Z$ induces a quasi-isomorphism $CH^{Y_\bullet}_\bullet(\Omega,\Omega)\to \Omega(M^Y)^\bullet$.

Denote by $\So^\bullet$ the singular cochain functor, and denote by $/[\Delta^k]:\So^\bullet(M^Y\times \Delta^k)\to \So^\bullet(M^Y)$ the slant product with the fundamental cycle of $\Delta^k$. Consider the diagram,
$$ \xymatrix{
\Omega (M^{Y})^\bullet \ar[r]^{ev^*\quad}\ar[d]& \Omega^\bullet(M^Y\times \Delta^k) \ar[r]^{\quad\int_{\Delta^k}}\ar[d]& \Omega^\bullet(M^Y)\ar[d]\\
\So(M^Y)^\bullet \ar[r]^{ev^*\quad} & \So^\bullet(M^Y\times \Delta^k) \ar[r]^{\quad/[\Delta^k]}  & \So^\bullet(M^Y) } $$
which commutes after taking homology. Since $\Omega^\bullet$ and $\So^\bullet$ are naturally equivalent, and using Lemma \ref{L:CD}, we see furthermore that the vertical maps are isomorphisms on homology. Recalling from \cite[Corollary 2]{PT}, that the bottom line is a quasi-isomorphism, we conclude that the iterated integral also induces a quasi-isomorphism, which is the claim of the proposition.
\end{proof}
\begin{remk}
An alternative proof of the above proposition may be obtained by following the ideas of Getzler Jones Petrack \cite[Theorem 3.1]{GJP}, via induction on the simplicial skeletal degree (cf. \cite[p. 8]{GJ}).
\end{remk}

Let $A^*=Hom(A,k)$ be the (graded) dual of $A$.
If we denote the graded dual of $CH_\bullet^{Y_\bullet}(A,A)$ by $CH^\bullet_{Y_\bullet}(A,A^*):=\prod_{k\geq 0}(A^*\otimes (A^*)^{\otimes y_k})$, then we also have the following dual statement to Proposition \ref{qi}.
\begin{cor}\label{qi*}
Under the assumptions from Lemma \ref{L:CD}, we have a quasi-isomorphism $(\Ch^{Y_\bullet})^*:C_\bullet(\Map(Y,M))\to CH^\bullet_{Y_\bullet}(\Omega,\Omega^*)$.
\end{cor}

\begin{remk}
In the above discussion, we did not include the Chen space $\Chen(M^Y)$, which, by definition, is given by the image of the iterated integral map $\Chen(M^Y)=Im(\Ch^{Y_\bullet}: CH_\bullet^{Y_\bullet}(\Omega,\Omega)\to \Omega^\bullet(M^Y))\subset \Omega^\bullet(M^Y)$. Chen showed in the case of the circle $Y_\bullet=S^1_\bullet$ ({\it cf.} \cite{C2}), that $\Chen(M^{S^1})$ is in fact quasi-isomorphic to $\Omega^\bullet(M^Y)$ by showing that its kernel $Ker(\It^{S^1_\bullet})$ is acyclic. In case of a general simplicial set $Y_\bullet$, this task turns out to become quite more elaborate, as the kernel $Ker(\It^{Y_\bullet})$ contains many non-trivial combinatorial restrictions, coming from the combinatorics induced by $Y_\bullet$.

Let us illustrate this by the example of a simplicial graph $G_\bullet$, having $v$ vertices and $e$ edges as its only non-degenerate simplicies. Combining the models for the interval and the point as in Example~\ref{E:wedge}, we may assume that $\#G_k=v+e\cdot k$. Then for given functions $f_1,\dots, f_v, g_1, \dots, g_e:M\to \mathbb R$, which we associate to the vertices and edges of $G_\bullet$, we can define a degree $0$ element $x\in CH^{G_1}_0(\Omega,\Omega)\cong \Omega^{\otimes v}\otimes \Omega^{\otimes e}[1]$, by setting 
$$ x=f_1\otimes \dots\otimes f_v\otimes d_{DR}(g_1\otimes \dots \otimes g_e)\in (\Omega^{\otimes v})^0\otimes (\Omega^{\otimes e})^1 [1]$$ where $d_{DR}$ is the De Rham differential.
A computation then shows that $x\in Ker(\Ch^{G_\bullet})$ exactly when for every vertex $w$ of $G_\bullet$, the product of the functions on the incoming edges $g_{i^1_w}, \dots,g_{i^{r_w}_w}$ at $w$ is equal to the product of the functions on the outgoing edges $g_{j^1_w}, \dots,g_{j^{s_w}_w}$ at $w$ up to a constant $c_w$, and these constants multiply to $1$,
$$ \forall w: \prod_{k} g_{i^k_w}=c_w\cdot \prod_{k} g_{j^k_w}, \quad\text{ and } \quad \prod_w c_w=1.$$

The conclusion is, that the explicit identification of the kernel $Ker(\It^{Y_\bullet})$ for a general simplicial set $Y_\bullet$, as well as the proof of its acyclicity, require considerably more effort. However, we conjecture that the kernel is acyclic.
\end{remk}

\section{String topology product for surfaces mapping spaces}\label{S:surfaceproduct}

Beside the cup product on the cohomology of the mapping space, there is also a ``string topology" type product on the homology of certain mapping spaces. We now demonstrate how this string topology product may be modeled via the generalized Hochschild cohomology. In particular, we look at the case of surfaces $\Sigma^g$ of various genus $g$. The string topology product for this is then expressed as a map $H_\bullet(M^{\Sigma^{g_1}})\otimes H_\bullet(M^{\Sigma^{g_2}})\to H_\bullet(M^{\Sigma^{g_1+g_2}})$.

\subsection{A Hochschild model for the surface of genus $g$} \label{S:surfacemodel}
We start by giving a Hochschild model of the mapping space $\Map(\Sigma^g,M)$ from the surface of genus $g$ to a $2$-connected, compact, and oriented manifold $M$.

Recall, that the surface $\Sigma^g$ of genus $g\geq 1$ can be represented as a $4g$-gon, where the boundary is identified via the word $$ a_1 b_1a_2b_2\dots a_g b_g a_g^{-1} b_g^{-1}\dots a_2^{-1} b_2^{-1} a_1^{-1}b_1^{-1}. $$
We choose a subdivision of this polygon that fits with the string topology product. For this, we use a subdivision into $g^2$ squares, and further subdivide the off diagonal squares further into two triangles, so that for instance for $g=3$ we obtain

\begin{equation} \label{E:sigmag}
\begin{pspicture}(0,0)(4,4)
 \psline(0.5,.5)(0.5,3.5) \psline(.5,.5)(3.5,.5)
 \psline(1.5,.5)(1.5,3.5) \psline(.5,1.5)(3.5,1.5)
 \psline(2.5,.5)(2.5,3.5) \psline(.5,2.5)(3.5,2.5)
 \psline(3.5,.5)(3.5,3.5) \psline(.5,3.5)(3.5,3.5)
 \psline(.5,.5)(1.5,1.5) \psline(2.5,2.5)(3.5,3.5)
 \psline(.5,1.5)(2.5,3.5) \psline(1.5,.5)(3.5,2.5)
 \rput(0,3){$a_1$}  \rput(0,2){$b_1$} \rput(0,1){$a_2$}  \rput(1,0){$b_2$}
 \rput(2,0){$a_3$}  \rput(3,0){$b_3$} \rput(4,1){$a_3^{-1}$}  \rput(4,2){$b_3^{-1}$}
 \rput(4,3){$a_2^{-1}$}  \rput(3,4){$b_2^{-1}$} \rput(2,4){$a_1^{-1}$}  \rput(1,4){$b_1^{-1}$}
\end{pspicture}
\end{equation}
Each of the diagonal squares is represented via the simplicial model of the square $I^2_\bullet$, with $|I^2_k|=k^2+4k+4$. Also, for the square build out of two triangles, we glue two squares along an edge and collapse the opposite sides to a point,
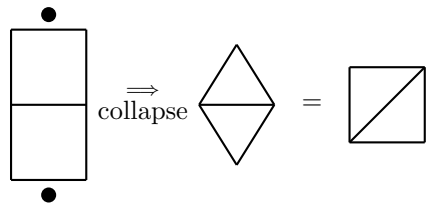
\[
\begin{pspicture}(0,0)(6.5,3)
 \psline(.5,.5)(1.5,.5)  \psline(.5,1.5)(1.5,1.5) \psline(.5,2.5)(1.5,2.5)
 \psline(1.5,.5)(1.5,2.5) \psline(.5,.5)(.5,2.5)
 \pscircle*(1,.3){.1} \pscircle*(1,2.7){.1}
 \rput(2.25,1.5){$\stackrel{\Longrightarrow}{\text{collapse}}$}
 \psline(3,1.5)(4,1.5) \psline(3,1.5)(3.5,2.3) \psline(3,1.5)(3.5,0.7)
 \psline(4,1.5)(3.5,2.3) \psline(4,1.5)(3.5,0.7)
 \rput(4.5,1.5){$=$}
 \psline(5,1)(5,2) \psline(6,1)(6,2) \psline(5,1)(6,1) \psline(5,2)(6,2)
 \psline(5,1)(6,2)
\end{pspicture}
\]
This is a model, which has in simplicial degree $k$ exactly $2k^2+5k+4$ elements. Gluing all these squares together by identifying the corresponding edges and identifying all vertices, we obtain a simplicial model $(\Sigma^g)_\bullet$ for the surface of genus $g$, with
\begin{equation}\label{E:kgenus} \#(\Sigma^g)_k=(2g^2-g)\cdot k^2+(3g^2-g)\cdot k+1+(g-1)^2. \end{equation}
We denote $\sigma^g_k=\#(\Sigma^g)_k-1$. 

\begin{notation}
We write $\Triangle_\bullet = {I^2_\bullet}/\sim$ for the simplicial model of a triangle obtained as a quotient of a square where a side is collapsed to a point.
\end{notation}

For genus $g=0$, we use the simplicial model $(\Sigma^0)_\bullet=S^2_\bullet$ of the sphere introduced in Example~\ref{E:2sphere}.
If $M$ is a $2$-connected, compact, oriented manifold, then Proposition \ref{qi} and Corollary \ref{qi*} imply that $\Ch^{(\Sigma^g)_\bullet}:CH^{(\Sigma^g)_\bullet}_\bullet(\Omega,\Omega)\to C^\bullet(M^{\Sigma^g})$ and $(\Ch^{(\Sigma^g)_\bullet})^*:C_\bullet(M^{\Sigma^g})\to CH_{(\Sigma^g)_\bullet}^\bullet(\Omega,\Omega^*)$ induce isomorphisms on (co)homology.

\medskip

The reason for studying this particular model of a surface of genus $g$ is that it comes with a simplicial description of  pinching maps. 
Pinching maps are obtained by collapsing to a point a circle, which contains the basepoint,  on a surface $\Sigma^{n}$ yielding a wedge $\Sigma^g\vee \Sigma^{h}$ (for any decomposition $n=g+h$).

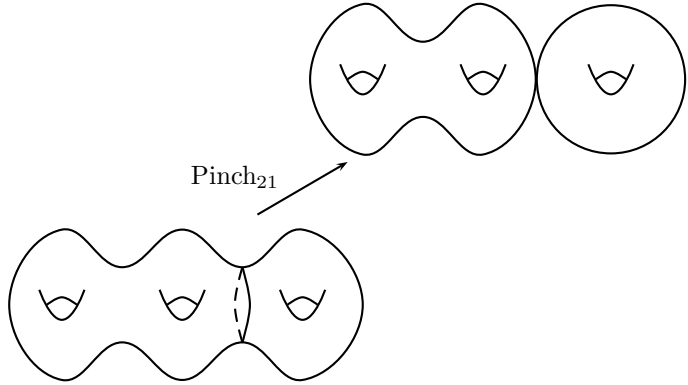
\[
\begin{pspicture}(1,0)(12.7,5)
  \psline[arrows=->](5.3,2.2)(6.5,2.9) \rput(5,2.7){$\pinch{21}$}
\psccurve(6,4)(6.7,5)(7.5,4.5)(8.3,5)(9,4)(8.3,3)(7.5,3.5)(6.7,3)
\pscurve(6.4,4.2)(6.7,3.8)(7,4.2) \pscurve(6.5,4)(6.7,4.1)(6.9,4)
\pscurve(8,4.2)(8.3,3.8)(8.6,4.2) \pscurve(8.1,4)(8.3,4.1)(8.5,4)
\pscircle(10,4){1}
\pscurve(9.7,4.2)(10,3.8)(10.3,4.2) \pscurve(9.8,4)(10,4.1)(10.2,4)

 \psccurve(2,1)(2.7,2)(3.5,1.5)(4.3,2)(5.1,1.5)(5.9,2)(6.7,1)(5.9,0)(5.1,.5)(4.3,0)(3.5,.5)(2.7,0)
 \pscurve(2.4,1.2)(2.7,0.8)(3,1.2) \pscurve(2.5,1)(2.7,1.1)(2.9,1)
 \pscurve(4,1.2)(4.3,0.8)(4.6,1.2) \pscurve(4.1,1)(4.3,1.1)(4.5,1)
 \pscurve(5.6,1.2)(5.9,0.8)(6.2,1.2) \pscurve(5.7,1)(5.9,1.1)(6.1,1)
 \pscurve(5.1,.5)(5.2,1)(5.1,1.5) \pscurve[linestyle=dashed](5.1,.5)(5,1)(5.1,1.5)
\end{pspicture}
\]

 This is realized on our simplicial model as follows. For $n=g+h$, we can consider 4 different regions in the model for $\Sigma^{g+h}$, 
namely we can consider the  top left square build out of $g^2$-squares (labelled $a_1,b_1,a_2, \cdots$ on the left and $b_1^{-1},a_1^{-1},\cdots$ on the top), the lower right square build out of $h^2$-squares (labelled $\cdots ,a_{g+h},b_{g+h}$ on the bottom  and $a_{g+h}^{-1},b_{g+g}^{-1}..$ on the right), and the two off diagonals rectangles denoted $R_b$, $R_t$. 
Note that all squares in the off diagonals regions $R_b$, $R_t$ are subdivided into triangles.    
Let \begin{equation} \label{eq:pinch} \pinch{g,h}:\Sigma^{g+h}_\bullet \to \Sigma^g_\bullet\vee \Sigma^h_\bullet \end{equation} be the map defined by identifying all the points in all triangles in $R_b$ and $R_t$ which belongs to a same parallel to the hypothenuse of the triangle (that is the edge parallel to the one which has been collapsed in the model). 
In other words, $\pinch{g,h}$ collapses, along the anti-diagonal, the two off diagonal regions $R_b$ and $R_t$ to the boundary of the top left and bottom right square.  For instance $\pinch{2,1}:\Sigma^{3}_\bullet \to \Sigma^2_\bullet\vee\Sigma^1_\bullet$ is given by the diagram
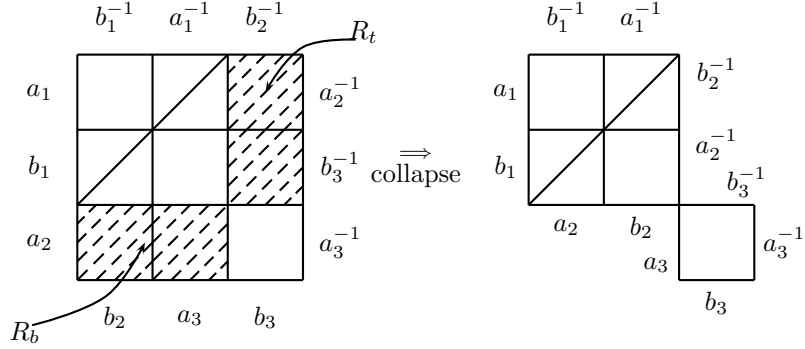
\[
\begin{pspicture}(-.5,-.5)(10,4.5)
 \rput(-.2, -.2){$R_b$}\rput(4.3, 3.8){$R_t$} \pscurve[arrows=->](-.1,-.1)(.9,.3)(1.4,1) 
 \pscurve[arrows=->](4.2,3.7)(3.5,3.6)(3,3)
 \psline(0.5,.5)(0.5,3.5) \psline(.5,.5)(3.5,.5)
 \psline(1.5,.5)(1.5,3.5) \psline(.5,1.5)(3.5,1.5)
\psline[linestyle=dashed](2.75,1.5)(3.5,2.25) \psline[linestyle=dashed](2.5,1.75)(3.25,2.5)
\psline[linestyle=dashed](3,1.5)(3.5,2) \psline[linestyle=dashed](2.5,2)(3,2.5)
\psline[linestyle=dashed](3.25,1.5)(3.5,1.75) \psline[linestyle=dashed](2.5,2.25)(2.75,2.5)
\psline[linestyle=dashed](1.5,0.5)(3.5,2.5)
\psline[linestyle=dashed](1.75,.5)(2.5,1.25) \psline[linestyle=dashed](1.5,.75)(2.25,1.5)
\psline[linestyle=dashed](2,.5)(2.5,1) \psline[linestyle=dashed](1.5,1)(2,1.5)
\psline[linestyle=dashed](2.25,.5)(2.5,.75)  \psline[linestyle=dashed](1.5,1.25)(1.75,1.5)
 \psline(2.5,.5)(2.5,3.5) \psline(.5,2.5)(3.5,2.5)
 \psline(3.5,.5)(3.5,3.5) \psline(.5,3.5)(3.5,3.5)
 \psline[linestyle=dashed](.5,.5)(1.5,1.5)
\psline[linestyle=dashed](.75,.5)(1.5,1.25) \psline[linestyle=dashed](.5,.75)(1.25,1.5)
\psline[linestyle=dashed](1,.5)(1.5,1)   \psline[linestyle=dashed](.5,1)(1,1.5)
\psline[linestyle=dashed](1.25,.5)(1.5,0.75) \psline[linestyle=dashed](.5,1.25)(.75,1.5)
 \psline[linestyle=dashed](2.5,2.5)(3.5,3.5)
\psline[linestyle=dashed](2.75,2.5)(3.5,3.25) \psline[linestyle=dashed](2.5,2.75)(3.25,3.5)
\psline[linestyle=dashed](3,2.5)(3.5,3) \psline[linestyle=dashed](2.5,3)(3,3.5)
\psline[linestyle=dashed](3.25,2.5)(3.5,2.75) \psline[linestyle=dashed](2.5,3.25)(2.75,3.5)
 \psline(.5,1.5)(2.5,3.5) 
 \rput(0,3){$a_1$}  \rput(0,2){$b_1$} \rput(0,1){$a_2$}  \rput(1,0){$b_2$}
 \rput(2,0){$a_3$}  \rput(3,0){$b_3$} \rput(4,1){$a_3^{-1}$}  \rput(4,2){$b_3^{-1}$}
 \rput(4,3){$a_2^{-1}$}  \rput(3,4){$b_2^{-1}$} \rput(2,4){$a_1^{-1}$}  \rput(1,4){$b_1^{-1}$}

\rput(5,2){$\stackrel{\Longrightarrow}{\text{collapse}}$}

\psline(6.5,1.5)(6.5,3.5) \psline(8.5,0.5)(9.5,.5) \psline(6.5,1.5)(8.5,3.5)
 \psline(7.5,1.5)(7.5,3.5) \psline(6.5,1.5)(9.5,1.5)
 \psline(8.5,.5)(8.5,3.5) \psline(6.5,2.5)(8.5,2.5)
 \psline(9.5,.5)(9.5,1.5) \psline(6.5,3.5)(8.5,3.5)
\rput(6.2,3){$a_1$}  \rput(6.2,2){$b_1$} \rput(7,1.2){$a_2$}  \rput(8,1.2){$b_2$}
 \rput(8.2,0.7){$a_3$}  \rput(9,0.2){$b_3$} \rput(9.9,1){$a_3^{-1}$}  
\rput(9.4,1.8){$b_3^{-1}$}
 \rput(9,2.3){$a_2^{-1}$}  \rput(9,3.25){$b_2^{-1}$} \rput(8,4){$a_1^{-1}$}  \rput(7,4){$b_1^{-1}$}
\end{pspicture} 
\]
where all elements in the same dashed line are identified, \emph{i.e.} collapsed to the same point.
 
\begin{lem}\label{L:modelglueing}
The map $\pinch{g,h}: \Sigma^{g+h}_\bullet \to \Sigma^g_\bullet\vee \Sigma^h_\bullet $ is simplicial.
\end{lem}
 \begin{proof}
The map  $\pinch{g,h}$ is obtained by wedging along an edge, or a vertex, maps such as the identity $id:I^2_\bullet \to I^2_\bullet$, $id:\Triangle_\bullet \to \Triangle_\bullet$ and collapse $\Triangle_\bullet \to pt_\bullet$ of a triangle onto a point or collapse $\Triangle_\bullet\cong {I^2_\bullet}/\sim \stackrel{p_j}\to I_\bullet$ of a triangle onto one edge (which has not been identified to a point). Now it follows from Example~\ref{E:wedge} and Lemma~\ref{L:collapse} that $\pinch{g,h}$ is a map of (pointed) simplicial sets.
 \end{proof}

\begin{remk}
It is crucial to use the simplicial model $\Sigma^g_\bullet$ described in Lemma~\ref{L:modelglueing}. For instance, if one uses a model where the off diagonal squares are not subdivided, the Lemma~\ref{L:modelglueing} above is no longer true.
\end{remk}
The following Lemma is a straightforward.
\begin{lem}
The simplicial pinching map is associative, \emph{i.e.} the following diagram is commutative   $$\xymatrix{\Sigma^{g+h+k}_\bullet  \ar[rr]^{\pinch{g+h,k}} \ar[d]_{\pinch{g,h+k}} && \Sigma^{g+h}_\bullet \vee \Sigma^k_\bullet \ar[d]^{\pinch{g,h}\vee id_{\Sigma^k_\bullet}} \\
\Sigma^{g}_\bullet \vee \Sigma^{h+k}_\bullet \ar[rr]^{ id_{\Sigma^k_\bullet}\vee \pinch{h,k}}&& \Sigma^{g}_\bullet \vee \Sigma^h_\bullet  \vee \Sigma^k_\bullet} $$
\end{lem}

\subsection{The ``string topology" product for surfaces}\label{thom}
In this section, we recall the ``string topology" type operation adapted for surfaces, and then apply this to the model for the surface mapping space  given in the previous subsection. We start by recalling this operation, which was originally given for the mapping space of a circle by Moira Chas and Dennis Sullivan in \cite{CS}, see also the description of the Cohen-Jones map generalized to surfaces as given in \cite[Section 5.2]{CV} for the $k$-sphere.

\medskip
In this section, we use the model $\Sigma^{g}$ ($=|\Sigma^g_\bullet|$) for a (compact oriented) surface of genus $g$ introduced in Section~\ref{S:surfacemodel} with its basepoint $\bpt$. 

Denote by $\Map(\Sigma^g,M)$ the space of (continuous, non pointed) maps from a surface $\Sigma^g$  to the manifold $M$, which we assume to be compact and oriented. For two such surfaces $\Sigma^{g}$ and $\Sigma^{h}$, denote by $\Sigma^g\vee \Sigma^{h}$ their wedge product, \emph{i.e.} their disjoint union  modulo the identification of the two basepoints. The space $\Map(\Sigma^g\vee \Sigma^h,M)$ denotes the corresponding mapping space from $\Sigma^{g}\vee \Sigma^{h}$ to $M$. Then there are induced maps
$$ \Map(\Sigma^g,M)\times \Map(\Sigma^{h},M) \stackrel {\rho_{in}} \leftarrow \Map(\Sigma^{g}\vee \Sigma^{h},M) \stackrel {\rho_{out}}\rightarrow \Map(\Sigma^{g+h},M), $$
where $\rho_{in}$ is given by including to the first and second component in $\Sigma^{g}\vee \Sigma^{h}$. For surfaces with positive genus,  $\rho_{out}$ is induced by the pinching map $\pinch{g,h}:\Sigma^{g+h}\to \Sigma^{g}\vee \Sigma^{h}$ given by the geometric realization of the simplicial map~\eqref{eq:pinch} $\pinch{g,h}: \Sigma^{g+h}_\bullet \to \Sigma^g_\bullet\vee \Sigma^h_\bullet $.  If $g=0$, another model for $\Sigma^h$ is given by gluing 4 squares 
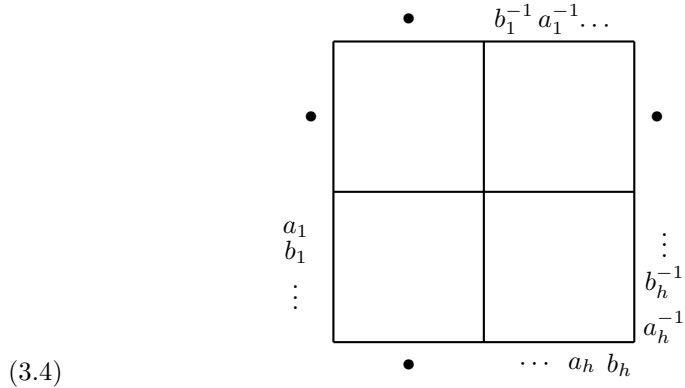
\begin{equation}\label{E:sigma0h}
\begin{pspicture}(0,0)(5,5)
 \psline(0.5,.5)(0.5,4.5) \psline(.5,.5)(4.5,.5)
 \psline(4.5,4.5)(4.5,0.5) \psline(.5,4.5)(4.5,4.5)
 \psline(2.5,.5)(2.5,4.5) \psline(.5,2.5)(4.5,2.5)
 \rput(1.5,4.8){$\bullet$}  \rput(1.5,0.2){$\bullet$} \rput(4.8,3.5){$\bullet$}  \rput(0.2,3.5){$\bullet$}\rput(0,2,3){$a_1$}  \rput(0,1.7){$b_1$} \rput(0,1.2){$\vdots$} \rput(3.2,0.2){$\cdots$} \rput(3.8,0.2){$a_h$} \rput(4.3,0.2){$b_h$} 
\rput(4.9,0.7){$a_h^{-1}$} \rput(4.9,1.3){$b_h^{-1}$} \rput(4.9,1.9){$\vdots$}
\rput(2.9,4.8){$b_1^{-1}$}  \rput(3.5,4.8){$a_1^{-1}$} \rput(4,4.7){$\cdots$}
\end{pspicture}
\end{equation}
where the bulleted edges are collapsed to a point and the other boundary edges are identified with the word  $a_1 b_1a_2b_2\dots a_h b_h a_h^{-1} b_h^{-1}\dots a_2^{-1} b_2^{-1} a_1^{-1}b_1^{-1}$ as for the model~\eqref{E:sigmag}. The pinching map $\pinch{0,h}:\Sigma^h \to \Sigma^0\vee \Sigma^h={\Sigma}^{0}\vee \Sigma^{h}$ is given by 

\begin{equation}\label{E:sigma0hcollapse}
\begin{pspicture}(0,0)(12,5)
 \psline(0.5,.5)(0.5,4.5) \psline(.5,.5)(4.5,.5)
 \psline(4.5,4.5)(4.5,0.5) \psline(.5,4.5)(4.5,4.5)
\psline[linestyle=dashed](0.5,0.75)(2.5,0.75)\psline[linestyle=dashed](0.5,1)(2.5,1)
\psline[linestyle=dashed](0.5,1.3)(2.5,1.3) \psline[linestyle=dashed](0.5,1.6)(2.5,1.6)
\psline[linestyle=dashed](0.5,1.9)(2.5,1.9) \psline[linestyle=dashed](0.5,2.2)(2.5,2.2)
 \psline(2.5,.5)(2.5,4.5) \psline(.5,2.5)(4.5,2.5)
\psline[linestyle=dashed](2.75,2.5)(2.75,4.5)\psline[linestyle=dashed](3,2.5)(3,4.5)
\psline[linestyle=dashed](3.3,2.5)(3.3,4.5) \psline[linestyle=dashed](3.6,2.5)(3.6,4.5)
\psline[linestyle=dashed](3.9,2.5)(3.9,4.5) \psline[linestyle=dashed](4.2,2.5)(4.2,4.5)

 \rput(1.5,4.8){$\bullet$}  \rput(1.5,0.2){$\bullet$} \rput(4.8,3.5){$\bullet$}  \rput(0.2,3.5){$\bullet$}\rput(0,2,3){$a_1$}  \rput(0,1.7){$b_1$} \rput(0,1.2){$\vdots$} \rput(3,0.2){$\cdots$} \rput(3.8,0.2){$a_h$} \rput(4.3,0.2){$b_h$} 
\rput(4.9,0.7){$a_h^{-1}$} \rput(4.9,1.3){$b_h^{-1}$} \rput(4.9,1.9){$\vdots$}
\rput(2.9,4.8){$b_1^{-1}$}  \rput(3.5,4.8){$a_1^{-1}$} \rput(4,4.7){$\cdots$} 

\rput(6.2,2,3){$\stackrel{\Longrightarrow}{\text{collapse}}$}

\psline(7.5,2.5)(7.5,4.5) \psline(9.5,.5)(11.5,.5)
 \psline(11.5,2.5)(11.5,0.5) \psline(7.5,4.5)(9.5,4.5)
 \psline(9.5,.5)(9.5,4.5) \psline(7.5,2.5)(11.5,2.5)
 \rput(8.5,4.8){$\bullet$}  \rput(8.5,2.2){$\bullet$} \rput(9.8,3.5){$\bullet$}  \rput(7.2,3.5){$\bullet$}\rput(9.1,2,3){$a_1$}  \rput(9.1,1.7){$b_1$} \rput(9.1,1.2){$\vdots$} \rput(10,0.2){$\cdots$} \rput(10.8,0.2){$a_h$} \rput(11.3,0.2){$b_h$} 
\rput(11.9,0.7){$a_h^{-1}$} \rput(11.9,1.3){$b_h^{-1}$} \rput(11.9,1.9){$\vdots$}
\rput(9.9,2.8){$b_1^{-1}$}  \rput(10.5,2.8){$a_1^{-1}$} \rput(11,2.7){$\cdots$} 
\end{pspicture}
\end{equation}
where all elements in the same dashed line are identified, \emph{i.e.} collapsed to the same point. There is a similar pinching map $\pinch{0,h}:\Sigma^h \to  \Sigma^h\vee \Sigma^0$ given by
 \begin{equation}\label{E:sigmag0collapse}
\begin{pspicture}(0,0)(12,5)
 \psline(0.5,.5)(0.5,4.5) \psline(.5,.5)(4.5,.5)
 \psline(4.5,4.5)(4.5,0.5) \psline(.5,4.5)(4.5,4.5)
\psline[linestyle=dashed](2.5,2.75)(4.5,2.75)\psline[linestyle=dashed](2.5,3)(4.5,3)
\psline[linestyle=dashed](2.5,3.3)(4.5,3.3) \psline[linestyle=dashed](2.5,3.6)(4.5,3.6)
\psline[linestyle=dashed](2.5,3.9)(4.5,3.9) \psline[linestyle=dashed](2.5,4.2)(4.5,4.2)
 \psline(2.5,.5)(2.5,4.5) \psline(.5,2.5)(4.5,2.5)
\psline[linestyle=dashed](.75,.5)(.75,2.5)\psline[linestyle=dashed](1,.5)(1,2.5)
\psline[linestyle=dashed](1.3,.5)(1.3,2.5) \psline[linestyle=dashed](1.6,.5)(1.6,2.5)
\psline[linestyle=dashed](1.9,.5)(1.9,2.5) \psline[linestyle=dashed](2.2,.5)(2.2,2.5)

 \rput(3.5,4.8){$\bullet$}  \rput(3.5,0.2){$\bullet$} \rput(4.8,1.5){$\bullet$}  \rput(0.2,1.5){$\bullet$}\rput(0,4.3){$a_1$}  \rput(0,3.8){$b_1$} \rput(0,3.2){$\vdots$} \rput(1,0.2){$\cdots$} \rput(1.8,0.2){$a_h$} \rput(2.3,0.2){$b_h$} 
\rput(4.9,2.8){$a_h^{-1}$} \rput(4.9,3.4){$b_h^{-1}$} \rput(4.9,4.1){$\vdots$}
\rput(0.9,4.8){$b_1^{-1}$}  \rput(1.5,4.8){$a_1^{-1}$} \rput(2,4.7){$\cdots$} 

\rput(6.2,2,3){$\stackrel{\Longrightarrow}{\text{collapse}}$}

\psline(7.5,2.5)(7.5,4.5) \psline(9.5,.5)(11.5,.5)
 \psline(11.5,2.5)(11.5,0.5) \psline(7.5,4.5)(9.5,4.5)
 \psline(9.5,.5)(9.5,4.5) \psline(7.5,2.5)(11.5,2.5)
 \rput(10.5,2.8){$\bullet$}  \rput(10.5,0.2){$\bullet$} \rput(11.8,1.5){$\bullet$}  \rput(9.2,1.5){$\bullet$}\rput(7.1,4.3){$a_1$}  \rput(7.1,3.8){$b_1$} \rput(7.1,3.2){$\vdots$} \rput(8,2.2){$\cdots$} \rput(8.7,2.2){$a_h$} \rput(9.2,2.2){$b_h$} 
\rput(9.9,2.8){$a_h^{-1}$} \rput(9.9,3.4){$b_h^{-1}$} \rput(9.9,4.2){$\vdots$}
\rput(7.9,4.8){$b_1^{-1}$}  \rput(8.5,4.8){$a_1^{-1}$} \rput(9,4.7){$\cdots$} 
\end{pspicture}
\end{equation}
By collapsing all boundary edges to a point in the model~\eqref{E:sigma0h} and in definition of the map $\pinch{0,h}$~\eqref{E:sigma0hcollapse} yields the usual pinching map $\pinch{0,0}:\Sigma^0\to \Sigma^0\vee \Sigma^0$ for the dimension 2 sphere $S^2=\Sigma^0$. The pinching maps $\pinch{0,\bullet}$, $\pinch{\bullet,0}$ above induce 
$\rho_{out}$ when one of the surfaces has genus zero.

\medskip

Note that the map $\rho_{in}$ is given as a pullback of diagrams
\begin{equation*}
\xymatrix{
  \Map(\Sigma^{g}\vee \Sigma^{h},M) \ar[r]^{\qquad \rho_{in}\quad\quad\quad\quad} \ar[d] 
  & \Map(\Sigma^g,M)\times \Map(\Sigma^{h},M) \ar[d] \\
 M  \ar[r]^{\text{diagonal}\quad} & M\times M }
\end{equation*}
In particular, $\rho_{in}$ is an embedding of infinite dimensional manifolds with finite codimension equal to the dimension of $M$, $codim(\rho_{in})=dim(M)$ and the associated normal bundle is of dimension $dim(M)$ and oriented (since $M$ is). Thus, if we denote by $\Map(\Sigma^{g}\vee \Sigma^{h},M)^{-TM}$ the Thom space of this embedding, there is a  Thom class in $H^{m}(\Map(\Sigma^{g}\vee \Sigma^{h},M)^{-TM})$ inducing  the Thom isomorphism $$t:H_\bullet(\Map(\Sigma^{g}\vee \Sigma^{h},M)^{-TM})\stackrel{\cong} {\rightarrow}H_{\bullet-m}(\Map(\Sigma^{g}\vee \Sigma^{h},M)),$$ where $m=dim(M)$.  Together with the Thom collapse map $\tau:\Map(\Sigma^g,M)\times \Map(\Sigma^{h},M)\to \Map(\Sigma^{g}\vee \Sigma^{h},M)^{-TM}$,  we obtain the following Umkehr map,
\begin{multline*}
(\rho_{in})_! : H_\bullet(\Map(\Sigma^g,M))\otimes H_{\bullet}(\Map(\Sigma^{h},M))\\
\cong H_\bullet(\Map(\Sigma^g,M)\times \Map(\Sigma^{h},M))
 \\ \stackrel{\tau_*}{\longrightarrow} H_\bullet(\Map(\Sigma^{g}\vee \Sigma^{h},M)^{-TM})\stackrel{t}{\longrightarrow}H_{\bullet-m}(\Map(\Sigma^{g}\vee \Sigma^{h},M)).
\end{multline*}
\begin{defn} \label{D:surfaceproduct} With this, we define the product $\scup := (\rho_{out})_* \circ(\rho_{in})_!$,
\begin{multline*}
 \scup: H_\bullet(\Map(\Sigma^g,M))\otimes H_\bullet(\Map(\Sigma^{h},M)) \\ \stackrel{(\rho_{in})_!}{\longrightarrow} H_\bullet(\Map(\Sigma^{g}\vee \Sigma^{h},M))\stackrel{(\rho_{out})_*}{\longrightarrow} H_\bullet(\Map(\Sigma^{g+h},M))
\end{multline*} that we call the \emph{surface product}.
\end{defn}
We denote  $\mathbb{H}_{\bullet}(\Map(\Sigma^g,M))$ the shifted homology groups $H_{\bullet+dim(M)}(\Map(\Sigma^g,M))$. This shifting makes the surface product a degree zero map.
\begin{thm}\label{T:surfaceproduct}
Let $M$ be a compact oriented manifold. Then the surface product  $\scup:\mathbb{H}_{\bullet}(\Map(\Sigma^g,M))\otimes \mathbb{H}_{\bullet}(\Map(\Sigma^{h},M)) \to \mathbb{H}_{\bullet}(\Map(\Sigma^{g+h},M))$ is associative.
\end{thm}
\begin{proof}
 It is well-known that $\pinch{0,0}:\Sigma^0\to \Sigma^0\vee \Sigma^0$ is homotopy associative. From there and Lemma~\ref{L:modelglueing} it follows that $(\rho_{out})_*:H_\bullet(\Map(\Sigma^{g}\vee \Sigma^{h},M))\stackrel{(\rho_{out})_*}{\longrightarrow} H_\bullet(\Map(\Sigma^{g+h},M))$ is associative. Now the Theorem follows from naturality property of Umkehr maps: $(\rho_{out})_*\circ (\rho_{in})_!=(\rho_{in})_!\circ (\rho_{out})_* $ as in the usual string topology case (see~\cite{CJ, BGNX} for details). 
\end{proof}

\begin{remk}
Note that the surface product gives a structure of associative graded (with respect to the genus)  algebra with unit (see Proposition~\ref{P:bimodule}) to $$\mathbb{H}_\bullet\big(\mathop{\bigsqcup}_{g\geq 0} \Map\big(\Sigma^{g},M\big)\big)\cong \bigoplus_{g\geq 0} \mathbb{H}_{\bullet}(\Map(\Sigma^{g},M)).$$
\end{remk}

There is an obvious embedding $i_g$ of $M$ into $\Map(\Sigma^g, M)$ as constant functions. Thus, for any $g\geq 0$, the fundamental class of $M$ yields a class $$[M]_g=i_g([M]) \in \mathbb{H}_0(\Map(\Sigma^g,M)).$$ 

\smallskip

Also note that for genus zero, $\Sigma^0\cong S^2$, the surface product restricts to a product $\mathbb{H}_{\bullet}(\Map(\Sigma^0,M))\otimes \mathbb{H}_{\bullet}(\Map(\Sigma^{0},M)) \to \mathbb{H}_{\bullet}(\Map(\Sigma^{0},M))$. This product is the usual (dimension 2) Brane Topology product:
\begin{prop}\label{P:Brane}
The restriction of the surface product (see Definition~\ref{D:surfaceproduct}) to $\mathbb{H}_{\bullet}(\Map(\Sigma^0,M))$ coincides with the Brane topology product $\mathbb{H}_{\bullet}(\Map(S^2,M))^{\otimes 2} \to \mathbb{H}_{\bullet}(\Map(S^2,M))$ see~\cite{H, CV}. In particular it is graded commutative with $[M]_0$ as a unit.  
\end{prop}
\begin{proof}
The Brane Topology product (for dimension $2$-spheres) as defined in~\cite[Section 5]{CV} is induced by a structure of algebra over the homology $H_\bullet(\mathfrak{Cac})$ of the (2-dimensional) cactus operad $\mathfrak{Cac}$ on $\mathbb{H}_{\bullet}(\Map(\Sigma^0,M))$. By definition,  an element $c\in \mathfrak{Cac}(2)$ is  a map $c:S^2\to S^2\vee S^2$. Thus it induces a map $\rho_{in}(c): \Map(S^2\vee S^2,M)\to \Map(S^2,M)$. The Brane Topology product~\cite[Section 5.2]{CV} is then given by the composition $\rho_{out}\circ (\rho_{in}(c))_!$ for any cactus $c\in \mathfrak{Cac}(2)$. The result follows by choosing $c=\pinch{0,0}$.
\end{proof}

By Theorem~\ref{T:surfaceproduct}, $\mathbb{H}_\bullet(\Map(\Sigma^g,M))$ inherits a left  $\mathbb{H}_\bullet(\Map(\Sigma^0,M))$-module structure $\mathbb{H}_\bullet(\Map(\Sigma^0,M))\otimes \mathbb{H}_\bullet(\Map(\Sigma^g,M))\stackrel{\scup}\to \mathbb{H}_\bullet(\Map(\Sigma^g,M))$ as well as a right module structure. 
\begin{prop}\label{P:bimodule}
 $\mathbb{H}_\bullet(\Map(\Sigma^g,M))$ is a (graded) symmetric $\mathbb{H}_\bullet(\Map(\Sigma^0,M))$-bimodule, \emph{i.e.} for any $x\in\mathbb{H}_\bullet(\Map(\Sigma^g,M))$, $y\in \mathbb{H}_\bullet(\Map(\Sigma^0,M))$, one has
$$[M]_0\scup x=x \qquad \mbox{ and } \qquad x\scup y= (-1)^{|y|\cdot |x|} y\scup x.$$
\end{prop}
\begin{proof}
Note that there is a commutative diagram of pullbacks
$$
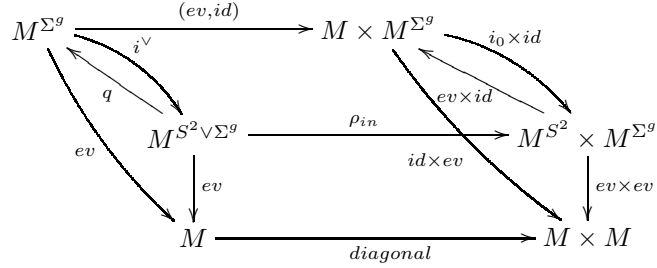
\xymatrix{M^{\Sigma^g}\ar[rr]^{(ev,id)}\ar@/^.8pc/[rd]^{i^\vee}\ar@/_/[rdd]_{ev}& & M\times M^{\Sigma^g}  \ar@/^1pc/[rd]^{i_0\times id}   \ar@/_/[rdd]_{id\times ev}& \\ 
  & M^{S^2\vee \Sigma^g} \ar[rr]^{\rho_{in}\qquad} \ar[ul]^{q}\ar[d]^{ev}& & M^{S^2}\times M^{\Sigma^g}\ar[d]^{ev\times ev}  \ar[ul]^{ev\times id}\\
  & M \ar[rr]_{diagonal}& &M\times M}
 $$ 
 where $q$ is induced by the inclusion $\Sigma^g\hookrightarrow S^2\vee \Sigma^g$ and $i^\vee$ is induced by $S^2\vee \Sigma^g\to \Sigma^g$, which collapses the $S^2$ component to a point. Since $[M]_0=i_0([M])$, it follows that $(\rho_{in})^!([M]_0\times x)=i^\vee_*\big((ev,id)^!([M]\times x)\big)= i^\vee(x)$ and the identity $[M]_0\scup x=x $ follows, since $\pinch{0,g}\circ i^\vee$ is homotopic to the identity. 

\smallskip

It remains to show that $\pinch{0,g}$ and $\pinch{g,0}$ are homotopic maps $\Sigma^g \to \Sigma^0\vee \Sigma^g$. For each $t\in[0,1]$, there is a parametrization of $\Sigma^g$ obtained by attaching standard squares and rectangles and identifying some  boundary components as in the following figure:

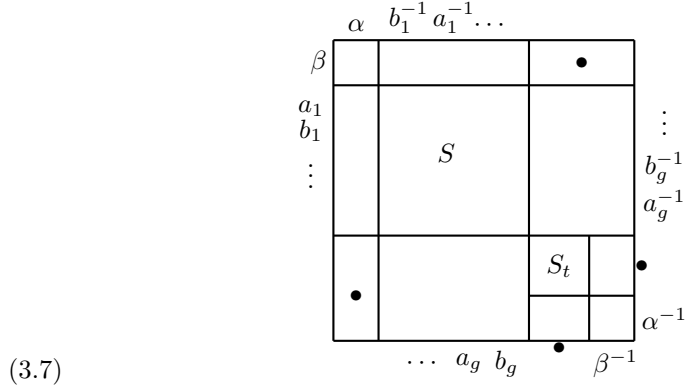
\begin{equation} \label{eq:sigmagt}
\begin{pspicture}(0,0)(5,5)
 \psline(0.5,.5)(0.5,4.5) \psline(.5,.5)(4.5,.5)
 \psline(4.5,4.5)(4.5,0.5) \psline(.5,4.5)(4.5,4.5)
 \psline(1.1,.5)(1.1,4.5) \psline(3.1,.5)(3.1,4.5)
 \psline(.5,1.9)(4.5,1.9) \psline(.5,3.9)(4.5,3.9)
 \psline(3.1,1.1)(4.5,1.1) \psline(3.9,.5)(3.9,1.9)
\rput(0.2,3.6){$a_1$}   \rput(0.2,3.3){$b_1$} \rput(0.2,2.8){$\vdots$} \rput(1.7,0.2){$\cdots$} \rput(2.3,0.2){$a_g$} \rput(2.8,0.2){$b_g$} 
\rput(4.9,2.3){$a_g^{-1}$} \rput(4.9,2.8){$b_g^{-1}$} \rput(4.9,3.5){$\vdots$}
\rput(1.5,4.8){$b_1^{-1}$}  \rput(2.1,4.8){$a_1^{-1}$} \rput(2.6,4.7){$\cdots$}
\rput(0.3,4.2){$\beta$} \rput(0.8,4.7){$\alpha$}
\rput(4.25,0.2){$\beta^{-1}$} \rput(4.9,.8){$\alpha^{-1}$}
\rput(3.5,.4){$\bullet$}  \rput(4.6,1.5){$\bullet$} 
\rput(.8,1.1){$\bullet$} \rput(3.8,4.2){$\bullet$} 
\rput(2,3){$S$} \rput(3.5,1.5){$S_{t}$}
\end{pspicture}
\end{equation}
More precisely, the big central square $S$ in figure~\eqref{eq:sigmagt} is a standard square $[0,1]^2$ while the small central square $S_t$ has edges of length $t$ (thus  $S_t=[0,t]^2$). The edges labelled by $a_i$s, $b_j$s and their inverses are identified just in the usual model for $\Sigma^g$, see Figures~\eqref{E:sigma0h} and~\eqref{E:sigmag}. The edges $\alpha$ and $\beta$ are of length $(1-t)/2$. The two bulleted edges are identified to the base point and the (top right and bottom left) bulleted rectangles are entirely collapsed to the  base point as well. Overall, the parametrization pictured by figure~\eqref{eq:sigmagt} is the square $[0,2]^2$ with some boundary elements identified and two sub-rectangles collapsed to the base point.

We now define a pinching map $Pinch(t,-)$ from $\Sigma^g$ to $\Sigma^0\vee \Sigma^g$. Similar to the pinching maps $\pinch{0,g}$ and $\pinch{g,0}$, $Pinch(t,-)$ is obtained by collapsing some elements in the above parametrization~\eqref{eq:sigmagt} of $\Sigma^g$ to the base point, as shown in the following picture:

\begin{equation*} 
\begin{pspicture}(0,0)(12,5.2)
 \psline(0.5,.5)(0.5,4.5) \psline(.5,.5)(4.5,.5)
 \psline(4.5,4.5)(4.5,0.5) \psline(.5,4.5)(4.5,4.5)
 \psline(1.1,.5)(1.1,4.5) \psline(3.1,.5)(3.1,4.5)
 \psline(.5,1.9)(4.5,1.9) \psline(.5,3.9)(4.5,3.9)
 \psline(3.1,1.1)(4.5,1.1) \psline(3.9,.5)(3.9,1.9)
\rput(0.2,3.6){$a_1$}   \rput(0.2,3.3){$b_1$} \rput(0.2,2.8){$\vdots$} \rput(1.7,0.2){$\cdots$} \rput(2.3,0.2){$a_g$} \rput(2.8,0.2){$b_g$} 
\rput(4.9,2.3){$a_g^{-1}$} \rput(4.9,2.8){$b_g^{-1}$} \rput(4.9,3.5){$\vdots$}
\rput(1.5,4.8){$b_1^{-1}$}  \rput(2.1,4.8){$a_1^{-1}$} \rput(2.6,4.7){$\cdots$}
\rput(0.3,4.2){$\beta$} \rput(0.8,4.7){$\alpha$}
\rput(4.25,0.2){$\beta^{-1}$} \rput(4.9,.8){$\alpha^{-1}$}
\rput(3.5,.4){$\bullet$}  \rput(4.6,1.5){$\bullet$} 
\rput(.8,1.1){$\bullet$} \rput(3.8,4.2){$\bullet$}  

\psline[linestyle=dashed](.5,2.3)(1.1,2.3)\psline[linestyle=dashed](.5,2.7)(1.1,2.7)
\psline[linestyle=dashed](.5,3.1)(1.1,3.1)\psline[linestyle=dashed](.5,3.5)(1.1,3.5)
\psline[linestyle=dashed](3.1,2.3)(4.5,2.3)\psline[linestyle=dashed](3.1,2.7)(4.5,2.7)
\psline[linestyle=dashed](3.1,3.1)(4.5,3.1)\psline[linestyle=dashed](3.1,3.5)(4.5,3.5)

\psline[linestyle=dashed](1.5,4.5)(1.5,3.9)\psline[linestyle=dashed](1.9,4.5)(1.9,3.9)
\psline[linestyle=dashed](2.3,4.5)(2.3,3.9)\psline[linestyle=dashed](2.7,4.5)(2.7,3.9)
\psline[linestyle=dashed](1.5,.5)(1.5,1.9)\psline[linestyle=dashed](1.9,.5)(1.9,1.9)
\psline[linestyle=dashed](2.3,.5)(2.3,1.9)\psline[linestyle=dashed](2.7,.5)(2.7,1.9)

\rput(6.3,2.3){$\stackrel{\stackrel{Pinch(t,-)}\Longrightarrow}{\text{collapse}}$}

 \psline(7.5,3.9)(7.5,4.5) 
  \psline(7.5,4.5)(8.1,4.5)
 \psline(8.1,1.9)(8.1,4.5) \psline(10.1,0.5)(10.1,3.9)
 \psline(8.1,1.9)(11.5,1.9) \psline(7.5,3.9)(10.1,3.9)
 \psline(10.1,0.5)(11.5,0.5) \psline(11.5,0.5)(11.5,1.9)
 \psline(10.9,0.5)(10.9,1.1) \psline(10.9,1.1)(11.5,1.1)
\rput(7.8,3.5){$a_1$}   \rput(7.8,3.2){$b_1$} \rput(7.8,2.7){$\vdots$} \rput(8.6,1.6){$\cdots$} \rput(9.2,1.6){$a_g$} \rput(9.7,1.6){$b_g$} 
\rput(10.5,2.4){$a_g^{-1}$} \rput(10.5,2.9){$b_g^{-1}$} \rput(10.5,3.6){$\vdots$}
\rput(8.6,4.2){$b_1^{-1}$}  \rput(9.2,4.2){$a_1^{-1}$} \rput(9.7,4.1){$\cdots$}
\rput(7.3,4.2){$\beta$} \rput(7.8,4.7){$\alpha$}
\rput(11.25,0.2){$\beta^{-1}$} \rput(11.9,.8){$\alpha^{-1}$}
\rput(7.75,3.8){$\bullet$}  \rput(8.2,4.3){$\bullet$} 
\rput(10,1.1){$\bullet$}  \rput(10.9,2){$\bullet$} 
\rput(10.5,.4){$\bullet$}  \rput(11.6,1.5){$\bullet$} 
\end{pspicture}
\end{equation*}
Here all elements in the same dashed line get identified by the pinching map $Pinch(t,-)$, \emph{i.e.} they get collapsed to the same point, and all bulleted rectangles and egdes get collapsed to a point.

\smallskip

Note that $Pinch(1,-) =\pinch{g,0}$.
Thus the map $Pinch(-,-): [0,1]\times \Sigma^g\to \Sigma^0\vee \Sigma^g$ is an homotopy between $\pinch{g,0}$ and $Pinch(0,-)$ which is the collapse map given by the following figure :
\begin{equation*} 
\begin{pspicture}(0,0)(12,5)
 \psline(0.5,.5)(0.5,4.5) \psline(.5,.5)(4.5,.5)
 \psline(4.5,4.5)(4.5,0.5) \psline(.5,4.5)(4.5,4.5)
 \psline(1.5,.5)(1.5,4.5) \psline(3.5,.5)(3.5,4.5)
 \psline(.5,1.5)(4.5,1.5) \psline(.5,3.5)(4.5,3.5)
 
\rput(0.2,3.2){$a_1$}   \rput(0.2,2.9){$b_1$} \rput(0.2,2.4){$\vdots$} \rput(2.1,0.2){$\cdots$} \rput(2.7,0.2){$a_g$} \rput(3.2,0.2){$b_g$} 
\rput(4.9,1.9){$a_g^{-1}$} \rput(4.9,2.4){$b_g^{-1}$} \rput(4.9,3.1){$\vdots$}
\rput(1.9,4.8){$b_1^{-1}$}  \rput(2.5,4.8){$a_1^{-1}$} \rput(3,4.7){$\cdots$}
\rput(0.3,4.1){$\beta$} \rput(0.9,4.7){$\alpha$}
\rput(4.1,0.2){$\beta^{-1}$} \rput(4.9,.9){$\alpha^{-1}$}
\rput(1,1){$\bullet$} \rput(4,4){$\bullet$}

\psline[linestyle=dashed](.5,1.9)(1.5,1.9)\psline[linestyle=dashed](.5,2.3)(1.5,2.3)
\psline[linestyle=dashed](.5,2.7)(1.5,2.7)\psline[linestyle=dashed](.5,3.1)(1.5,3.1)
\psline[linestyle=dashed](3.5,1.9)(4.5,1.9)\psline[linestyle=dashed](3.5,2.3)(4.5,2.3)
\psline[linestyle=dashed](3.5,2.7)(4.5,2.7)\psline[linestyle=dashed](3.5,3.1)(4.5,3.1)

\psline[linestyle=dashed](1.9,4.5)(1.9,3.5)\psline[linestyle=dashed](2.3,4.5)(2.3,3.5)
\psline[linestyle=dashed](2.7,4.5)(2.7,3.5)\psline[linestyle=dashed](3.1,4.5)(3.1,3.5)
\psline[linestyle=dashed](3.1,.5)(3.1,1.5)\psline[linestyle=dashed](1.9,.5)(1.9,1.5)
\psline[linestyle=dashed](2.3,.5)(2.3,1.5)\psline[linestyle=dashed](2.7,.5)(2.7,1.5)

\rput(6.3,2.3){$\stackrel{\stackrel{Pinch(0,-)}\Longrightarrow}{\text{collapse}}$}

 \psline(7.5,3.5)(7.5,4.5) 
  \psline(7.5,4.5)(8.5,4.5) \psline(7.5,3.5)(10.5,3.5)
 \psline(8.5,1.5)(8.5,4.5) \psline(10.5,.5)(10.5,3.5)
 \psline(8.5,1.5)(11.5,1.5) \psline(11.5,.5)(11.5,1.5) \psline(10.5,.5)(11.5,.5)
\rput(8.2,3.1){$a_1$}   \rput(8.2,2.8){$b_1$} \rput(8.2,2.3){$\vdots$} \rput(8.9,1.2){$\cdots$} \rput(9.5,1.2){$a_g$} \rput(10,1.2){$b_g$} 
\rput(11,2.1){$a_g^{-1}$} \rput(11,2.6){$b_g^{-1}$} \rput(11,3.3){$\vdots$}
\rput(9.2,3.8){$b_1^{-1}$}  \rput(9.8,3.8){$a_1^{-1}$} \rput(10.3,3.7){$\cdots$}
\rput(7.2,4.1){$\beta$} \rput(7.9,4.7){$\alpha$}
\rput(11,.2){$\beta^{-1}$} \rput(11.9,1){$\alpha^{-1}$}
\rput(7.9,3.4){$\bullet$}  \rput(8.6,4){$\bullet$} 
\rput(11,1.6){$\bullet$}  \rput(10.4,1){$\bullet$} 
\end{pspicture}
\end{equation*}
There is a similar homotopy $\widetilde{Pinch}(-,-): [0,1]\times \Sigma^g\to \Sigma^0\vee \Sigma^g$ with  $\widetilde{Pinch}(0,-)=Pinch(0,-)$ and $\widetilde{Pinch}(1,-)=\pinch{0,g}$, which is obtained by taking a parametrization similar to~\eqref{eq:sigmagt} but with the small center square above and on the left of the big center square, \emph{i.e}, a parametrization ``symmetric" to~\eqref{eq:sigmagt} with respect to the anti-diagonal. The composition of the homotopies $Pinch(-,-)$ and $\widetilde{Pinch}(-,-)$ yields the desired homotopy between $\pinch{0,g}$ and $\pinch{g,0}$.
\end{proof}

\begin{remk}
Note that the surface product is \emph{not} (graded) commutative in general. For instance if $M=S^{3}$, the 3-dimensional sphere, then the center of $\big(\bigoplus_{g\geq 0}\mathbb{H}_\bullet(\Map(\Sigma^g,S^3)),\scup \big)$ is $\mathbb{H}_\bullet(\Map(\Sigma^0,S^3))$, see Example~\ref{E:oddspheres}.
\end{remk}

For any genus $g>0$-surface, there is a map $\pi_g:\Sigma^g\to \Sigma^0\cong S^2$ obtained by collapsing all edges $a_1,b_1,\dots$ of the $4g$-gon to a point. By pullback it yields a map $\pi^g:\Map(S^2,M)\to \Map(\Sigma^g,M)$. Hence, a linear morphism $\pi^g_*:H_\bullet(\Map(S^2,M))\to H_\bullet(\Map(\Sigma^g,M))$.
\begin{prop}\label{P:pig}
 Let $M$ be a compact oriented manifold. Then, for  $g> 0$ and $h> 0$, \begin{itemize}
\item[i)] the map $\pi^g_*:\mathbb{H}_\bullet(\Map(S^2,M))\to \mathbb{H}_\bullet(\Map(\Sigma^g,M))$ is an $\mathbb{H}_\bullet(\Map(S^2,M))$-module morphism and satisfies 
$$\pi^g_*(x)\scup \pi^h_*(y)=\pi^{g+h}_*(x\scup y); $$ \item[ii)]   $\pi^g_*(x)=x\scup [M]_g$ for any $x\in \mathbb{H}_\bullet(\Map(S^2,M))$.
\end{itemize}\end{prop} 
\begin{proof}
 From the definitions of $\pinch{g,h}$ for $g,h\geqslant 0$ (See Lemma~\ref{L:modelglueing} and the arrows~\eqref{E:sigma0hcollapse} and~\eqref{E:sigmag0collapse} we easily get that the two maps
$ (\pi_g\vee \pi_h)\circ \pinch{g,h}$ and $\pinch{0,0}\circ (\pi_{g+h})$ are homotopic and further that three maps $ (\pi_g\vee {id})\circ \pinch{g,0}$, $(id\vee \pi_g)\circ \pinch{0,g}$ and $\pinch{0,0}\circ (\pi_{g})$ are homotopic to each other.
Now the  claim~{\bf i)} follows from the naturality of  Umkehr maps: $$((\pi^g\vee \pi^h)')_*\circ (\rho_{in})_!=(\rho_{in})_!\circ (\pi^g\times \pi^h)_* $$ where $(\pi^g\vee \pi^h)': \Map(\Sigma^0\vee \Sigma^0, M) \to  \Map(\Sigma^g\vee \Sigma^h, M)$ is the natural map induced by $\pi_g\vee \pi_h:\Sigma^g\vee \Sigma^h\to \Sigma^0\vee \Sigma^0$.

\smallskip
 
By~{\bf i)} and Proposition~\ref{P:bimodule}, it is sufficient to prove claim {\bf ii)} for $x=[M]_0$ the unit of $\mathbb{H}_\bullet (\Map(S^2,M))$. That is to prove that $\pi^g_*([M]_0)=[M]_g$ which follows since $\pi^g \circ i_0= i_g$ where $i_g:M\hookrightarrow \Map(\Sigma^g,M)$ is the canonical embedding of $M$ as constant maps.
\end{proof}
\begin{ex}
By Proposition~\ref{P:pig} applied to the unit $[M]_0$, we get, for any $g,h>0$  $$[M]_g\scup [M]_{h} = \pi^g_*([M]_0)\scup \pi^h_*([M]_0) =\pi^{g+h}_{*}([M]_0)= [M]_{g+h}.$$In particular, $[M]_g$ and $[M]_h$ commute. 
\end{ex}

\subsection{Surface Hochschild cup product} \label{Surface Hochschild cup product}
In this section we give an analogue of the surface product defined in higher Hochschild cohomology over surfaces. Similarly to the  Hochschild homology over a simplicial set, there are Hochschild cochain complexes associated to any pointed simplicial set $Y_\bullet$ (see~\cite{G}) defined as follows. Let $(A,d)$ be a differential graded commutative algebra and $(M,d)$ an $A$-module viewed as a symmetric bimodule. We define 
 $$ CH^n_{Y_\bullet}(A,M):=Hom_{k}\big(\bigoplus_{k\geq 0}  A^{\otimes y_k}, M \big)^{n-k} $$ where the upper index $n-k$ is the total degree of a map $A^{\otimes y_k}\to A$. A map of pointed sets $\gamma:Y_k\to Y_l$ and a linear map $f:A^{\otimes y_l}\to M$, yields a map $\gamma^* f:A^{\otimes y_k}\to M$ given, for $a_1\otimes \dots\otimes a_{y_k}\in A^{\otimes y_k}$, by $$\gamma^* f(a_1\otimes \dots\otimes a_{y_k})= \pm b_0\cdot f(b_1\otimes \cdots\otimes b_{y_l})$$ where  $b_{j\geq 1}=\prod_{i\in f^{-1}(j)} a_i$ and $b_0=\prod_{0\neq i\in f^{-1}(0)} a_i$. The sign $\pm$ is the total Koszul sign obtained as the sum of  $(-1)^{|x|\cdot |y|}$ whenever $y$ moves across $x$ as in Definition~\ref{D:Hoch}.  Note that $CH^\bullet_{Y_\bullet}(A,M)$ is thus a cosimplicial vector space, with cosimplicial structure induced by the boundaries $d_i$ and degeneracies $s_j$ of $Y_\bullet$. The differential on $CH_{Y_\bullet}^\bullet(A,M)$ is given, for $f: A^{\otimes y_k}\to M$, by the sum $D(f)=(-1)^k d_f+ b_f$, where $d_f:A^{\otimes y_k}\to M$ is given, for $a_1\otimes \dots\otimes a_{y_k}\in A^{\otimes y_k}$,  by 
$$d_f(a_1\otimes \dots\otimes a_{y_k})= d\big(f(a_1\otimes \dots\otimes a_{y_k})\big) + \sum_{i=1}^{y_k}\pm f(a_1\otimes \dots\otimes d(a_i)\otimes \dots\otimes a_{y_k})$$ and $b_f: A^{\otimes y_{k+1}}\to M$ is given, for $a_1\otimes \dots\otimes a_{y_{k+1}}\in A^{\otimes y_{k+1}}$,  by
$$b_f(a_1\otimes \dots\otimes a_{y_{k+1}})=\sum_{i=0}^{k+1}(-1)^i (d_i^* f)(a_1\otimes \dots\otimes a_{y_{k+1}}). $$ Again the $\pm$ sign is the Koszul sign as in Definition~\ref{D:Hoch}.

\smallskip

As for homology, the cosimplicial identities imply $D^2=0$. We call $CH^\bullet_{Y_\bullet}(A,M)$ the Hochschild cochain complex for $Y_\bullet$ of $A$ with value in $M$. We denote $HH^n_{Y_\bullet}(A,M)$ its cohomology groups. Let $X_\bullet \stackrel{f}\to Y_\bullet$ be a morphism of pointed simplicial sets. Then, for any $k$, we have a map $f_k:X_k\to Y_k$, thus a map $f^*_k:Hom(A^{\otimes y_k},M)\to Hom(A^{\otimes x_k},M)$. Since, $f$ is simplicial, the map $f_k^*$ combines to give a  cochain complex morphism $f^*: \big(CH^\bullet_{Y_\bullet}(A,M),D\big)\to  \big(CH^\bullet_{X_\bullet}(A,M),D\big)$. The following Lemma follows from~\cite{P, G}:
\begin{lem}\label{L:quis}
The higher Hochschild cochain complex $CH^\bullet_{Y_\bullet}(A,M)$ is covariant in $M$, contravariant with respect to $A$ and $Y_\bullet$ and preserves homology equivalences, namely, if $f:A\to A'$, $g:M\to M'$ are quasi-isomorphisms and $\gamma:X_\bullet\to Y_\bullet$ induces an isomorphism in homology , then $f^*$, $g_*$ and $\gamma^*$ are all quasi-isomorphisms. 
\end{lem}

\smallskip \subsubsection{The cup product}

We now define a cup product  $$\cup: \HHS{g}(A,B)\otimes \HHS{h}(A,B) \to \HHS{g+h}(A,B)$$ for the Hochschild cohomology over surfaces, where  $B$ is a differential graded commutative and unital $A$-algebra, viewed as a symmetric bimodule. We are particularly interested in the case $B=A$. Henceforth, we use the simplicial model $\Sigma^g_\bullet$ of a surface $\Sigma^g$ of Section~\ref{S:surfacemodel}. 

\smallskip

We first consider the case $g,h\geq 1$. Since $\CHS{g}(A,B)$ is a cosimplicial complex, the tensor product $\CHS{g}(A,B)\otimes \CHS{h}(A,B)$ is bicosimplicial and we have the Alexander-Whitney quasi-isomorphisms
$$\CH{g}{i}(A,B)\otimes \CH{h}{j}(A,B) \stackrel{AW}\longrightarrow \CH{g}{i+j}(A,B)\otimes \CH{h}{i+j}(A,B)$$
where the right hand side is equipped with the diagonal cosimplicial structure. Recall that the Alexander-Whitney map is explicitly given by $AW=AW_{(1)}^*\otimes AW_{(2)}^*$ where $AW_{(1)}$ is the  map  $[i]\stackrel{\delta_{i+1}}\to[i+1]\stackrel{\delta_{i+2}}\to\dots  \stackrel{\delta_{i+j}}\to[i+j]$ (in the category $\Delta$ see Definition~\ref{D:Delta}) and  $AW_{(2)}$ is the  map $[j]\stackrel{\delta_{0}}\to  [1+j]\stackrel{\delta_0}\to \dots \stackrel{\delta_0}\to [i+j]$.

\smallskip

 Let $f\otimes g$ be in $\CH{g}{n}(A,B)\otimes \CH{h}{n}(A,B) $, then we define the ``wedge" $f\vee g \in \CHW{g}{h}{n}(A,B)$ of $f$ and $g$ by the formula
\begin{equation*}
f\vee g(a_1\otimes \cdots\otimes a_{\sigma^g_n} \otimes a_{\sigma^g_n+1} \otimes \cdots \otimes a_{\sigma^g_n+\sigma^h_n}) = f(a_1\otimes \cdots\otimes a_{\sigma^g_n}) \cdot g(a_{\sigma^g_n+1} \otimes \cdots \otimes a_{\sigma^g_n+\sigma^h_n})
\end{equation*} 
where $\cdot$ in the right hand side is the multiplication in the algebra $B$. 

\smallskip

In Section~\ref{S:surfacemodel} we defined  the pinching map $\pinch{g,h}: (\Sigma^g \vee \Sigma^h)_\bullet \to \Sigma^{g+h}_\bullet$~\eqref{eq:pinch}, which is a map of pointed simplicial sets. Composing the pinching map with the wedge and Alexander-Whitney maps, we make the following definition.
\begin{defn}\label{D:cupg} For $g,h\geq 1$, the cup-product  is the composition
\begin{multline*}
 \cup: \CH{g}{i}(A,B)\otimes \CH{h}{j}(A,B) \\ \stackrel{AW} 
\to  \CH{g}{i+j}(A,B)\otimes \CH{h}{i+j}(A,B) \\
\stackrel{\vee} \to\CHW{g}{h}{i+j}(A,B)  \stackrel{\pinch{g,h}^*}\to\CH{g+h}{i+j}(A,B)
\end{multline*}
\end{defn}
\begin{prop}\label{P:cup} Let $B$ be a (differential graded) commutative $A$-algebra. The cup product $\cup:\CHS{g}(A,B)\otimes \CHS{h}(A,B) \to \CHS{g+h}(A,B)$ is a map of cochain complexes and is associative. \end{prop}
\begin{proof}
It is straightforward to check that the wedge map $(f,g)\mapsto f\vee g$ is a morphism of simplicial modules. Since $\pinch{g,h}$ is a simplicial morphism and $AW$ a map of chain complexes, $\cup$ is a map of cochain complexes. Now the result follows from Lemma~\ref{L:modelglueing} and the associativity of $B$.
\end{proof}

We now turn to the genus zero case. Similarly to Section~\ref{thom}, there is a $\HHS{0}(A,B)$-module structure on $\bigoplus_{g\geq 1} \HHS{g}(A,B)$.
 However the module structure is a little bit more subtle since the standard model for the sphere is slightly different from the other genus models (we still assume that $B$ is an $A$-algebra). 
Note that $\CHS{0}(A,B)\cong CH_{S^2_\bullet}^\bullet(A,B)$ with simplicial structure described in Example~\ref{E:2sphere}. Thus, $CH_{S^2_k}^\bullet(A,B)=\{f:A^{\otimes k^2}\to B\}$.
Let $f\in CH_{S^2_p}^\bullet (A,B)$ and $g\in CH_{S^2_q}^\bullet(A,B)$. Then we define $f\cup g\in CH_{S^2_{p+q}}^\bullet(A,B)$  by the formula 
\begin{equation}\label{eq:cup00}
f\cup g\big((x_{i,j})\big) = \pm f\big((x_{i,j})_{i,j\leq p}\big)\cdot g\big((x_{i,j})_{p+1\leq i,j}\big) \cdot\hspace{-0.4cm} \prod_{\scriptsize \begin{array}{ll} i\geq p+1\\ \scriptsize  j\leq p \end{array}}\hspace{-0.4cm} x_{i,j} \cdot \hspace{-0.4cm}\prod_{\scriptsize \begin{array}{ll} i\leq p\\ \scriptsize  j\geq p+1\end{array}}\hspace{-0.4cm} x_{i,j}
\end{equation}
where $(x_{ij})$ stands for a tensor $x_{1,1}\otimes\cdots \otimes x_{p+q, p+q}$. It is straightforward to check that $(CH_{S^2_\bullet}^\bullet(A,B),\cup, D)$ is a differential graded associative algebra. In fact,
\begin{prop}[\cite{G} Proposition 3.2 and Remark 1] \label{P:cup00} The cup-product makes $HH_{S^2}^\bullet(A,B)$ a graded commutative algebra.
\end{prop}
We now define the cup-product $\cup:\CHS{0}(A,B)\otimes \CHS{g}(A,B)\to \CHS{g}(A,B)$. Later on, using the edgewise subdivision we will give another model for the cup-product in Section~\ref{S:subdivision} (see Definition~\ref{D:edgewisecup}) which will alow us to define equivalent cup-products for Hochschild cohomology over different simplicial models for the surfaces. 
\begin{defn}\label{D:cup0} 
Let $f\in\CH{0}{k}(A,B)$ and $g\in\CH{g}{l}(A,B)$, {\it i.e.} $f: A^{\otimes k^2}\to B$ and $g: A^{\otimes \sigma^g_{l}}\to B$, where $\sigma^g_{l}=\# \Sigma^g_{l}-1$. We will define  $f\cup g\in \CH{g}{k+l}(A,B)$ and $g\cup f\in \CH{g}{k+l}(A,B)$. The idea is to use the Alexander-Whitney diagonal in a slightly different way. Applying  $AW_{(2)}:[l]\to [k+l]$ from above induces a map $\Sigma^g_{k+l}\to \Sigma^g_l$ for the simplicial model described in \eqref{E:sigmag}, which is given by collapsing certain elements in $\Sigma^g_{k+l}$, 
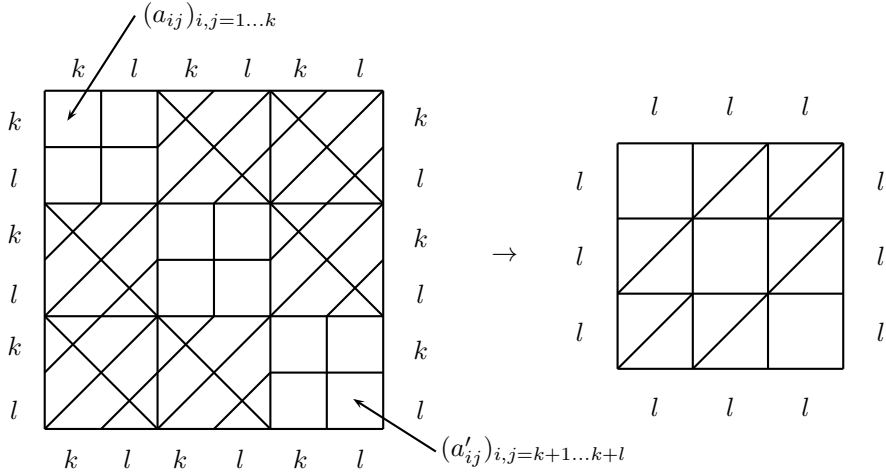
\[ \begin{pspicture}(0,0)(5.5,6.2)
 \psline(0.5,.5)(0.5,5) \psline(.5,.5)(5,.5) \psline(2,.5)(2,5) \psline(3.5,.5)(3.5,5)
\psline(5,.5)(5,5) \psline(.5,5)(5,5) \psline(0.5,2)(5,2) \psline(.5,3.5)(5,3.5)
 \psline(.5,4.25)(2,4.25) \psline(2,2.75)(3.5,2.75) \psline(3.5,1.25)(5,1.25) 
 \psline(1.25,3.5)(1.25,5) \psline(2.75,2)(2.75,3.5) \psline(4.25,.5)(4.25,2) 
\psline(.5,2)(2,.5) \psline(0.5,3.5)(3.5,.5) \psline(2,5)(5,2) \psline(3.5,5)(5,3.5)
\psline(.5,2)(3.5,5) \psline(0.5,.5)(2,2) \psline(3.5,3.5)(5,5) \psline(2,.5)(5,3.5)
\psline(.5,2.75)(1.25,3.5) \psline(2,4.25)(2.75,5) 
\psline(.5,1.25)(2,2.75) \psline(2.75,3.5)(4.25,5)
\psline(1.25,.5)(2.75,2) \psline(3.5,2.75)(5,4.25) 
\psline(2.75,.5)(3.5,1.25) \psline(4.25,2)(5,2.75)

 \rput(0.1,4.6){$k$} \rput(0.1,3.1){$k$} \rput(0.1,1.6){$k$}  
\rput(0.85,0.1){$k$} \rput(2.3,0.1){$k$}  \rput(3.9,0.1){$k$} 
\rput(5.5, 1.55){$k$}  \rput(5.5,3.05){$k$} \rput(5.5,4.65){$k$} 
 \rput(3.9,5.3){$k$} \rput(2.45,5.3){$k$}   \rput(0.95,5.3){$k$}

\rput(0.1,3.85){$l$}  \rput(0.1,2.3){$l$} \rput(0.1,.75){$l$}  
\rput(1.6,0.1){$l$} \rput(3.1,0.1){$l$}  \rput(4.7,0.1){$l$} 
\rput(5.5, .75){$l$}  \rput(5.5,2.3){$l$} \rput(5.5,3.85){$l$} 
 \rput(4.7,5.3){$l$} \rput(3.2,5.3){$l$}  \rput(1.7,5.3){$l$}

\psline{->}(1.7,6)(.8,4.6) \rput(2.7,6){$(a_{ij})_{i,j=1\dots k}$}
\psline{->}(5.7,.2)(4.6,.9) \rput(7,.2){$(a'_{ij})_{i,j=k+1\dots k+l}$}
\end{pspicture}
\
\begin{pspicture}(-2,-.8)(4,4)
 \psline(0.5,.5)(0.5,3.5) \psline(.5,.5)(3.5,.5)
 \psline(1.5,.5)(1.5,3.5) \psline(.5,1.5)(3.5,1.5)
 \psline(2.5,.5)(2.5,3.5) \psline(.5,2.5)(3.5,2.5)
 \psline(3.5,.5)(3.5,3.5) \psline(.5,3.5)(3.5,3.5)
 \psline(.5,.5)(1.5,1.5) \psline(2.5,2.5)(3.5,3.5)
 \psline(.5,1.5)(2.5,3.5) \psline(1.5,.5)(3.5,2.5)
 \rput(-1,2){$\to$}
 \rput(0,3){$l$}  \rput(0,2){$l$} \rput(0,1){$l$}  \rput(1,0){$l$}
 \rput(2,0){$l$}  \rput(3,0){$l$} \rput(4,1){$l$}  \rput(4,2){$l$}
 \rput(4,3){$l$}  \rput(3,4){$l$} \rput(2,4){$l$}  \rput(1,4){$l$}
\end{pspicture}
 \]
In particular, all the elements of coordinates $(i,j)_{i,j=1\dots k}$ in the top left square of $\Sigma_{k+l}^g$ are collapsed to the basepoint in $\Sigma^g_l$.

Let $a$ be a homogeneous element in $A^{\otimes \sigma^g_{k+l}}$, and denote by $a_{0j}$, $a_{i0}$, $a_{i,j}$ ($i,j=1\dots k+l$) the $(k+l)^2 +2(k+l)$ tensor factors of $a$ corresponding to the top left square in the simplicial model $\Sigma^g_{k+l}$. Then, in particular, the elements $(a_{ij})_{i,j=1\dots k}$ get multiplied to the basepoint under the induced $AW_{(2)}$, but there are also other elements, whose product we denote by $\prod c$. Then, we may express $AW_{(2)}^*(g)(a)$ as $$AW_{(2)}^*(g)(a)= \pm \big(\prod_{i,j=1\dots k}(a_{ij}) \big)\cdot g(b)\cdot \prod  c,$$ where $b, c$ are certain (products of) subtensors of $a\in A^{\otimes \sigma^g_{k+l}}$, determined by the mapping $AW_{(2)}$. With this notation, we define $f\cup g\in \CH{g}{k+l}(A,B)$ by
$$f\cup g\, (a)=\pm f((a_{ij})_{i,j=1\dots k})\cdot g(b)\cdot \prod c.$$
Similarly,  $AW_{(1)}:[k]\to [k+l]$ induces a map $\Sigma^g_{k+l}\to \Sigma^g_k$, which on the Hochschild cochain level may be expressed as
$$AW_{(1)}^*(g)(a)=\pm  g(b')\cdot \big(\prod_{i,j=k+1\dots k+l}(a'_{ij})\big)\cdot \prod  c',$$ where $a_{ i,j}'$ ($i,j=1\dots k+l$) are the tensor factors of $a$ corresponding to the lower right square of $\Sigma^g_{k+l}$, and $b',c'$ are determined by $AW_{(1)}$ similarly to the above. Define $g\cup f\in \CH{g}{k+l}(A,B)$ by
$$ g \cup f \, (a) =\pm g(b')\cdot f((a'_{ij})_{i,j=k+1\dots k+l}) \cdot\prod c'.$$
\end{defn}
\begin{ex} Assume the genus is $1$, and  $g\in \CH{1}{l}(A,B)$. We denote by $(a_{i,j})_{\scriptsize \begin{array}{c} i,j=0\dots l \\ ij\neq 0\end{array}}\in CH_{\bullet}^{\Sigma^1_{k+l}}(A,B)$ a generic element, {\it i.e.} $a_{0,0}\otimes\dots \otimes a_{l,l}\in A^{\otimes (l^2+2l +1)}$.  Then, for any $f \in \CH{0}{k}(A,B)$, one has
$$f \cup g (a_{i,j}) = \pm f\big((a_{i,j})_{i,j=1\dots k}\big) \cdot g\big((b_{i,j})_{\scriptsize \begin{array}{ll} i,j=k\dots k+l\\ ij\neq k^2\end{array}}\big)\cdot\prod_{i=1}^{k}c_{0i}\cdot c_{i,0}$$ where $b_{k,j}=a_{0,j}\cdots a_{k,j}$, $b_{i,k}=a_{i,0}\cdots a_{i,k}$ and $b_{i,j}=a_{i,j}$ for $i,j>k$, and $c_{0,i}=a_{0,i}, c_{i,0}=a_{i,0}$.
\end{ex}
\begin{remk}
For $g=0$, Definition~\ref{D:cup0} coincides with formula~\eqref{eq:cup00}. 
\end{remk}
Definition~\ref{D:cup0}  induces a right and a left action of $\CHS{0}(A,B)$ on $\CHS{g}(A,B)$. 
\begin{lem}\label{L:cupbimodule}
The cup product makes $\CHS{g}(A,B)$ a differential graded  $\CHS{0}(A,B)$-bimodule.
\end{lem}
\begin{proof}
It follows from the associativity of $A$, the fact that $B$ is an $A$-algebra and formula~\eqref{eq:cup00} that the cup-product makes $\CHS{0}(A,B)$ a unital associative algebra with unit $1_B\in B\cong CH^{0}_{\Sigma^0_0}(A,B)$ which is bigraded with respect to both simplicial degree and internal degree (of $A$ and $B$). Further, for $f\in CH^\bullet_{S^2_p}(A,B)$ and $g\in CH_{S^2_q}^\bullet(A,B)$, note that 
\begin{eqnarray*}\label{eq:Lbim}
d_{p+1}^*(f)\cup g\big((a_{i,j})\big) &= &\pm f\big((a_{i,j})_{i,j\leq p}\big)\cdot\hspace{-0.2cm} \prod_{\scriptsize \begin{array}{ll}\,\,\, i,j\leq p+1\\  i \mbox{ or } j=p+1 \end{array}}  \hspace{-0.2cm} a_{i,j} \\ &&\quad \cdot  g\big((a_{i,j})_{p+2\leq i,j}\big) \cdot\hspace{-0.4cm} \prod_{\scriptsize \begin{array}{ll} i\geq p+2\\ \scriptsize  j\leq p+1 \end{array}}\hspace{-0.4cm} a_{i,j} \cdot \hspace{-0.4cm}\prod_{\scriptsize \begin{array}{ll} i\leq p+1\\ \scriptsize  j\geq p+2\end{array}}\hspace{-0.4cm} a_{i,j}\\
&= &- f\cup d_{0}^*(g) \big((a_{i,j})\big).
\end{eqnarray*}
since  the total degree $|d_{p+1}^*(f)|=1+|f|$. Hence
$$b(f)\cup g + (-1)^{|f|} f\cup b(g) = \sum_{i=0}^{p} (d_i^* f) +d_{p+1}^*(f)\cup g + f\cup d_0^*(g) +\sum_{i=1}^{q+1} f\cup (d_i^* g)=b(f\cup g).$$ It follows that $\CHS{0}(A,B)$ is a unital differential graded associative algebra. Similarly, one proves $D(f\cup g)=D(f) \cup g \pm f\cup D(g)$. 

Now, for any $f\in CH^{\bullet}_{\Sigma^0_p}(A,B)$, $g\in CH^{\bullet}_{\Sigma^0_q}(A,B)$, and $h\in \CH{k}{r}(A,B)$, we may use Definition~\ref{D:cup0} and the fact that $AW_{(2)}:[r]\to [r+p+q]$ is equal to the composition $[r]\stackrel{AW_{(2)}}\to [r+q]\stackrel{AW_{(2)}}\to [r+p+q]$, to see that
\begin{eqnarray*} 
((f\cup g) \cup  h)(a) &=& \pm (f\cup g)((a_{ij})_{i,j=1\dots p+q})\cdot h(b)\cdot \prod c\\ 
&=& f((a_{ij})_{i,j=1\dots p})\cdot g((a_{ij})_{i,j=p+1\dots p+q}) \cdot \\ 
&& \hspace{-0.4cm} \prod_{\scriptsize \begin{array}{ll} \quad \,\,\,\, 1\leq i\leq p\\ \scriptsize  p+1\leq j\leq p+q \end{array}}\hspace{-0.4cm} a_{i,j} \cdot \prod_{\scriptsize \begin{array}{ll} p+1\leq i\leq p+q\\ \scriptsize \quad\,\,\,\, 1\leq j\leq p\end{array}}\hspace{-0.4cm} a_{i,j} \cdot h(b)\cdot \prod c\\
&=& (f\cup (g\cup h))(a)
\end{eqnarray*}
This is exactly the left module identity. Similarly, the right module identity is obtained by using the equality of $AW_{(1)}:[r]\to [r+p+q]$ with the composition $[r]\stackrel{AW_{(1)}}\to [r+p]\stackrel{AW_{(1)}}\to [r+p+q]$, whereas the compatibility of left and right module structure is obtained via the equality of $[r]\stackrel{AW_{(1)}}\to [r+p]\stackrel{AW_{(2)}}\to [r+p+q]$ and $[r]\stackrel{AW_{(2)}}\to [r+q]\stackrel{AW_{(1)}}\to [r+p+q]$.
\end{proof}
 This bimodule structure is not symmetric at the chain level but as we will discuss it will induce a symmetric bimodule structure after passing to homology.

\smallskip \subsubsection{Subdivision} \label{S:subdivision}

We now give another description of the bimodule structure of $\CHS{g}(A,B)$, by means of the edgewise subdivision. Recall the notations of Definition~\ref{D:Delta}. The edgewise subdivision~\cite{BHM, McC} is an endofunctor of the simplicial category $\Delta$ which associates, to any simplicial set $X_\bullet$, a simplicial set $sd_2(X_\bullet)$ whose realization is homeomorphic to the one of $X_\bullet$. One of its main properties is that the realization of the edgewise subdivision $|sd_2(\Delta^n_\bullet)|$  of the standard $n$-simplex $\Delta^n$  is a triangulation of $|\Delta^n_\bullet|$ by $2^n$ standard simplexes. The functor $sd_2:\Delta\to \Delta$ is defined by $sd_2([n-1])=[2n-1]$, and, for any map, $f: [n-1]\to [m-1]$, by $sd_2(f):[2n-1]\to [2m-1], sd_2(f):i+ n j\mapsto f(i) +mj$ where $0\leqslant i\leqslant n-1$ and $j\in \{0,1\}$, see~\cite{BHM}. The edgewise subdivision $sd_2 (X_\bullet)$ of a simplicial set $ X_\bullet$ is the composition $X_\bullet \circ sd_2$.

There is a natural homeomorphism $D:|sd_2(X_\bullet)|\stackrel{\sim}\to |X_\bullet|$(see~\cite[Lemma 1.1]{BHM}) induced by the maps $\Delta^{n-1}\times X_{2n-1} \to \Delta^{2n-1}\times X_{2n-1}$ defined by $(u,x)\mapsto ((u/2\oplus u/2),x)$ where $u=(u_0,\dots, u_{n-1})\in \R^{n}$ is such that $u_0+\dots +u_{n-1}=1$. In~\cite[Definition 3.3]{McC}  a natural chain map  $\mathcal{D}_\bullet(2): C(X_\bullet) \to C(sd_2(X_\bullet))$ is defined, where $C(Y_\bullet)$ is the chain complex associated to a simplicial set $Y_\bullet$.  More precisely, for any $x\in X_n$,
\begin{equation} \label{E:D(2)} \mathcal{D}_\bullet(2)(x) = \sum_{(\sigma,\eta) \in \mathcal{S}(2,n)} (-1)^{\sigma} X_\bullet(\varepsilon_{(\sigma,\eta)}) (x) \end{equation}  where $\mathcal{S}(2,n)$ is the set 
$$\mathcal{S}(2,n)=\{(\sigma,\eta) \in S_n\times Hom_{\Delta}([n-1],[1])\,|\, \sigma(i) >\sigma(i+1) \Rightarrow \eta(i-1)<\eta(i)\} $$
and $\varepsilon_{(\sigma,\eta)}:[2n+1] \to [n]$ is defined by \begin{equation} \label{E:epsilonsigmadelta} \varepsilon_{(\sigma,\eta)}^{-1}(\{j\})=\{\eta(j-1)\cdot(n+1)+\sigma(j),\dots,\eta(j)\cdot (n+1)+\sigma(j+1)-1\}.\end{equation}
McCarthy~\cite[Proposition 3.4 and Corollary 3.7]{McC} proved that $\mathcal{D}_\bullet(2)$ is a quasi-isomorphism realizing $D^{-1}$ in homology and passes to normalized chain complexes. 

\smallskip

The following Lemma is straightforward.
\begin{lem}\label{L:edgewisewedge} With the same notation as in Definition~\ref{D:wedge}, one has a natural isomorphism
$sd_2(X_\bullet\cup_{Z_\bullet} Y_\bullet)\cong sd_2(X_\bullet)\cup_{sd_2(Z_\bullet)} sd_2(Y_\bullet)$. 
\end{lem}

\begin{ex} \label{E:edgewise}
 Recall from Examples~\ref{S^1} and~\ref{E:wedge}  the pointed simplicial sets $S^1_\bullet$, $pt_\bullet$ and $I_\bullet$ for the circle, point and interval. Then $sd_2(S^1_n)=\{0,\dots,2n+1\}$, $sd_2(pt_n) =\{0\}$ and $sd_2(I_n)=\{0,\dots,2n+2\}$ and it is easy to see that $CH_\bullet^{sd_2(I_\bullet)}(A,M)= B(M,A,A)\otimes_{A} B(A,A,A)$ where $B(M,A,N)$ is the two sided bar construction and the tensor product uses the right (resp. left) $A$-module structure on $B(M,A,N)$ (resp. $B(A,A,A)$). In fact, $sd_2(I_\bullet)= I_\bullet \cup_{pt_\bullet} I_\bullet$. Further  $sd_2(S^1_\bullet)=sd_2(I_\bullet) \cup_{sd_2(pt_\bullet)} sd_2(pt_\bullet)$ where the two endpoints $0$ and $2n+2$ of $sd_2(I_n)=\{0,2n+2\}$ get collapsed.  In particular the Hochschild chain complex is $CH^{sd_2(S^1_n)}_\bullet(A,M) = M\otimes A^{\otimes n} \otimes A\otimes A^{\otimes n}$ with differential 
\begin{multline*}
D(a_0 \otimes\dots\otimes a_n\otimes a_{n+1} \otimes a_{n+2}\dots\otimes a_{2n+1}) = \sum_{i=0}^{2n+1} \pm a_0\otimes\dots\otimes d(a_i)\otimes \dots\otimes a_{2n+1}\\
+\sum_{i=0}^{n-1} \big(\pm a_0\otimes\dots\otimes(a_i\cdot a_{i+1})\otimes\dots\otimes a_{2n+1} \pm a_0\otimes\dots\otimes(a_{n+1+i}\cdot a_{n+i+2})\otimes\dots\otimes a_{2n+1}\big) \\
\pm (a_{2n+1}\cdot a_0)\otimes a_1\otimes\dots\otimes a_{2n} \pm a_0\otimes \dots\otimes (a_{n}\cdot a_{n+1})\otimes a_{n+2}\otimes \dots \otimes a_{2n+1}.
\end{multline*}
Similarly, the edgewise subdivision $sd_2(I^2_\bullet)$ of a square is canonically identified with the wedge 
\begin{equation*} 
\begin{pspicture}(0.5,0.5)(2.5,2.5)
 \psline(0.5,.5)(0.5,2.5) \psline(.5,.5)(2.5,.5)
 \psline(2.5,.5)(2.5,2.5) \psline(.5,2.5)(2.5,2.5)
 \psline(1.5,.5)(1.5,2.5) \psline(.5,1.5)(2.5,1.5)
\end{pspicture}
\end{equation*}
 of four standard squares $I^2_\bullet$.
\end{ex}

Using Section~\ref{S:surfacemodel} and Lemma~\ref{L:edgewisewedge}, we obtain that $sd_2(\Sigma^g_\bullet)$ is a wedge of $4g^2$ squares and $4g(g-1)$ triangles (where a model for a triangle is given by a square with an edge collapsed to a point). For instance, for a surface of genus 3, we obtain the following model
\begin{equation} \label{E:edgewisesigma3}
\begin{pspicture}(0,0)(5.5,5.5)
 \psline(0.5,.5)(0.5,5) \psline(.5,.5)(5,.5) \psline(2,.5)(2,5) \psline(3.5,.5)(3.5,5)
\psline(5,.5)(5,5) \psline(.5,5)(5,5) \psline(0.5,2)(5,2) \psline(.5,3.5)(5,3.5)
 \psline(.5,4.25)(2,4.25) \psline(2,2.75)(3.5,2.75) \psline(3.5,1.25)(5,1.25) 
 \psline(1.25,3.5)(1.25,5) \psline(2.75,2)(2.75,3.5) \psline(4.25,.5)(4.25,2) 
\psline(.5,2)(2,.5) \psline(0.5,3.5)(3.5,.5) \psline(2,5)(5,2) \psline(3.5,5)(5,3.5)
\psline(.5,2)(3.5,5) \psline(0.5,.5)(2,2) \psline(3.5,3.5)(5,5) \psline(2,.5)(5,3.5)
\psline(.5,2.75)(1.25,3.5) \psline(2,4.25)(2.75,5) 
\psline(.5,1.25)(2,2.75) \psline(2.75,3.5)(4.25,5)
\psline(1.25,.5)(2.75,2) \psline(3.5,2.75)(5,4.25) 
\psline(2.75,.5)(3.5,1.25) \psline(4.25,2)(5,2.75)

 \rput(0.1,4.6){$a'_1$} \rput(0.1,3.1){$b'_1$} \rput(0.1,1.6){$a'_2$}  
\rput(0.85,0.1){$b'_2$} \rput(2.3,0.1){$a'_3$}  \rput(3.9,0.1){$b'_3$} 
\rput(5.5, 1.55){${a'_3}^{-1}$}  \rput(5.5,3.05){${b'_3}^{-1}$} \rput(5.5,4.65){${a'_2}^{-1}$} 
 \rput(3.9,5.3){${b_2'}^{-1}$} \rput(2.45,5.3){${a_1'}^{-1}$}   \rput(0.95,5.3){${b_1'}^{-1}$}

\rput(0.1,3.85){$a_1$}  \rput(0.1,2.3){$b_1$} \rput(0.1,.75){$a_2$}  
\rput(1.6,0.1){$b_2$} \rput(3.1,0.1){$a_3$}  \rput(4.7,0.1){$b_3$} 
\rput(5.5, .75){$a_3^{-1}$}  \rput(5.5,2.3){$b_3^{-1}$} \rput(5.5,3.85){$a_2^{-1}$} 
 \rput(4.7,5.3){$b_2^{-1}$} \rput(3.2,5.3){$a_1^{-1}$}  \rput(1.7,5.3){$b_1^{-1}$}
\end{pspicture}
\end{equation}

\smallskip \subsubsection{Cup product via subdivision}

The reason for introducing the edgewise subdivison is that,
for any positive genus surface $\Sigma^g$, there is a   pinching map $P_{0,g}:sd_2(\Sigma^g_\bullet) \to \Sigma^0_\bullet\vee \Sigma^g_\bullet$, which is a simplicial.  The map $P_{0,3}:sd_2(\Sigma^3_\bullet) \to \Sigma^0_\bullet\vee \Sigma^3_\bullet$ is given by the following picture:
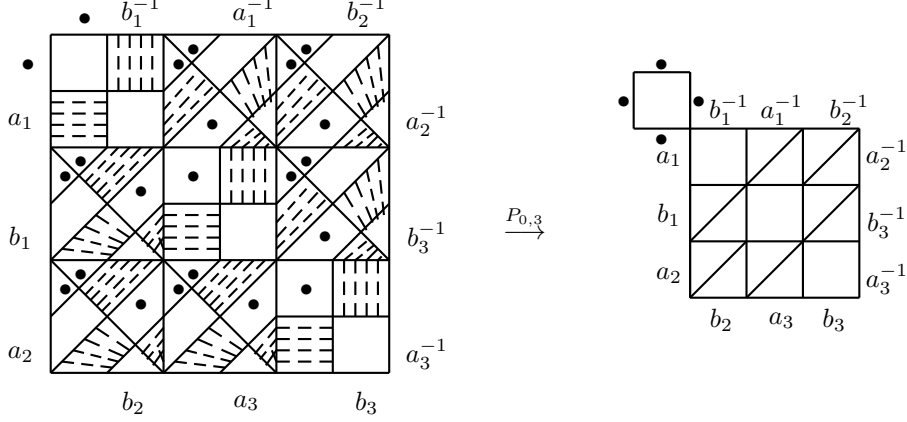
\begin{equation} \label{E:P03}
\begin{pspicture}(0,0)(11.5,5.5)
 \psline(0.5,.5)(0.5,5) \psline(.5,.5)(5,.5) \psline(2,.5)(2,5) \psline(3.5,.5)(3.5,5)
\psline(5,.5)(5,5) \psline(.5,5)(5,5) \psline(0.5,2)(5,2) \psline(.5,3.5)(5,3.5)
 \psline(.5,4.25)(2,4.25) \psline(2,2.75)(3.5,2.75) \psline(3.5,1.25)(5,1.25) 
 \psline(1.25,3.5)(1.25,5) \psline(2.75,2)(2.75,3.5) \psline(4.25,.5)(4.25,2) 
\psline(.5,2)(2,.5) \psline(0.5,3.5)(3.5,.5) \psline(2,5)(5,2) \psline(3.5,5)(5,3.5)
\psline(.5,2)(3.5,5) \psline(0.5,.5)(2,2) \psline(3.5,3.5)(5,5) \psline(2,.5)(5,3.5)
\psline(.5,2.75)(1.25,3.5) \psline(2,4.25)(2.75,5) 
\psline(.5,1.25)(2,2.75) \psline(2.75,3.5)(4.25,5)
\psline(1.25,.5)(2.75,2) \psline(3.5,2.75)(5,4.25) 
\psline(2.75,.5)(3.5,1.25) \psline(4.25,2)(5,2.75)

 \rput(0.2,4.6){$\bullet$}    \rput(0.95,5.2){$\bullet$}

\rput(0.1,3.85){$a_1$}  \rput(0.1,2.3){$b_1$} \rput(0.1,.75){$a_2$}  
\rput(1.6,0.1){$b_2$} \rput(3.1,0.1){$a_3$}  \rput(4.7,0.1){$b_3$} 
\rput(5.5, .75){$a_3^{-1}$}  \rput(5.5,2.3){$b_3^{-1}$} \rput(5.5,3.85){$a_2^{-1}$} 
 \rput(4.7,5.3){$b_2^{-1}$} \rput(3.2,5.3){$a_1^{-1}$}  \rput(1.7,5.3){$b_1^{-1}$}

\rput(2.2, 4.6){$\bullet$} \rput(2.4, 4.8){$\bullet$} \rput(2.65, 3.8){$\bullet$} 
\rput(3.7, 4.6){$\bullet$} \rput(3.9, 4.8){$\bullet$}\rput(4.15, 3.8){$\bullet$}
\rput(3.7, 3.1){$\bullet$} \rput(3.9, 3.3){$\bullet$} \rput(4.15, 2.3){$\bullet$}

\rput(.7, 3.1){$\bullet$} \rput(.9, 3.3){$\bullet$} \rput(1.7, 2.9){$\bullet$} 
\rput(.7, 1.6){$\bullet$} \rput(.9, 1.8){$\bullet$} \rput(1.7, 1.4){$\bullet$} 
\rput(2.2, 1.6){$\bullet$} \rput(2.4, 1.8){$\bullet$} \rput(3.2, 1.4){$\bullet$} 

\rput(2.4,3.1){$\bullet$} \rput(3.9,1.6){$\bullet$}

\psline[linestyle=dashed](.5,3.65)(1.25,3.65) \psline[linestyle=dashed](.5,3.8)(1.25,3.8) \psline[linestyle=dashed](.5,3.95)(1.25,3.95) \psline[linestyle=dashed](0.5,4.1)(1.25,4.1)
\psline[linestyle=dashed](1.4,4.25)(1.4,5)\psline[linestyle=dashed](1.55,4.25)(1.55,5)
\psline[linestyle=dashed](1.7,4.25)(1.7,5)\psline[linestyle=dashed](1.85,4.25)(1.85,5)

\psline[linestyle=dashed](2,2.15)(2.75,2.15) \psline[linestyle=dashed](2,2.3)(2.75,2.3) \psline[linestyle=dashed](2,2.45)(2.75,2.45) \psline[linestyle=dashed](2,2.6)(2.75,2.6)
\psline[linestyle=dashed](2.9,2.75)(2.9,3.5)\psline[linestyle=dashed](3.05,2.75)(3.05,3.5)
\psline[linestyle=dashed](3.2,2.75)(3.2,3.5)\psline[linestyle=dashed](3.35,2.75)(3.35,3.5)

\psline[linestyle=dashed](3.5,.65)(4.25,.65) \psline[linestyle=dashed](3.5,.8)(4.25,.8) \psline[linestyle=dashed](3.5,.95)(4.25,.95) \psline[linestyle=dashed](3.5,1.1)(4.25,1.1)
\psline[linestyle=dashed](4.4,1.25)(4.4,2)\psline[linestyle=dashed](4.55,1.25)(4.55,2)
\psline[linestyle=dashed](4.7,1.25)(4.7,2)\psline[linestyle=dashed](4.85,1.25)(4.85,2)

\psline[linestyle=dashed](2,3.65)(2.65,4.3) \psline[linestyle=dashed](2,3.8)(2.55,4.4) \psline[linestyle=dashed](2,4)(2.4,4.45)
\psline[linestyle=dashed](2.9,3.5)(3.2,3.75)\psline[linestyle=dashed](3.05,3.5)(3.3,3.7)
\psline[linestyle=dashed](3.2,3.5)(3.35,3.65)\psline[linestyle=dashed](3.35,3.5)(3.4,3.55)
\psline[linestyle=dashed](2.9,4.4)(3.2,3.9)\psline[linestyle=dashed](3.05,4.55)(3.28,4)
\psline[linestyle=dashed](3.2,4.7)(3.35,4.04)\psline[linestyle=dashed](3.35,4.85)(3.43,4.1)

\psline[linestyle=dashed](3.5,3.65)(4.15,4.3) \psline[linestyle=dashed](3.5,3.8)(4.05,4.4) \psline[linestyle=dashed](3.5,4)(3.9,4.45)
\psline[linestyle=dashed](4.4,3.5)(4.7,3.75)\psline[linestyle=dashed](4.55,3.5)(4.8,3.7)
\psline[linestyle=dashed](4.7,3.5)(4.85,3.65)\psline[linestyle=dashed](4.85,3.5)(4.9,3.55)
\psline[linestyle=dashed](4.4,4.4)(4.7,3.9)\psline[linestyle=dashed](4.55,4.55)(4.78,4)
\psline[linestyle=dashed](4.7,4.7)(4.85,4.04)\psline[linestyle=dashed](4.85,4.85)(4.93,4.1)

\psline[linestyle=dashed](3.5,2.15)(4.15,2.8) \psline[linestyle=dashed](3.5,2.3)(4.05,2.9) \psline[linestyle=dashed](3.5,2.5)(3.9,2.95)
\psline[linestyle=dashed](4.4,2)(4.7,2.25)\psline[linestyle=dashed](4.55,2)(4.8,2.2)
\psline[linestyle=dashed](4.7,2)(4.85,2.15)\psline[linestyle=dashed](4.85,2)(4.9,2.05)
\psline[linestyle=dashed](4.4,2.9)(4.7,2.4)\psline[linestyle=dashed](4.55,3.05)(4.78,2.5)
\psline[linestyle=dashed](4.7,3.2)(4.85,2.54)\psline[linestyle=dashed](4.85,3.35)(4.93,2.6)

\psline[linestyle=dashed](2,2.3)(1.86, 2.1) \psline[linestyle=dashed](2,2.45)(1.78, 2.18)
\psline[linestyle=dashed](2,2.6)(1.71,2.25) \psline[linestyle=dashed](2,2.15)(1.94, 2.3)
\psline[linestyle=dashed](1.4,3.5)(.96,3.04)\psline[linestyle=dashed](1.55,3.5) (1.03, 2.97)
\psline[linestyle=dashed](1.7,3.5)(1.11,2.89)\psline[linestyle=dashed](1.85,3.5)(1.18, 2.82)
\psline[linestyle=dashed](.65,2.15)(1.32,2.07)\psline[linestyle=dashed](.8, 2.3) (1.4, 2.15)
\psline[linestyle=dashed](.95,2.45)(1.47, 2.22)\psline[linestyle=dashed](1.1, 2.6)(1.55, 2.3)
  
\psline[linestyle=dashed](2,.8)(1.86, .6) \psline[linestyle=dashed](2,.95)(1.78, .68)
\psline[linestyle=dashed](2,1.1)(1.71,.75) \psline[linestyle=dashed](2,.65)(1.94, .8)
\psline[linestyle=dashed](1.4,2)(.96,1.54)\psline[linestyle=dashed](1.55,2) (1.03, 1.47)
\psline[linestyle=dashed](1.7,2)(1.11,1.39)\psline[linestyle=dashed](1.85,2)(1.18, 1.32)
\psline[linestyle=dashed](.65,.65)(1.32,.57)\psline[linestyle=dashed](.8, .8) (1.4, .65)
\psline[linestyle=dashed](.95,.95)(1.47, .72)\psline[linestyle=dashed](1.1, 1.1)(1.55, .8)
  
\psline[linestyle=dashed](3.5,.8)(3.36, .6) \psline[linestyle=dashed](3.5,.95)(3.28, .68)
\psline[linestyle=dashed](3.5,1.1)(3.21,.75) \psline[linestyle=dashed](3.5,.65)(3.44, .8)
\psline[linestyle=dashed](2.9,2)(2.46,1.54)\psline[linestyle=dashed](3.05,2) (2.53, 1.47)
\psline[linestyle=dashed](3.2,2)(2.61,1.39)\psline[linestyle=dashed](3.35,2)(2.68, 1.32)
\psline[linestyle=dashed](2.15,.65)(2.82,.57)\psline[linestyle=dashed](2.3, .8) (2.9, .65)
\psline[linestyle=dashed](2.45,.95)(2.97, .72)\psline[linestyle=dashed](2.6, 1.1)(3.05, .8)
\rput(6.8,2.5){$\stackrel{P_{0,3}}\longrightarrow$} 


\psline(8.25,4.5)(9,4.5) \psline(8.25,4.5)(8.25, 3.75) \psline(8.25, 3.75)(11.25,3.75) 
\psline(9,4.5)(9,1.5) \psline(9, 3)(11.25,3) \psline(9, 2.25)(11.25,2.25)
\psline(9, 1.5)(11.25,1.5) \psline(9,4.5)(9,1.5) 
\psline(9.75, 3.75)(9.75,1.5) \psline(10.5, 3.75)(10.5,1.5) \psline(11.25, 3.75)(11.25,1.5)
\psline(9,2.25)(10.5,3.75) \psline(9,1.5)(9.75,2.25) \psline(10.5,3)(11.25,3.75) \psline(9.75, 1.5)(11.25,3)\,\,
 
\rput(8.6,3.35){$a_1$}  \rput(8.6,2.6){$b_1$} \rput(8.6,1.8){$a_2$}  
\rput(9.3,1.2){$b_2$} \rput(10.1,1.2){$a_3$}  \rput(10.8,1.2){$b_3$} 
\rput(11.5, 1.75){$a_3^{-1}$}  \rput(11.5,2.5){$b_3^{-1}$} \rput(11.5,3.4){$a_2^{-1}$} 
 \rput(11,4){$b_2^{-1}$} \rput(10.1,4){$a_1^{-1}$}  \rput(9.4,4){$b_1^{-1}$}
\rput(8,4.1){$\bullet$} \rput(8.5,4.6){$\bullet$} \rput(9,4.1){$\bullet$} \rput(8.5,3.6){$\bullet$}
\end{pspicture}
\end{equation}
Here, the bulleted squares and triangles are all collapsed to a point, and all elements in the same dashed line are identified, \emph{i.e.} they are collapsed to the same point. Note that all the squares above the diagonal that are obtained by gluing two triangles, are collapsed by $P_{0,3}$ in the same way. Similarly, all squares below the diagonal that are obtained by gluing two triangles are collapsed by $P_{0,3}$ in the same way, which is symmetric (with respect to the diagonal) to the one above the diagonal.

For general $g>0$, the map $P_{0,g}:sd_2(\Sigma^g_\bullet)\to \Sigma^0\vee\Sigma^g_\bullet$ is defined similarly, using the same identifications for the diagonal squares and off-diagonal squares as for $P_{0,3}$.
 
\begin{lem}\label{L:P0g}
 $P_{0,g}:sd_2(\Sigma^g_\bullet)\to \Sigma^0\vee\Sigma^g_\bullet$ is a map of pointed simplicial sets.
\end{lem}
\begin{proof}
As in the proof of Lemma~\ref{L:modelglueing}, this follows from the fact that $P_{0,g}$ is obtained as a wedge along an edge or a vertex of collapse maps $\Triangle_\bullet \to pt_\bullet$ of a triangle to a point or of a triangle to an edge $\Triangle_\bullet \to I_\bullet$.
\end{proof}

\begin{remk}
The induced map ${P_{0,g}}_*:CH^{sd_2(\Sigma^g_n)}_\bullet(A,M) \to CH^{\Sigma^0_n\vee\Sigma^g_n}_\bullet(A,M)$ can be seen as follows. Recall from Examples~\ref{E:torus} and~\ref{E:2sphere} that each square in the model for $sd_2(\Sigma^g_\bullet)$ contributes to $(n+2)^2$-tensors in $CH^{sd_2(\Sigma^g_n)}_\bullet(A,M)$, which can be indexed as a $(n+2)\times (n+2)$-matrix. Similarly, each triangle contributes $(n+1)(n+2)+1$ tensors in $CH^{sd_2(\Sigma^g_n)}_\bullet(A,M)$, which can be indexed as $a\otimes M$, where $M$ is an $(n+2)\times (n+1)$-matrix. By construction, $sd_2(\Sigma^g)$ is obtained by gluing  subdivided squares $sd_2(I_\bullet^2)$ and triangles $sd_2(T_\bullet)=sd_2(I_\bullet)/\sim$ along edges and vertices. Then  ${P_{0,g}}_*:CH^{sd_2(\Sigma^g_n)}_\bullet(A,M) \to CH^{\Sigma^0_n\vee\Sigma^g_n}_\bullet(A,M)$ is the map which multiplies together the first $n+1$ columns and the first $n+1$ rows of the matrix corresponding to  each subdivided square (except for the top left square) or triangle. In other words, it is obtained by applying the $(n+1)$-th power $(d_0)^{\circ n}$ of the face map $d_0$ to each subdivided triangle or square (except for the top left square) in $sd_2(\Sigma^g_n)$.  
\end{remk}

We now define a left action of $\CH{0}{\bullet}(A,B)$ on $\CH{g}{\bullet}(A,B)$ which we will show to be equivalent to the one given in Definition~\ref{D:cup0} above.
\begin{defn}\label{D:edgewisecup}
 For $g\geq 1$, we define a cup-product $\tilde{\cup}$ as the composition
\begin{multline*}
 \tilde{\cup}: \CH{0}{i}(A,B)\otimes \CH{g}{j}(A,B) \\ \stackrel{AW} 
\to  \CH{0}{i+j}(A,B)\otimes \CH{g}{i+j}(A,B) \\
\stackrel{\vee} \to\CHW{0}{g}{i+j}(A,B) \\ \stackrel{P_{0,g}^*}\to CH_{sd_2(\Sigma^g_{i+j})}^\bullet (A,B)
\stackrel{\mathcal{D}_\bullet(2)^*}\to \CH{g}{i+j}(A,B).
\end{multline*}
\end{defn}

\begin{prop}\label{P:edgewisecup}
 The cup-product  $\tilde{\cup}:\CH{0}{\bullet}(A,B)\otimes \CH{g}{\bullet}(A,B)\to \CH{g}{\bullet}(A,B)$ is a cochain map. Furthermore, if $f\in \CH{0}{\bullet}(A,B)$ and $\alpha\in \CH{g}{\bullet}(A,B)$ are normalized cochains, then 
$f\cup \alpha= f\tilde{\cup} \alpha$. 
\end{prop}
In particular, Definition~\ref{D:edgewisecup} and Definition~\ref{D:cup0} coincide on normalized cochains and thus in cohomology.
\begin{proof}
 By Lemma~\ref{L:P0g}, ${P_{0,g}}^*$ is a morphism of cochain complexes. Since $AW$, $\vee$, and $\mathcal{D}_\bullet(2)$ are also chain maps, it follows that $\tilde{\cup}$ is a cochain map, too. 
 
Now, assume $f\in CH_{\Sigma^0_p}^{\bullet}(A,B)$ and $\alpha \in CH_{\Sigma^g_q}^{\bullet}(A,B)$ are  normalized cochains, and set $n=p+q$. Recall from Definition~\ref{normal}, that ``normalized" means that we divide the Hochschild chains $CH_\bullet^{Y_\bullet}(A,M)$ by the degeneracies, and dually we take the subcomplex of $CH^\bullet_{Y_\bullet}(A,M)$ vanishing on these degeneracies. In particular $f\big( (a_{ij})_{1\leq i,j\leq p}) \big)=0$ whenever there exists an $i$  such that $a_{i,j}=a_{j,i}=1$ for all $j$, {\it i.e.} if the matrix of the $(a_{ij})$ has only ones in the $i$-th column and the $i$-th row.

By definition of the edgewise subdivison functor, a cochain in $CH_{sd_2(\Sigma^g_{n})}^\bullet (A,B)$  is a linear map $A^{\otimes \sigma^g_{2n+1}}\to B$, where $\sigma^g_{2n+1}=\# \Sigma^g_{2n+1}-1$.
For any $x\in A^{\otimes \sigma^g_{2n+1}}$, note that $(f\tilde{\cup}\alpha)(x)$ is given by the composition,
\begin{equation*}
 (f\tilde{\cup}\alpha)(x)= (AW_{(1)}(f)\vee AW_{(2)}(\alpha))( (P_{0,g})_*\circ \mathcal D_\bullet (2)_*(x)).
 \end{equation*}
Here, $AW_{(1)}(f)\vee AW_{(2)}(\alpha):  A^{\otimes \# (\Sigma^0\vee \Sigma^g)_n -1}\cong A^{\otimes \sigma^0_n}\otimes A^{\otimes \sigma^g_n}\to B$ is given by mapping $x'\otimes x''\in A^{\otimes \sigma^0_n}\otimes A^{\otimes \sigma^g_n}$ to the product $AW_{(1)}(f)(x')\cdot AW_{(2)}(\alpha)(x'')$ in $B$.
Furthermore, by formula~\eqref{E:D(2)}, $(P_{0,g})_*\circ \mathcal D_\bullet (2)_*(x)\in A^{\otimes \sigma^0_n}\otimes A^{\otimes \sigma^g_n}$ is given by a sum of terms indexed by $(\sigma, \eta)\in \mathcal S(2,n)$,
$$  (P_{0,g})_*\circ \mathcal D_\bullet (2)_*(x)= \sum_{(\sigma,\eta)} \tilde{x}_{(\sigma,\eta)} \in A^{\otimes \sigma^0_n}\otimes A^{\otimes \sigma^g_n}. $$
We claim that $AW_{(1)}(f)\vee AW_{(2)}(\alpha)$ applied to $\tilde{x}_{(\sigma,\eta)}$ vanishes for all $(\sigma,\eta)$ except in one case $(\bar \sigma,\bar \eta)$, where $\bar\sigma=id_{\{1,\dots ,p+q\}}$ and $\bar \eta(i)= \begin{cases} 0\text{, for }i\leq p\\ 1\text{, for } i> p\end{cases}$. In fact, $\tilde x_{(\sigma,\eta)}=(P_{0,g})_* \circ \epsilon_{(\sigma,\eta)}^*(x)$,  where the map $\epsilon_{(\sigma,\eta)}:\Delta\to \Delta$ is defined by formula~\eqref{E:epsilonsigmadelta}. From formulas~\eqref{E:D(2)} and~\eqref{E:epsilonsigmadelta} we see that when $\sigma(1)\neq 1$ or $\eta(1)\neq 0$, we need to apply a degeneracy $(s_0)_*$ to $x$, so that the first row and the first column of the $A^{\otimes \sigma^0_n}$ factor of $\tilde{x}_{(\sigma,\eta)}$ are ones, and thus the normalized cochain $AW_{(1)}(f)$ is applied to a generate element, making the term vanishing. Similar arguments apply to $\sigma(i)\neq 2$ or $\eta(i)\neq 0$, for $i=2,\dots, p$. For $i>p$, and $\sigma(i)\neq i$ or $\eta(i)\neq 1$, we obtain a degenerate element in $A^{\otimes \sigma^g_q}$, vanishing on the $AW_{(2)}(\alpha)$ factor.

It is now straightforward to check, that $ (AW_{(1)}(f)\vee AW_{(2)}(\alpha))((P_{(0,g)})_*(\tilde x_{(\bar\sigma, \bar\eta)}))$ multiplies the tensor factors of $x$ exactly as in Definition~\ref{D:cup0}, showing that this is equal to $(f\cup g) (x).$
\end{proof}

In order to give a similar right action of $\CH{0}{\bullet}(A,B)$ on $\CH{g}{\bullet}(A,B)$, we define a pinching map 
 $P_{g,0}:sd_2(\Sigma^g_\bullet) \to \Sigma^0_\bullet\vee \Sigma^g_\bullet$.  The map $P_{3,0}:sd_2(\Sigma^3_\bullet) \to \Sigma^0_\bullet\vee \Sigma^3_\bullet$ is given by the following picture:
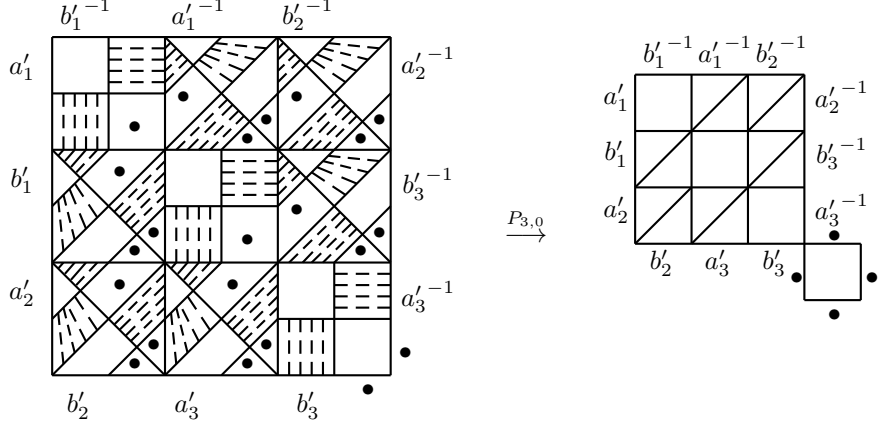
\begin{equation} \label{E:P30}
\begin{pspicture}(0,0)(11.5,5.5)
 \psline(0.5,.5)(0.5,5) \psline(.5,.5)(5,.5) \psline(2,.5)(2,5) \psline(3.5,.5)(3.5,5)
\psline(5,.5)(5,5) \psline(.5,5)(5,5) \psline(0.5,2)(5,2) \psline(.5,3.5)(5,3.5)
 \psline(.5,4.25)(2,4.25) \psline(2,2.75)(3.5,2.75) \psline(3.5,1.25)(5,1.25) 
 \psline(1.25,3.5)(1.25,5) \psline(2.75,2)(2.75,3.5) \psline(4.25,.5)(4.25,2) 
\psline(.5,2)(2,.5) \psline(0.5,3.5)(3.5,.5) \psline(2,5)(5,2) \psline(3.5,5)(5,3.5)
\psline(.5,2)(3.5,5) \psline(0.5,.5)(2,2) \psline(3.5,3.5)(5,5) \psline(2,.5)(5,3.5)
\psline(.5,2.75)(1.25,3.5) \psline(2,4.25)(2.75,5) 
\psline(.5,1.25)(2,2.75) \psline(2.75,3.5)(4.25,5)
\psline(1.25,.5)(2.75,2) \psline(3.5,2.75)(5,4.25) 
\psline(2.75,.5)(3.5,1.25) \psline(4.25,2)(5,2.75)

\rput(5.2,0.8){$\bullet$} \rput(4.7,0.3){$\bullet$}
 \rput(0.1,4.6){$a'_1$} \rput(0.1,3.1){$b'_1$} \rput(0.1,1.6){$a'_2$}  
\rput(0.85,0.1){$b'_2$} \rput(2.3,0.1){$a'_3$}  \rput(3.9,0.1){$b'_3$} 
\rput(5.5, 1.55){${a'_3}^{-1}$}  \rput(5.5,3.05){${b'_3}^{-1}$} \rput(5.5,4.65){${a'_2}^{-1}$} 
 \rput(3.9,5.3){${b_2'}^{-1}$} \rput(2.45,5.3){${a_1'}^{-1}$}   \rput(0.95,5.3){${b_1'}^{-1}$}

\rput(2.25, 4.2){$\bullet$} \rput(3.1, 3.65){$\bullet$} \rput(3.35, 3.9){$\bullet$} 
\rput(3.75, 4.2){$\bullet$} \rput(4.6, 3.65){$\bullet$} \rput(4.85, 3.9){$\bullet$} 
\rput(3.75, 2.7){$\bullet$} \rput(4.6, 2.15){$\bullet$} \rput(4.85, 2.4){$\bullet$} 

\rput(1.6, 2.15){$\bullet$} \rput(1.85, 2.4){$\bullet$}  \rput(1.4, 3.2){$\bullet$} 
\rput(1.6, .65){$\bullet$} \rput(1.85, .9){$\bullet$}  \rput(1.4, 1.7){$\bullet$} 
\rput(3.1, .65){$\bullet$} \rput(3.35, .9){$\bullet$}  \rput(2.9, 1.7){$\bullet$} 

\rput(1.6,3.8){$\bullet$} \rput(3.1,2.3){$\bullet$}

\psline[linestyle=dashed](.65,3.5)(.65,4.25) \psline[linestyle=dashed](.8,3.5)(.8,4.25) \psline[linestyle=dashed](.95,3.5)(.95,4.25) \psline[linestyle=dashed](1.1,3.5)(1.1,4.25)
\psline[linestyle=dashed](1.25,4.4)(2,4.4)\psline[linestyle=dashed](1.25,4.55)(2,4.55)
\psline[linestyle=dashed](1.25,4.7)(2,4.7)\psline[linestyle=dashed](1.25,4.85)(2,4.85)

\psline[linestyle=dashed](2.15,2)(2.15,2.75) \psline[linestyle=dashed](2.3,2)(2.3,2.75) \psline[linestyle=dashed](2.45,2)(2.45,2.75) \psline[linestyle=dashed](2.6,2)(2.6,2.75)
\psline[linestyle=dashed](2.75,2.9)(3.5,2.9)\psline[linestyle=dashed](2.75,3.05)(3.5,3.05)
\psline[linestyle=dashed](2.75,3.2)(3.5,3.2)\psline[linestyle=dashed](2.75,3.35)(3.5,3.35)

\psline[linestyle=dashed](3.65,.5)(3.65,1.25) \psline[linestyle=dashed](3.8,.5)(3.8,1.25) \psline[linestyle=dashed](3.95,.5)(3.95,1.25) \psline[linestyle=dashed](4.1,.5)(4.1,1.25)
\psline[linestyle=dashed](4.25,1.4)(5,1.4)\psline[linestyle=dashed](4.25,1.55)(5,1.55)
\psline[linestyle=dashed](4.25,1.7)(5,1.7)\psline[linestyle=dashed](4.25,1.85)(5,1.85)

\psline[linestyle=dashed](2,4.4)(2.3,4.7) \psline[linestyle=dashed](2,4.55)(2.23,4.77) \psline[linestyle=dashed](2,4.7)(2.15,4.85) \psline[linestyle=dashed](2,4.85)(2.07,4.92)
\psline[linestyle=dashed](2.45,4.7)(2.9,4.4)\psline[linestyle=dashed](2.52,4.77)(3.05,4.55)
\psline[linestyle=dashed](2.6,4.85)(3.2,4.7)\psline[linestyle=dashed](2.67,4.92)(3.35,4.85)
\psline[linestyle=dashed](2.15,3.5)(2.83,4.17)\psline[linestyle=dashed](2.3,3.5)(2.9,4.1)
\psline[linestyle=dashed](2.45,3.5)(2.98,4.02)\psline[linestyle=dashed](2.6,3.5)(3.05,3.95)

\psline[linestyle=dashed](3.5,4.4)(3.8,4.7) \psline[linestyle=dashed](3.5,4.55)(3.73,4.77) \psline[linestyle=dashed](3.5,4.7)(3.65,4.85) \psline[linestyle=dashed](3.5,4.85)(3.57,4.92)
\psline[linestyle=dashed](3.95,4.7)(4.4,4.4)\psline[linestyle=dashed](4.02,4.77)(4.55,4.55)
\psline[linestyle=dashed](4.1,4.85)(4.7,4.7)\psline[linestyle=dashed](4.17,4.92)(4.85,4.85)
\psline[linestyle=dashed](3.65,3.5)(4.33,4.17)\psline[linestyle=dashed](3.83,3.5)(4.4,4.1)
\psline[linestyle=dashed](3.95,3.5)(4.48,4.02)\psline[linestyle=dashed](4.1,3.5)(4.55,3.95)

\psline[linestyle=dashed](3.5,2.9)(3.8,3.2) \psline[linestyle=dashed](3.5,3.05)(3.73,3.27) \psline[linestyle=dashed](3.5,3.2)(3.65,3.35) \psline[linestyle=dashed](3.5,3.35)(3.57,3.42)
\psline[linestyle=dashed](3.95,3.2)(4.4,2.9)\psline[linestyle=dashed](4.02,3.27)(4.55,3.05)
\psline[linestyle=dashed](4.1,3.35)(4.7,3.2)\psline[linestyle=dashed](4.17,3.42)(4.85,3.35)
\psline[linestyle=dashed](3.65,2)(4.33,2.67)\psline[linestyle=dashed](3.83,2)(4.4,2.6)
\psline[linestyle=dashed](3.95,2)(4.48,2.52)\psline[linestyle=dashed](4.1,2)(4.55,2.45)

\psline[linestyle=dashed](0.65,3.5)(0.58, 3.42) \psline[linestyle=dashed](0.8,3.5)(.67,3.35)
\psline[linestyle=dashed](0.95,3.5)(.75,3.27) \psline[linestyle=dashed](1.1,3.5)(.82, 3.2)
\psline[linestyle=dashed](2,3.35)(1.33,2.68)\psline[linestyle=dashed](2,3.2) (1.4, 2.61)
\psline[linestyle=dashed](2,3.05)(1.48,2.53)\psline[linestyle=dashed](2,2.9)(1.55, 2.46)
\psline[linestyle=dashed](.65,2.15)(.58,2.83)\psline[linestyle=dashed](.8, 2.3) (.65, 2.9)
\psline[linestyle=dashed](.95,2.45)(.73, 2.98)\psline[linestyle=dashed](1.1, 2.6)(.8, 3.05)
  
\psline[linestyle=dashed](0.65,2)(0.58, 1.92) \psline[linestyle=dashed](0.8,2)(.67,1.85)
\psline[linestyle=dashed](0.95,2)(.75,1.77) \psline[linestyle=dashed](1.1,2)(.82, 1.7)
\psline[linestyle=dashed](2,1.85)(1.33,1.18)\psline[linestyle=dashed](2,1.7) (1.4, 1.11)
\psline[linestyle=dashed](2,1.55)(1.48,1.03)\psline[linestyle=dashed](2,1.4)(1.55, .96)
\psline[linestyle=dashed](.65,.65)(.58,1.33)\psline[linestyle=dashed](.8, .8) (.65, 1.4)
\psline[linestyle=dashed](.95,.95)(.73, 1.38)\psline[linestyle=dashed](1.1, 1.1)(.8, 1.55)
  
\psline[linestyle=dashed](2.15,2)(2.08, 1.92) \psline[linestyle=dashed](2.3,2)(2.17,1.85)
\psline[linestyle=dashed](2.45,2)(2.25,1.77) \psline[linestyle=dashed](2.6,2)(2.32, 1.7)
\psline[linestyle=dashed](3.5,1.85)(2.83,1.18)\psline[linestyle=dashed](3.5,1.7) (2.9, 1.11)
\psline[linestyle=dashed](3.5,1.55)(2.98,1.03)\psline[linestyle=dashed](3.5,1.4)(3.05, .96)
\psline[linestyle=dashed](2.15,.65)(2.08,1.33)\psline[linestyle=dashed](2.3, .8) (2.15, 1.4)
\psline[linestyle=dashed](2.45,.95)(2.23, 1.38)\psline[linestyle=dashed](2.6, 1.1)(2.3, 1.55)
\rput(6.8,2.5){$\stackrel{P_{3,0}}\longrightarrow$} 


\psline(8.25,4.5)(10.5,4.5) \psline(8.25,4.5)(8.25, 2.25) \psline(8.25, 3.75)(10.5,3.75) 
\psline(10.5, 1.5)(10.5,4.5)\psline(11.25,2.25)(8.25, 2.25) \psline(8.25, 3)(10.5,3)
\psline(9.75, 2.25)(9.75,4.5) \psline(8.25, 3)(9.75,4.5) \psline(8.25, 2.25)(9,3) 
\psline(9.75, 3.75)(10.5,4.5) \psline(10.5,1.5)(11.25,1.5) \psline(11.25,1.5)(11.25,2.25)
\psline(9,4.5)(9,2.25) \psline(9,2.25)(10.5,3.75)

\rput(8,4.2){$a'_1$}  \rput(8,3.45){$b'_1$} \rput(8,2.7){$a'_2$}  
\rput(8.6,2){$b'_2$} \rput(9.35,2){$a'_3$}  \rput(10.1,2){$b'_3$} 
\rput(11, 2.7){${a'_3}^{-1}$}  \rput(11,3.45){${b'_3}^{-1}$} \rput(11,4.2){${a'_2}^{-1}$} 
 \rput(10.2,4.8){${b'_2}^{-1}$} \rput(9.45,4.8){${a'_1}^{-1}$}  \rput(8.7,4.8){${b'_1}^{-1}$}
\rput(10.4,1.8){$\bullet$} \rput(10.9,1.3){$\bullet$} \rput(11.4,1.8){$\bullet$} \rput(10.9,2.35){$\bullet$}
\end{pspicture}
\end{equation}
Again, the bulleted squares and triangles are all collapsed to a point, and all elements in the same dashed line are identified, \emph{i.e.} they are collapsed to the same point. Note that all the squares above the diagonal that are obtained by gluing two triangles are collapsed by $P_{3,0}$ in the same way. And similarly all  the squares below the diagonal that are obtained by gluing two triangles are collapsed by $P_{3,0}$ in the same way, which is symmetric (with respect to the diagonal) to the one  above the diagonal.

For general $g>0$, the map $P_{g,0}:sd_2(\Sigma^g_\bullet)\to \Sigma^0_\bullet\vee\Sigma^g_\bullet$ is defined similarly, using the same identifications for the diagonal squares and off-diagonal squares as for $P_{3,0}$. Note that the identifications on the squares describing $P_{g,0}$ are symmetric to those describing $P_{0,g}$.

\smallskip

\begin{defn}\label{D:edgewisecupright}
 For $g\geq 1$, we define a right action by $\tilde{\cup}$ as the composition
\begin{multline*}
 \tilde{\cup}: \CH{g}{i}(A,B)\otimes \CH{0}{j}(A,B) \\ \stackrel{AW} 
\to  \CH{g}{i+j}(A,B)\otimes \CH{0}{i+j}(A,B) \\
\stackrel{\vee} \to\CHW{g}{0}{i+j}(A,B) \\ \stackrel{P_{g,0}^*}\to CH_{sd_2(\Sigma^g_{i+j})}^\bullet (A,B)
\stackrel{\mathcal{D}_\bullet(2)^*}\to \CH{g}{i+j}(A,B).
\end{multline*}
\end{defn}

An argument similar to the one of  Proposition~\ref{P:edgewisecup} shows that
\begin{prop}\label{P:edgewisecupright}
 The cup-product  $\tilde{\cup}:\CH{g}{\bullet}(A,B)\otimes \CH{0}{\bullet}(A,B)\to \CH{g}{\bullet}(A,B)$ is a cochain map. Furthermore, if $f\in \CH{0}{\bullet}(A,B)$ and $\alpha\in \CH{g}{\bullet}(A,B)$ are normalized cochains, then 
$\alpha\cup f= \alpha\tilde{\cup} f$. 
\end{prop}
In particular, Definition~\ref{D:edgewisecupright} and Definition~\ref{D:cup0} coincide on normalized cochains and therefore also in cohomology.

\smallskip \subsubsection{Properties of the cup product}

The cup product is not symmetric on cochains. However, for $B=A$, and passing to cohomology, we obtain 
\begin{prop}\label{P:cupbimodule}
 $\HHS{g}(A,A)$ is a (graded) symmetric $\HHS{0}(A,A)\cong HH^\bullet_{S^2}(A,A)$-bimodule. 
\end{prop} 
\begin{proof}
 By Lemma~\ref{L:cupbimodule} and Proposition~\ref{P:cup00}, we only need to prove that $f\cup \alpha =\alpha \cup f\in \HHS{g}(A,A)$ for any $f\in \HHS{0}(A,A)$ and $\alpha \in \HHS{g}(A,A)$. We are going to use an argument similar to the one from Proposition~\ref{P:bimodule}. To do so, we use the Hochschild cochain complexes $CH_{S_\bullet(|\Sigma^h_\bullet|)}^\bullet(A,A)$ of $A$ over the simplicial set $S_\bullet(|\Sigma^h_\bullet|)$ (see Definitions~\ref{Chen-map} and~\ref{Singular}). The natural map $\eta:\Sigma^h_\bullet\to S_\bullet(|\Sigma^h_\bullet|)$ induces the cochain map 
$$\eta^*: CH_{S_\bullet(|\Sigma^h_\bullet|)}^\bullet(A,A)\to CH_{\Sigma^h_\bullet}^\bullet(A,A)$$ which is a quasi-isomorphism by Lemma~\ref{L:quis}. 
Similarly the natural inclusion $S_\bullet(|\Sigma^h_\bullet|)\vee S_\bullet(|\Sigma^g_\bullet|)\stackrel{i}\hookrightarrow S_\bullet(|\Sigma^h_\bullet|\vee |\Sigma^g_\bullet|)$ induced 
by the canonical maps $\Sigma^g\hookrightarrow \Sigma^g\vee \Sigma^h$ and   $\Sigma^h\hookrightarrow \Sigma^g\vee \Sigma^h$ yields a 
quasi-isomorphism $$CH_{S_\bullet(|\Sigma^h_\bullet|\vee |\Sigma^g_\bullet|)}^\bullet(A,A)\stackrel{i^*}\to  CH_{S_\bullet(|\Sigma^h_\bullet|)\vee S_\bullet(|\Sigma^g_\bullet|)}^\bullet(A,A).$$ We define the map $\mu_{0,h}:HH_{S_\bullet(|\Sigma^0_\bullet|)}^\bullet(A,A) \otimes HH_{S_\bullet(|\Sigma^h_\bullet|)}^\bullet(A,A)\to HH_{S_\bullet(|\Sigma^h_\bullet|)}^\bullet(A,A)$ to be the composition
\begin{multline*}
 \mu_{0,h}:   HH_{S_\bullet(|\Sigma^0_\bullet|)}^\bullet(A,A) \otimes HH_{S_\bullet(|\Sigma^h_\bullet|)}^\bullet(A,A)\stackrel{\vee \circ AW} 
\longrightarrow  HH_{S_\bullet(|\Sigma^0_\bullet|\vee |\Sigma^h_\bullet|)}^\bullet(A,A)  \\
\stackrel{(i^*)^{-1}} \longrightarrow HH_{S_\bullet(|\Sigma^0_\bullet|)\vee S_\bullet(|\Sigma^h_\bullet|)}^\bullet(A,A)  \stackrel{\pinch{0,h}^*}\longrightarrow HH_{S_\bullet(|\Sigma^h_\bullet|)}^\bullet(A,A)
\end{multline*}
where the wedge map $\vee$ and Alexander-Whitney map $AW$ are defined as in Definition~\ref{D:cupg} and $\pinch{0,h}$ is the map~\eqref{E:sigma0hcollapse} defined in Section~\ref{S:surfacemodel}. Similarly we define the map $\mu_{h,0}:HH_{S_\bullet(|\Sigma^h_\bullet|)}^\bullet(A,A)\otimes HH_{S_\bullet(|\Sigma^0_\bullet|)}^\bullet(A,A) \to HH_{S_\bullet(|\Sigma^h_\bullet|)}^\bullet(A,A)$ as the composition
\begin{multline*}
 \mu_{h,0}:    HH_{S_\bullet(|\Sigma^h_\bullet|)}^\bullet(A,A)\otimes HH_{S_\bullet(|\Sigma^0_\bullet|)}^\bullet(A,A)  \stackrel{\vee \circ AW} 
\longrightarrow  HH_{S_\bullet(|\Sigma^h_\bullet|\vee |\Sigma^0_\bullet|)}^\bullet(A,A)  \\
\stackrel{(i^*)^{-1}} \longrightarrow HH_{S_\bullet(|\Sigma^h_\bullet|)\vee S_\bullet(|\Sigma^0_\bullet|)}^\bullet(A,A)  \stackrel{\pinch{h,0}^*}\longrightarrow HH_{S_\bullet(|\Sigma^h_\bullet|)}^\bullet(A,A).
\end{multline*}
Since  $\eta:X_\bullet\to S_\bullet(|X_\bullet|)$ is the natural map which sends any element $x\in X_n$ to the map $\eta(x):\Delta^n \stackrel{id\times x}\longrightarrow \coprod_{i\in \mathbb{N}} \Delta^i\times X_i \to |X|$ (see Definition~\ref{Chen-map}), there is a natural factorization 
\begin{eqnarray*}
 \xymatrix{X_\bullet\vee Y_\bullet \ar[rr]^{\eta_{\vee}} \ar[rrd]_{\eta} && S_\bullet(| X_\bullet|)\vee S_\bullet(|Y_\bullet|) \ar[d]^{i} \\
&&  S_\bullet(| X_\bullet|\vee |Y_\bullet | )}
\end{eqnarray*}
and furthermore the following diagrams are 
 commutative
\begin{eqnarray*}
 \xymatrix{CH_{S_\bullet(|\Sigma^0_\bullet|)\vee S_\bullet(|\Sigma^h_\bullet|)}^\bullet(A,A)  \ar[d]_{\eta^*_{\vee}}& CH_{S_\bullet(|\Sigma^0_\bullet|\vee |\Sigma^h_\bullet|)}^\bullet(A,A) \ar[d]_{\eta^*} \ar[l]_{i^*} \ar[r]^{ |P_{0,h}|^*}& CH_{S_\bullet(|sd_2(\Sigma^h_\bullet)|)}^\bullet(A,A)\ar[d]_{\eta^*}\\ 
CH_{(\Sigma^0\vee \Sigma^h)_\bullet}^\bullet(A,A)  & CH_{(\Sigma^0\vee \Sigma^h)_\bullet}^\bullet(A,A) \ar[l]_{id} \ar[r]^{ P_{0,h}^*}& CH_{sd_2(\Sigma^h_\bullet)}^\bullet(A,A),}\\
\xymatrix{CH_{S_\bullet(|\Sigma^0_\bullet|)}^\bullet(A,A)\otimes CH_{S_\bullet(|\Sigma^h_\bullet|)}^\bullet(A,A)  \ar[d]_{\eta^* \otimes \eta^*} \ar[rr]^{\;\;\vee \circ AW}&& CH_{S_\bullet(|\Sigma^0_\bullet|)\vee S_\bullet(|\Sigma^h_\bullet|)}^\bullet(A,A)\ar[d]_{\eta^*_\vee} \\ 
CH_{\Sigma^0_\bullet}^\bullet(A,A) \otimes CH_{\Sigma^h_\bullet}^\bullet(A,A) \ar[rr]_{\vee \circ AW} && CH_{(\Sigma^0\vee \Sigma^h)_\bullet}^\bullet(A,A) .}
\end{eqnarray*}
Now it follows from Proposition~\ref{P:edgewisecup} and the fact that $|P_{0,h}|\circ D^{-1}:\Sigma^h\to \Sigma^0\vee \Sigma^h$ is homotopic to $\pinch{0,h}^*$ that, for any $f\in   HH_{S_\bullet(|\Sigma^0_\bullet|)}^\bullet(A,A)$  and $\alpha\in HH_{S_\bullet(|\Sigma^h_\bullet|)}^\bullet(A,A)$ one has
\begin{equation}
 \eta^*(f)\cup \eta^*(\alpha) = \eta^*(\mu_{0,h}(f,\alpha)) \quad \mbox{ in }HH_{\Sigma^h_\bullet}^\bullet(A,A).
\end{equation}
In other words, $\eta$ is a map of left modules. 
Similarly, using Proposition~\ref{P:edgewisecupright} and the Pinching map $P_{h,0}$~\eqref{E:P30} instead of $P_{0,h}$, one proves that
\begin{equation}
 \eta^*(\alpha)\cup \eta^*(f) = \eta^*(\mu_{h,0}(\alpha,f)) \quad \mbox{ in }HH_{\Sigma^h_\bullet}^\bullet(A,A).
\end{equation}
Thus  $\eta$ is also a map of right modules and it is sufficient to prove that $\mu_{0,h}=\mu_{h,0}$ which easily follows from the fact that $\pinch{0,h}$ and $\pinch{h,0}$ are homotopic as in Proposition~\ref{P:bimodule}.
\end{proof}

We can now state the main result of this section.
\begin{thm}\label{T:cup} Let $(A,d_A)$ be a differential graded commutative algebra.
\begin{itemize} \item[i)] The cup product (Definition~\ref{D:cupg} and~\ref{D:cup0}) makes $\bigoplus_{g\geq 0}\HHS{g}(A,A)$ into an associative algebra which is bigraded with respect to the total degree grading and  the genus of the surfaces. Furthermore, $\bigoplus_{g\geq 0}\HHS{g}(A,A)$ is unital with unit being the cohomology class  $[1_A]\in \HH{0}{0}(A,A)\cong H^0(A,d_A)$. 
\item[ii)]   $\HHS{0}(A,A)$ lies in the center of $\bigoplus_{g\geq 0}\HHS{g}(A,B)$.\end{itemize}
\end{thm}
Note that, by construction, $\bigoplus_{g\geq 0}\HHS{g}(A,B)$ is also graded with respect to the cosimplicial degree and thus is in fact trigraded.
\begin{proof}
i) By Proposition~\ref{P:cup} and Lemma~\ref{L:cupbimodule} we are left to prove that for any $\alpha, \beta \in \HH{\bullet>0}{\bullet}(A,A)$ and $f\in \HH{0}{\bullet}(A,A)$ one has
\begin{eqnarray}\label{E:gh0} \alpha \cup ( \beta \cup f) &=&  \big(\alpha\cup \beta\big)\cup f, \\
\label{E:0gh} (f\cup \alpha) \cup \beta &=& f\cup \big(\alpha\cup \beta\big) \qquad \mbox{ and }\\ \label{E:g0h} \big(\alpha \cup f\big) \cup \beta   &=& \alpha \cup \big(f\cup \beta\big) .\end{eqnarray}
It is straightforward to check that the two first identities~\eqref{E:0gh} and~\eqref{E:gh0}  hold already for cochains. It follows from Proposition~\ref{P:cupbimodule} and identities~\eqref{E:gh0} and~\eqref{E:0gh} that
$$\big(\alpha \cup f\big) \cup \beta = \big(f\cup \alpha\big) \cup \beta = f\cup \big(\alpha\cup \beta\big) =\big(\alpha\cup \beta\big)\cup f=\alpha \cup \big(f\cup \beta\big)$$
hence identity~\eqref{E:g0h} holds.

\smallskip

According to its definition, the cup-product is graded with respect to the cosimplicial degree, total degree and genus degree on cochains, and hence in cohomology. Let $a\in CH_{\Sigma^0_0}^{\bullet}(A,A)\cong A$. Then for any $\alpha \in CH_{\Sigma^{g}_n}^{\bullet}(A,A)$ (with $g,n\geq 0$), one has $a\cup \alpha = a\cdot \alpha $ (where $\cdot$ is the multiplication in $A$). Similarly $\alpha \cup a =\alpha \cdot a$. In particular, $[1_A]$ is a unit for the cup-product and statement i) follows. 

\smallskip

ii) is an obvious corollary of Proposition~\ref{P:cupbimodule}. 
\end{proof}

\begin{remk}
Neither Theorem~\ref{T:cup} (i) nor part (ii) hold at the cochain complex level: $\bigoplus_{g\geq 0}\CHS{g}(A,B)$ is not  associative. In fact for any $f \in \CH{0}{k\geq 1}$, $\beta\in  \CH{g\geq 1}{l\geq 1}$, $\gamma\in  \CH{h\geq 1}{m\geq 1}$, a straightforward inspection shows that
$$(\beta \cup f)\cup \gamma \neq \pm \beta \cup (f \cup \gamma), \text{ and }  \beta\cup f\neq \pm f\cup \beta.$$
Also note that Theorem~\ref{T:cup} (i) can be proved by an argument similar to the one of Proposition~\ref{P:cupbimodule}, namely by using the homotopy associativity  of the maps $\pinch{h,0}$ and $\pinch{0,g}$ and the singular model $CH_{S_\bullet(|\Sigma^g_\bullet|)}^\bullet(A,A)$ for the Hochschild cohomology modeled on a surface of genus $g$.
\end{remk}

The cup product is natural and homotopy invariant.
\begin{prop}\label{P:cupinvariance} Let $B$ be a commutative $A$-algebra. \begin{itemize}
\item If $B\stackrel{f}\to B'$ is a quasi-isomorphism of $A$-algebras, then $f_*: \bigoplus_{g\geq 0} \HHS{g}(A,B)\to \bigoplus_{g\geq 0} \HHS{g}(A,B')$ is an isomorphism of algebras.
\item If $A'\stackrel{g}\to A$  is a quasi-isomorphism of algebras, then $g_*: \bigoplus_{g\geq 0} \HHS{g}(A,B)\to \bigoplus_{g\geq 0} \HHS{g}(A',B)$ is an isomorphism of algebras.
\end{itemize}
\end{prop}
\begin{proof}
This follows from Lemma~\ref{L:quis}. 
\end{proof}
Since quasi-isomorphic differential graded commutative algebras are connected by a zigzag of quasi-isomorphism of algebras, Proposition~\ref{P:cupinvariance} has an immediate Corollary.
\begin{cor}\label{C:cupinvariance} Let $A$ and $A'$ be quasi-isomorphic differential graded commutative algebras. Then $\bigoplus_{g\geq 0}\HH{g}{\bullet}(A,A)$ and $\bigoplus_{g\geq 0}\HHS{g}(A',A')$ are naturally isomorphic as algebras. 
\end{cor}

\smallskip

There is a (pointed) simplicial map $\pi^g_\bullet: \Sigma^g_\bullet \to S^2_\bullet$ obtained by collapsing all but the top left square in the simplicial model  $\Sigma^g_\bullet$ (see picture~\eqref{E:sigmag}) to a point. Note that in particular it collapses the boundary of this top left square to a point. Similarly to the topological situation (Proposition~\ref{P:pig}), this yields a map $\HHS{0}(A,B)\stackrel{(\pi^g_\bullet)^*}\longrightarrow \HHS{g}(A,B)$. 
\begin{prop}\label{P:Hpig} Let $B$ be a commutative $A$-algebra. Then
\begin{itemize}
\item The map $(\pi^g_\bullet)^*$ is an $\HHS{0}(A,B)$-module morphism.
\item If $B$ is unital, then $(\pi^g_\bullet)^*(\alpha)= \alpha \cup  [1_B]_g$ where $[1_B]_g\in \HH{g}{0}(A,B)$ is the class of $1_B$.
\end{itemize}
\end{prop}
\begin{proof}
 Let $\alpha\in CH_{\Sigma^0_p}^{\bullet}(A,B)$ and $\beta \in CH_{\Sigma^0_q}^{\bullet}(A,B)$ and $x\in A^{\otimes \sigma^g_{p+q}}$ be a homogeneous element, where $\sigma^g_{p+q}=\# \Sigma^g_{p+q} -1$. We can write $x= \big( \bigotimes_{\scriptsize i,j\leq p+q} a_{i,j} \big)\otimes y$, where the $a_{i,j}$'s are the tensor factors of $x$ corresponding to the top left square of $\Sigma^g_{p+q}$. Furthermore, $y$ can be written as a tensor $y=y_1\otimes \cdots \otimes y_{s^g_{p+q}}$ where $s^g_{p+q}=\sigma^g_{p+q}-(p+q)^2$. Formula~\eqref{eq:cup00} and Definition~\ref{D:cup0} imply that
\begin{eqnarray*}
(\pi^g_\bullet)^*(\alpha\cup \beta) (x)&=&(\alpha\cup \beta) \left((a_{i,j})_{\scriptsize i,j=1\dots p+q}\right) \cdot \prod_{k=1}^{s^g_{p+q}} y_k  \\ 
&=& \alpha\big((a_{i,j})_{\scriptsize i,j\leq p}\big)\cdot \beta\big((a_{i,j})_{\scriptsize i,j\geq p+1}\big) \cdot \hspace{-0.3cm} \prod_{\scriptsize \begin{array}{l}i\leq p \\
j\leq q\end{array}} \hspace{-0.2cm} a_{i,p+j} \cdot a_{j+p,i}\cdot \prod_{k=1}^{s^g_{p+q}} y_k\\
&=&\alpha \cup (\pi^g_\bullet)^*(\beta)(x). 
\end{eqnarray*}
Note that $1_B\in CH_{\Sigma^g_0}^{0}(A,B)$ has cosimplicial degree $0$. Since $AW_{(2)}: [p]\to [0]$ is the unique map to $\{0\}$, we get from Definition~\ref{D:cup0} that for any $x=\big(\bigotimes\limits_{\scriptsize i,j\leq p} a_{i,j}\big)\otimes y\in A^{\otimes \sigma^g_p}$, one has $$\alpha \cup  [1_B]_g(x)= \alpha \big((a_{i,j})_{\scriptsize i,j\leq p}\big) \cdot \prod_{k=1}^{s^g_{p}} y_k = (\pi^g_\bullet)^*(\alpha\cup 1_B)(x).$$
\end{proof}

\begin{remk}\label{r:modulehom}
Note that there is a simplicial inclusion $inc_\bullet:pt_\bullet\to\Sigma ^g_\bullet$ and projection $proj_\bullet:\Sigma^g_\bullet\to pt_\bullet$ between the point and the surface, with $proj_\bullet\circ inc_\bullet=id_{pt_\bullet}$. Since $B\cong CH^\bullet_{P_0}(A,B)$, we see that $B$ becomes a subcomplex of $\CH{g}{\bullet}(A,B)$ with a natural splitting induced by $inc_\bullet$ and $proj_\bullet$. Thus, $H^\bullet(B)$ is a direct summand of $\HH{g}{0}(A,B)$. 
\end{remk}

\begin{remk}
 Let $M$ be a differential graded $A$-module. Since the pinching maps are pointed, one can extend Definition~\ref{D:cupg}, Definition~\ref{D:cup0}, the results of Theorem~\ref{T:cup}, and Proposition~\ref{P:cupinvariance} to give to $\bigoplus_{g\geq 0} \HH{g}{\bullet}(A,M)$ the structure of a $\bigoplus_{g\geq 0} \HH{g}{\bullet}(A,A)$-bimodule, which is natural and homotopy invariant.
\end{remk}

\subsection{Topological identification of the cup product}

Let $M$ be a simply connected compact manifold and denote $\Omega=\Omega^\bullet M$ its de Rham cochain algebra and $\Omega^*=Hom(\Omega, k)$ its dual. By Theorem~\ref{T:surfaceproduct}, $(\mathbb{H}_\bullet(\Map(\Sigma^\bullet,M)),\scup)$ is an associative bigraded algebra. So is $(\HHS{\bullet}(\Omega, \Omega),\cup)$ by Theorem~\ref{T:cup}. In this section, we show that, similarly to the situation in string topology~\cite{CJ, FTV, FT}, the algebraic and topological constructions coincide. First notice that
\begin{lem}\label{L:PD}
There are natural ``Poincar\'e duality'' isomorphisms
$$\mathcal{P}:   \HH{g}{\bullet-\dim(M)}(\Omega, \Omega^*)\stackrel{\simeq} \to  \HH{g}{\bullet}(\Omega, \Omega), \quad \mathcal{P}: HH^{\Sigma^g_\bullet}_{\bullet}(\Omega, \Omega) \stackrel{\simeq} \to   HH^{\Sigma^g_\bullet}_{\bullet-\dim(M)}(\Omega, \Omega^*)$$ which are functorial with respect to smooth oriented maps between manifolds of the same dimension. 
\end{lem}
\begin{proof}
The lemma follows since the natural map $\int:\Omega \to \Omega^*$, $\omega \mapsto \int \omega\wedge -$ is a bimodule quasi-isomorphism. 
\end{proof}

Using Section~\ref{S:Chen-map}, we have the Chen iterated integral morphisms $(\Ch^{\Sigma^g_\bullet})^*:H_\bullet (\Map(\Sigma^g,M)) \to \HH{g}{-\bullet}(\Omega,\Omega^*)$ which is an isomorphism if $M$ is 2-connected, see Corollary~\ref{qi*}. Composing the iterated integral map with Poincar\'e duality from Lemma~\ref{L:PD}, yields a linear map
\begin{equation}\label{eq:ChenPoincare}
\Ch^{\Sigma^\bullet}: \bigoplus_{g\geq 0} \mathbb{H}_\bullet(\Map(\Sigma^g,M)) \stackrel{\oplus (\Ch^{\Sigma^g_\bullet})^*}\longrightarrow \bigoplus_{g\geq 0} \HH{g}{-\bullet-\dim(M)}(\Omega, \Omega^*) \stackrel{\oplus\mathcal{P}}\longrightarrow \bigoplus_{g\geq 0} \HH{g}{-\bullet}(\Omega, \Omega)
\end{equation}
that we call the \dualized{.}
\begin{thm}\label{T:surface=cup} Let $M$ be a 2-connected compact manifold. The \dualized{ } 
 $\Ch^{\Sigma^\bullet}:(\bigoplus_{g\geq 0 } \mathbb{H}_\bullet(\Map(\Sigma^g,M)),\scup)  \to (\bigoplus_{g\geq 0 }\HH{g}{-\bullet}(\Omega,\Omega),\cup) $ is an  isomorphism of algebras.
\end{thm}
The proof of Theorem~\ref{T:surface=cup} is given in Section~\ref{S:surface=cup} below. 
\begin{cor} \label{C:homotopy-inv} Let $M,N$ be 2-connected compact manifolds with equal dimensions, and let $i:M\to N$ be a homotopy equivalence. Then $$i_*: \big(\bigoplus_{g\geq 0}\mathbb{H}_\bullet(\Map(\Sigma^g,M)),\scup\big) \to  \big(\bigoplus_{g\geq 0}\mathbb{H}_\bullet(\Map(\Sigma^g,N)),\scup\big)$$ is an isomorphism of algebras.
\end{cor}
In particular, the surface product is  homotopy invariant for 2-connected manifolds.
\begin{remk}
The evaluation map $e^g:\Map(\Sigma^g,M) \to M$ has a section $i^g:M\to \Map(\Sigma^g,M)$ given by the constant surfaces at a point. It follows that $H_\bullet(\Map(\Sigma^g,M))$ contains $M$ as a direct summand. It is easy to check that this direct summand coincides with the summand $H^\bullet(CH^\bullet_{pt_\bullet}(\Omega,\Omega))$ from Remark~\ref{r:modulehom} under the isomorphism $\Ch^{\Sigma^\bullet}$. In particular, $\Ch^{\Sigma^\bullet}([1_{\Omega}])=[M]_0$ and it follows from Proposition~\ref{P:pig} and Proposition~\ref{P:Hpig}, that $\pi_g$ coincides with $\pi^g_\bullet$ under the dualized iterated integral map.
\end{remk}

\subsection{Proof of Theorem~\ref{T:surface=cup}} \label{S:surface=cup}
We follow an idea of F\'elix-Thomas~\cite{FT}, using rational homotopy theory techniques. To do so, we need to consider dual analogues of the surface product and cup product. 

\medskip

The construction of the surface product is easily dualized. Similarly to Section~\ref{thom}, the embedding $\rho_{in}: \Map(\Sigma^{g}\vee \Sigma^{h},M)\to \Map(\Sigma^g,M)\times \Map(\Sigma^{h},M) $  of  codimension $dim(M)$ induces an  Umkehr map in cohomology 
\begin{align*}
(\rho_{in})^! : H^{\bullet}(\Map(\Sigma^{g}\vee \Sigma^{h},M) \to &H^{\bullet+m}(\Map(\Sigma^g,M)\otimes \Map(\Sigma^h,M)) \\ & \cong \big(H^{\bullet}(\Map(\Sigma^g,M))\otimes H_{\bullet}(\Map(\Sigma^{h},M))\big)^{\bullet +m},
\end{align*} dual to $(\rho_{in})_!$.
Thus, for $k=g+h$, we can define a  linear map 
$$\cop{g,h}:  H^{\bullet-dim(M)}(\Map(\Sigma^{k},M))\to H^\bullet (\Map(\Sigma^g,M))\otimes H^{\bullet}(\Map(\Sigma^{h},M))$$ as the composition
\begin{multline*} 
 \cop{g,h}: H^{\bullet-dim(M)}(\Map(\Sigma^{k},M)) \stackrel{(\rho_{out})^*}{\longrightarrow}  H^{\bullet-dim(M)}(\Map(\Sigma^{g}\vee \Sigma^{h},M))   \\ \stackrel{(\rho_{in})^!}{\longrightarrow} H^{\bullet}(\Map(\Sigma^g,M))\otimes H^{\bullet}(\Map(\Sigma^{h},M)).
\end{multline*} 

\begin{lem}\label{L:surfaceduality}
The surface product $\scup: H_{\bullet}(\Map(\Sigma^g,M))\otimes H_{\bullet}(\Map(\Sigma^{h},M))\to  H_{\bullet-dim(M)}(\Map(\Sigma^{k},M))$ is the dual of the map $\cop{g,h}$.
\end{lem}

We now want to dualize the Hochschild cup product for surfaces. Since $M$ is a Poincar\'e duality space, by the main result of~\cite{LS}, there exists a differential graded commutative algebra $(A,d)$, weakly equivalent to $(\Omega^\bullet M, d_{dR})$, which is simply connected, finite dimensional and is equipped with a trace $A^{dim(M)}\stackrel{\epsilon}\to \R$ such that:\begin{itemize}
 \item the pairings $A^{i}\otimes A^{dim(M)-i}\stackrel{\cdot}\to A^{dim(M)}\stackrel{\epsilon}\to\R$ are non degenerate (where the first map is the multiplication in $A$);
\item $\epsilon \circ d=0$;
\item the induced pairing on  cohomology  $\langle\cdot \, , \cdot\rangle: H^\bullet(A)\otimes H^{dim(M)-\bullet}(A)\to \R$ coincides with the Poincar\'e duality pairing of $H^\bullet(\Omega^\bullet M)\cong H^\bullet(M)$ through the weak-equivalence between $A$ and $\Omega$.
\end{itemize}
It follows that the map $a\mapsto \langle a,\cdot\rangle$ is a linear isomorphism of symmetric $A$-bimodules $\Xi:A^\bullet\to (A^*)^{\bullet -dim(M)}$ and furthermore the composition
\begin{equation} \label{eq:mu}
\mu: A^* \otimes A^* \stackrel{\Xi^{-1} \otimes \Xi^{-1}}\longrightarrow A\otimes A \stackrel{\cdot} \longrightarrow A \stackrel{\Xi}\longrightarrow A^*
\end{equation}
is a degree $+\dim(M)$ graded commutative and associative multiplication on $(A^*,d^*)$. Note that $\mu$ is a model for the umkehr map  $H_\bullet(M)\otimes H_\bullet(M)\cong H_\bullet(M\times M) \stackrel{diag^!} \longrightarrow H_{\bullet+dim(M)}(M)$ of the diagonal $M\stackrel{diag}\longrightarrow M\times M$. By Proposition~\ref{P:cupinvariance}, there is a natural isomorphism of algebras $$(\bigoplus_{g\geq 0} \HHS{g}(A,A),\cup)\cong (\bigoplus_{g\geq 0}\HHS{g}(\Omega,\Omega),\cup).$$
The map $\mu$ from~\eqref{eq:mu} has the composition
\begin{equation} \label{eq:nabla}
\nabla : A \stackrel{\Xi}\longrightarrow    A^*    \stackrel{\cdot^*=(\mu^*)^*} \longrightarrow A^*   \otimes A^* \stackrel{\Xi^{-1} \otimes \Xi^{-1}}\longrightarrow A\otimes A 
\end{equation}
as a dual map. Clearly, $\nabla$ is a model for the umkehr map $H^{\bullet+dim(M)}(M) \stackrel{diag_!}\longrightarrow H^\bullet(M\times M) \cong H^\bullet(M)\otimes H^\bullet(M)$.
Since $\Xi:A\to A^*$ is an isomorphism of $A$-bimodules (of degree $-dim(M)$), $\Xi_*: \HHS{g}(A,A)\to  \HH{g}{\bullet-dim(M)}(A,A^*)$ is an isomorphism, hence there is a duality isomorphism  \begin{equation} \label{eq:Theta}\Theta:\HH{g}{\bullet+dim(M)}(A,A)\cong\HH{g}{\bullet}(A,A^*)\cong \big(HH^{\Sigma^g_\bullet}_{\bullet}(A,A)\big)^*.\end{equation}

\smallskip

Further, since $A$ is commutative, the multiplication $A\otimes A \stackrel{\cdot}\to A$ makes $A$ an $A\otimes A$-module and, since $\Xi$ an isomorphism of $A$-bimodules, the map $\nabla:A\to A\otimes A$ above~\eqref{eq:nabla} is a map of $A\otimes A$-modules.  For any $k\in \Delta$, $\nabla$ induces a linear map 
\begin{multline}\label{eq:nablatilde}
\tilde\nabla:CH^{(\Sigma^g\vee\Sigma^h)_k}_{\bullet+dim(M)}({A,A}) \stackrel{\nabla_*}\to  CH_{\bullet}^{(\Sigma^g\vee\Sigma^h)_k}({A,A\otimes A}) \cong CH^{\Sigma^g_k}_\bullet(A,A) \otimes CH^{\Sigma^h_k}_\bullet(A, A)\end{multline} where the last isomorphism follows as in Example~\ref{E:wedge} and $\nabla_*$ is the result of applying $\nabla$ to the sole module in the Hochschild complex (not the algebra). 
\begin{lem} \label{L:nablatilde}The map $\tilde\nabla:CH^{(\Sigma^g\vee\Sigma^h)_\bullet}({A,A})\to CH^{\Sigma^g_\bullet}_\bullet(A,A) \times CH^{\Sigma^h_\bullet}_\bullet(A, A)$, where the right hand side is the tensor product equipped with the diagonal simplicial structure ({\it cf. Definition~\ref{cross}}), is a morphism of the underlying chain complexes.
\end{lem} 
\begin{proof}
Note that that there is a canonical identification $CH^{\Sigma^g_\bullet}_\bullet(A,A) \times CH^{\Sigma^h_\bullet}_\bullet(A, A)\cong CH^{(\Sigma^g\sqcup\Sigma^h)_\bullet}_\bullet(A,A)$ and furthermore that  $CH^{pt_\bullet}_\bullet(A,A)$ and $CH^{(pt\coprod pt)_\bullet}_\bullet (A,A)$ are  the constant simplicial algebras $A$ and $A\otimes A$ respectively, see Example~\ref{E:wedge}. Hence it follows from Lemma~\ref{L:pushout} that $\tilde\nabla$ is the composition 
\begin{multline*} CH^{(\Sigma^g\vee\Sigma^h)_\bullet}_\bullet(A,A) \cong A\mathop{\otimes}\limits_{A\otimes A} CH^{(\Sigma^g\sqcup\Sigma^h)_\bullet}_\bullet(A,A)\\
\stackrel{\nabla \mathop{\otimes}\limits_{A\otimes A} id}\longrightarrow (A\otimes A)\mathop{\otimes}\limits_{A\otimes A}CH^{(\Sigma^g\sqcup\Sigma^h)_\bullet}_\bullet(A,A) \\ \cong CH^{\Sigma^g_\bullet}_\bullet(A,A) \times CH^{\Sigma^h_\bullet}_\bullet(A, A).\end{multline*} Since $\nabla$ is a map of $A\otimes A$-modules, the result follows.
\end{proof}

\subsubsection{Positive genus surfaces} 
We first consider the case of surfaces of positive genus.

\smallskip

Since $\Delta$ is the dual (through the duality isomorphism $\Xi$) of the multiplication $A\otimes A \stackrel{\cdot} \to A$, we deduce from Lemma~\ref{L:nablatilde} and the definition~\eqref{eq:Theta} of $\Theta$ the following Lemma.
\begin{lem}\label{L:Aduality} For $g,h>0$,
the duality isomorphism $\Theta$ (given by~\eqref{eq:Theta}) identifies the cup product $ \HHS{g}(A,A)\otimes   \HHS{h}(A,A)\stackrel{\cup}\to \HHS{g+h}(A,A)$ with the composition
\begin{multline*}
\Delta^{g,h}: HH^{\Sigma^{g+h}_\bullet}_{\bullet+dim(M)}(A,A) \stackrel{(\pinch{g,h})_*}{\longrightarrow} HH^{(\Sigma^{g}\vee \Sigma^{h})_\bullet}_{\bullet+dim(M)}(A,A) \\
\qquad \stackrel{\tilde\nabla}{\longrightarrow}  H_{\bullet}\big(CH^{\Sigma^g_\bullet}_\bullet(A,A) \times CH^{\Sigma^h_\bullet}_\bullet(A, A)\big) \\ \stackrel{AW}{\longrightarrow}  HH^{\Sigma^g_\bullet}_\bullet(A,A) \otimes HH^{\Sigma^h_\bullet}_\bullet(A, A)
\end{multline*}
\end{lem}

By~\cite{P}, there is a natural weak equivalence $CH_\bullet^{\Sigma^g}(A,A)\cong CH_\bullet^{\Sigma^g}(\Omega,\Omega)$ (for any genus $g$). Since $CH_\bullet^{\Sigma^g_\bullet}(A,A)$ is a semi-free model (see~\cite[Section 7]{FHT}) of $A$ as an $A$-bimodule, and $\Ch^{\Sigma^g}:CH_\bullet^{\Sigma^g_\bullet}(\Omega,\Omega)\to C^\bullet(M^{\Sigma^g})$ is a quasi-isomorphism (Proposition~\ref{qi}), it follows that $CH_\bullet^{\Sigma^g}(A,A)$ is a cochain model for $\Map(\Sigma^g,M)$. 
\begin{prop}\label{P:dualsurface=cup} If $g,h>0$, then the map $\Delta^{g,h}: HH^{\Sigma^{g+h}_\bullet}_{\bullet+dim(M)}(A,A)\to HH^{\Sigma^g_\bullet}_\bullet(A,A) \otimes HH^{\Sigma^h_\bullet}_\bullet(A, A)$ (defined in Lemma~\ref{L:Aduality}) is a cochain model for the map $\cop{g,h}:  H^{\bullet-dim(M)}(\Map(\Sigma^{k},M))\to H^\bullet (\Map(\Sigma^g,M))\otimes H^{\bullet}(\Map(\Sigma^{h},M)).$
\end{prop}
\begin{proof} 
The Alexander-Whitney map of simplicial modules and K\"unneth formula yield an isomorphism $ H_{\bullet}\big(CH^{\Sigma^g_\bullet}_\bullet(A,A) \times CH^{\Sigma^h_\bullet}_\bullet(A, A)\big)\cong  HH^{\Sigma^g_\bullet}_\bullet(A,A) \otimes HH^{\Sigma^h_\bullet}_\bullet(A, A)$. Since $CH_\bullet^{\Sigma^g}(A,A)$, $CH_\bullet^{\Sigma^g}(A,A)$ are models for $\Map(\Sigma^{h},M))
$ and $\Map(\Sigma^{h},M)$, we  are left to prove that the maps $(\pinch{g,h})_*$ and $\tilde \nabla$ in Lemma~\ref{L:Aduality} are respectively cochain models of $(\rho_{out})^*$ and $(\rho_{in})^!$. Thus, the result follows from Lemma~\ref{L:pinchdual} and Lemma~\ref{L:umkehrdual} below.
\end{proof}

The next Lemma gives a model for the evaluation map $ev:\Map(\Sigma^g,M)\to M$. There is a canonical quasi-isomorphism of differential graded algebras $(A,d_A)=(CH^{pt_0}_\bullet(A,A),d_A)\hookrightarrow (CH^{pt_\bullet}_\bullet(A,A),D)$, see Example~\ref{E:wedge}. By composition with the unique pointed map  $pt_\bullet\to \Sigma^g_\bullet$, it yields  the map  $e:(A,d_A)\to  (CH^{\Sigma^g_\bullet}_\bullet(A,A),D)$ which is a map of differential graded algebras and thus a map of $A$-modules.  Clearly the action of $A$  on $CH^{\Sigma^g_\bullet}_\bullet(A,A)$ is by multiplication on the module tensor of the Hochschild complex.
\begin{lem}\label{L:semifree} The map $e:A\to CH^{\Sigma^g_\bullet}_\bullet(A,A)$ (for any $g$) is a semi-free model (see~\cite{FHT}) for the evaluation map. The same holds with $\Sigma^{g}\vee \Sigma^{h}$ in place of $\Sigma^g$.
\end{lem}
\begin{proof} It is immediate that $CH^{\Sigma^g_\bullet}_\bullet(A,A)$ is $A$-semi-free (since $A$ is free over $\R$) and that $e$ is $A$-linear.
 Then, by functoriality of the iterated Chen integral, we have a commutative diagram $$\xymatrix{H^\bullet(A)\cong HH^{pt_\bullet}_\bullet(A,A) \ar[d]_{(pt_\bullet\hookrightarrow \Sigma^g_\bullet)_*)} \ar[r]^{\qquad \simeq} &  HH^{pt_\bullet}_{\bullet}(\Omega,\Omega) \ar[d]_{(pt_\bullet\hookrightarrow \Sigma^g_\bullet)_*)}\ar[r]^{\Ch^{pt_\bullet}} & H^*(M) \ar[d]^{ev^*}\\
\qquad HH^{\Sigma^g_\bullet}_\bullet(A,A) \ar[r]^{\qquad \simeq} & HH^{\Sigma^g_\bullet}_\bullet(\Omega,\Omega) \ar[r]^{\!\!\!\!\Ch^{\Sigma^g_\bullet}} & H^\bullet(\Map(\Sigma^g,M)) }$$ and the result follows. The argument for $\Sigma^{g}\vee \Sigma^{h}$ is the same.
\end{proof}

We now need the following  fact from rational homotopy theory~\cite[Section 7]{FHT}: given a pullback diagram  
\begin{equation*}
\xymatrix{
  X\times_Z Y  \ar[r]^{\quad i} \ar[d] 
  & Y \ar[d]^{p} \\
 X  \ar[r]^{j} & Z }
\end{equation*}
where $p:Y\to Z$ is a fibration, $A_Z$ a cochain algebra model for $Z$ and $B_X$, $B_Y$ two $A_Z$-semi free models for $X,Y$, then a model for the (homotopy) pullback $X\times_Z Y$ is given by the pushout $A_X \otimes_{A_Z} A_Y$. Furthermore, if $\widetilde{j}: A_Z \to B_X$ is a ($A_Z$-linear) model for $j:X\to Y$, then $\widetilde{j}\otimes_{A_Z} id_{B_Y}:B_Y\cong A_Z\otimes_{A_Z}B_Y\to A_X \otimes_{A_Z} A_Y$ is a model for $i:X\to Y$.

\begin{lem}\label{L:pinchdual}
For $g,h>0$, the following diagram is commutative: $$\xymatrix{ HH^{\Sigma^{g+h}_\bullet}_{\bullet}(A,A) \ar[r]^{(\pinch{g,h})_*} \ar[d]_{\simeq}&  HH^{(\Sigma^{g}\vee \Sigma^{h})_\bullet}_{\bullet}(A,A) \ar[d]^{\simeq} \\
 H^{\bullet}(\Map(\Sigma^{g+h},M)) \ar[r]^{(\rho_{out})^*}&  H^{\bullet}(\Map(\Sigma^{g}\vee \Sigma^{h},M)) } $$
\end{lem} 
\begin{proof}
Since the pinching map $\pinch{g,h}$ is simplicial and  Hochschild homology over  simplicial sets is a covariant functor on the algebras and on the simplicial sets, it is sufficient to prove the result with $\Omega^\bullet M$ in place of $A$. Now the result follows from the functoriality of the iterated Chen integral $\Ch^{Y_\bullet}:HH^{Y_\bullet}_\bullet (\Omega^\bullet M,\Omega^\bullet M)\to C^\bullet(M^Y)$ with respect to $Y$. 
\end{proof}

\begin{lem}\label{L:umkehrdual} For any $g,h$,
 $ \tilde \nabla\!:\! CH^{(\Sigma^{g}\vee \Sigma^{h})_\bullet}_{\bullet+dim(M)}(A,A)\to CH^{\Sigma^g_\bullet}_\bullet(A,A) \times CH^{\Sigma^h_\bullet}_\bullet(A, A)$ is a semi-free model for the Umkehr map $(\rho_{in})^! : H^{\bullet-dim(M)}(\Map(\Sigma^{g}\vee \Sigma^h,M)) \to H^{\bullet}(\Map(\Sigma^g,M)\times \Map(\Sigma^h,M))$ \emph{i.e.} the following diagram commutes
$$\xymatrix{HH^{(\Sigma^{g}\vee \Sigma^{h})_\bullet}_{\bullet-dim(M)}(A,A)\ar[r]^{\tilde \nabla\qquad} \ar[d]_{\simeq} & H_\bullet\big(CH^{\Sigma^g_\bullet}_\bullet(A,A) \times CH^{\Sigma^h_\bullet}_\bullet(A, A)\big) \ar[d]^{\simeq}\\
H^{\bullet-dim(M)}(\Map(\Sigma^{g}\vee \Sigma^{h},M)) \ar[r]^{\;\;\; (\rho_{in})^!\quad} &H^{\bullet}(\Map(\Sigma^g,M)\times \Map(\Sigma^h,M))
} $$ 
\end{lem}
\begin{proof}
We can assume that $M$ is equipped with a Riemannian metric and the mapping spaces $\Map(\Sigma^g,M)$ ($g\geq 0$) are equipped with a Fr\'echet manifold structure.
We have a cartesian square of fibrations
\begin{equation*}
\xymatrix{
  \Map(\Sigma^{g}\vee \Sigma^{h},M) \ar[r]^{\qquad \rho_{in}\quad\quad\quad\quad} \ar[d] 
  & \Map(\Sigma^g,M)\times \Map(\Sigma^{h},M) \ar[d]^{ev\times ev} \\
 M  \ar[r]^{\text{diagonal}\quad} & M\times M }
\end{equation*}
where the evaluation maps on the right are furthermore submersions. A tubular neighborhood $Tub(M)\subset M\times M$ of the diagonal of $M$ can be identified to the normal bundle of the diagonal. The pullback $(ev\times ev)^{-1}(Tub(M))$  by the submersion $ev\times ev:\Map(\Sigma^g,M)\times \Map(\Sigma^{h},M)\to M\times M$ can be identified with a tubular neighborhood $Tub(\Map(\Sigma^{g}\vee \Sigma^{h},M))$ of $\rho_{in}$ and thus with a normal bundle of $\rho_{in}$. The corresponding Thom spaces $M^{-TM}$ and $\Map(\Sigma^{g}\vee \Sigma^{h},M)^{-TM}$ are obtained by collapsing all the complements of the tubular neighborhhood to a point. They are disk bundles over, respectively $M$, and $\Map(\Sigma^{g}\vee \Sigma^{h},M)$. Hence, we have a diagram of pullback squares
\begin{equation*}
\xymatrix{   \Map(\Sigma^g,M)\times \Map(\Sigma^{h},M) \ar[r]^{\qquad \text{collapse}}\ar[d]_{ev\times ev} & \Map(\Sigma^{g}\vee \Sigma^{h},M)^{-TM} \ar[r]^{\pi} \ar[d]^{ev\times ev} & \Map(\Sigma^{g}\vee \Sigma^{h},M)  \ar[d]^{ev} \\
  M\times M   \ar[r]^{\text{collapse}\quad} &  M^{-TM} \ar[r]^{\pi}& M}
\end{equation*}
where the vertical arrows are  fibrations. In particular, the Thom class of $\rho_{in}$  is the pullback $(ev\times ev^*)(th(M))\in H^{dim(M)}(\Map(\Sigma^{g}\vee \Sigma^{h},M)^{-TM})$ of the Thom class $th(M)\in H^{dim(M)}(M^{-TM})$ of $M\to M\times M$. 
Since the Gysin map $(\rho_{in})^!$ is the composition $(\text{collapse})^*\circ (-\cap ev^*(th(M))) \circ \pi^* $,  it follows from Lemma~\ref{L:semifree}, and the discussion above,  that the Gysin map $(\rho_{in})^{!}$ can be modeled by the tensor product
 \begin{multline*} d^{!}\!\mathop{\otimes}_{A\otimes A} id : A \mathop{\otimes}_{A\otimes A}\big( CH_{\bullet}^{\Sigma^g_\bullet}(A,A)\otimes  CH_{\bullet}^{\Sigma^h_\bullet}(A,A)\big)\cong CH_{\bullet}^{(\Sigma^{g}\vee \Sigma^{h})_\bullet}\!\!(A,A) \\
\longrightarrow CH_{\bullet}^{\Sigma^g_\bullet}(A,A)\otimes CH_{\bullet}^{\Sigma^h_\bullet}(A,A) 
\end{multline*} 
where the $A\otimes A$-bimodule structure on $A$ is given by the multiplication and $d^!:A\to A\otimes A$ is a model for the Gysin morphism. 
Since $\Xi:A\to A^*$ is a model for the Poincar\'e duality isomorphism, we can choose the composition
$\nabla=d^! $ where $\nabla$ is defined in~\eqref{eq:nabla}.
\end{proof}
For the case of positive genus, Theorem~\ref{T:surface=cup} follows from Lemma~\ref{L:Aduality}, Lemma~\ref{L:surfaceduality} and Proposition~\ref{P:dualsurface=cup}. 
\subsubsection{Genus zero surfaces}
Now, if one of the surfaces has genus zero, we need to modify the previous arguments. First we need to define the dual of the cup product $\HHS{0}(A,A)\otimes \HHS{g}(A,A)\to \HHS{g}(A,A)$.

We denote by $a_{00}\otimes (a_{ij})_{i,j=1\dots k}$ an homogeneous element in $CH^{\Sigma^0_k}_\bullet(A,A)\cong A^{\otimes (k^2+1)}$. There is also a decomposition $CH^{\Sigma^g_k}_\bullet(A,A)\cong A^{\otimes (k^2+2k+3)}\otimes B_k^g(A)$ where $A^{\otimes (k^2+2k+3)}$ are the tensors corresponding to the top left square in the simplicial set $\Sigma^g_\bullet$ (without the bottom and right open edges, see diagram~\eqref{E:sigmag}) and $B_k^g(A)$ is the tensor power of other factors. Let $(a_{ij})_{i,j=0\dots k}$ be a generic (homogeneous) element in $A^{\otimes (k^2+2k+3)}$, and let $(b_s)$ be an element in $B_k^g(A)$. Note that there is an obvious isomorphism of vector spaces $ CH_\bullet^{(\Sigma^0\vee \Sigma^g)_k}(A,A)\cong CH_\bullet^{\Sigma^0_k}(A,A)\otimes_{A} CH^{ \Sigma^g_k}_\bullet(A,A)$ where $A$ acts on the module factors $A\cong (s_0)^k(CH_\bullet^{\Sigma^\bullet_0}(A,A))$ of $CH_\bullet^{\Sigma^\bullet_k}(A,A)$, \emph{i.e.}, the action is induced by the canonical map $A\hookrightarrow CH^{pt_\bullet}_\bullet(A,A)\to CH^{\Sigma^g_\bullet}_\bullet(A,A)$. 

 Let $pch_p: CH^{ \Sigma^g_k}_\bullet(A,A) \to CH_\bullet^{\Sigma^0_k}(A,A)\otimes_{A} CH^{ \Sigma^g_k}_\bullet(A,A)$, be the map given, for $(a_{ij})\in A^{\otimes (k^2+2k+3)}$, $(b_s)\in B^g_k(A)$, by 
\begin{multline} \label{eq:dual0cup}
 pch_p((a_{ij})\otimes (b_s)) =  \big(\hspace{-0.4cm} \mathop{\prod}_{\scriptsize \begin{array}{ll}i=0, j\leq p\\j=0, i\leq p \end{array} }\hspace{-0.4cm} a_{ij}\,  \otimes (a_{ij})_{i,j\leq p}\otimes (1) \big)\otimes_{A} \left( \big((1) \otimes (a_{i,j})_{ i \text{ or } j>p}\big) \otimes (b_s)\right)
\end{multline} where $(1)$ stands for the tensor products $1\otimes 1 \otimes \cdots$. The formula is displayed for genus $2$ in the following diagram:
\[
\begin{pspicture}(0,0.5)(10,5.5)

\psline(.5,1.5)(3.5,1.5) \psline(.5,1.5)(.5,4.5)
\psline(.5,4.5)(3.5,4.5) \psline(3.5,4.5)(3.5,1.5)
\psline(2,1.5)(2,4.5) \psline(.5,3)(3.5,3)
\psline(.5,1.5)(3.5,4.5)
\rput(1.2,3.7){$(a_{i,j})$} \rput(2.6,2.2){$(b_s)$}
 
\rput(4.2,3){$\stackrel{{\text{pch}_p}}\longmapsto $}

\psline(6.5,.5)(9.5,.5) \psline(6.5,.5)(6.5,3.5)
\psline(6.5,3.5)(9.5,3.5) \psline(9.5,3.5)(9.5,.5)
\psline(8,.5)(8,3.5)     \psline(6.5,2)(9.5,2)
\psline(6.5,.5)(9.5,3.5)
\psline(6.5,3.5)(5,3.5) \psline(5,3.5)(5,5)
\psline(5,5)(6.5,5) \psline(6.5,5)(6.5,3.5)
\psline(5,4.1)(5.8,4.1) \psline(5.8,4.1)(5.8,5)
\psline(6.5,2.6)(7.4,2.6) \psline(7.4,2.6)(7.4,3.5)
\rput(5.42,4.5){$(a_{\scriptsize{i,j}})$}
\rput(6,3.8){$(1)$} \rput(6.8,3.1){$(1)$} \rput(7.5,2.3){$(a_{i,j})$}\rput(8.4,1.4){$(b_s)$}
\rput(5,5.2){$0$} \rput(5.8,5.2){$p$} \rput(6.5,5.2){$k$}

\end{pspicture} 
\]
\begin{defn}
We define $(\pinch{0,g})_*:CH^{ \Sigma^g_k}_\bullet(A,A) \to CH_\bullet^{\Sigma^0_k}(A,A)\otimes_{A} CH^{ \Sigma^g_k}_\bullet(A,A)$ to be the map $\sum_{p=0}^k AW_{(1)}^p\otimes_A AW_{(2)}^p \circ pch_p$.
\end{defn}
Roughly speaking, the morphism $(\pinch{0,g})_*$ conists of removing the first $p^2$ tensors in $(a_{i,j})$ and tensoring them with the result of applying the second component of the Alexander-Whitney map $AW_{(2)}$ to $CH^{ \Sigma^g_\bullet}_\bullet(A,A)$, where the removed tensors have been replaced by $1$s.  

\begin{lem}\label{L:pinch0cup}
The map  $(\pinch{0,g})_*:CH^{ \Sigma^g_k}_\bullet(A,A) \to CH_\bullet^{\Sigma^0_k}(A,A)\otimes_{A} CH^{ \Sigma^g_k}_\bullet(A,A)$ is a chain map of complexes. Furthermore,
the composition 
\begin{multline*}
\Delta^{0,g}: HH^{\Sigma^{g}_\bullet}_{\bullet+dim(M)}(A,A) \stackrel{(\pinch{0,g})_*}{\longrightarrow} H^\bullet( CH_\bullet^{\Sigma^0_\bullet}(A,A)\otimes_{A} CH^{ \Sigma^g_\bullet}_\bullet(A,A)) \\
\qquad \stackrel{\tilde\nabla}{\longrightarrow}  H_{\bullet}\big(CH^{\Sigma^0_\bullet}_\bullet(A,A) \otimes CH^{\Sigma^g_\bullet}_\bullet(A, A)\big) \\ \cong HH^{\Sigma^0_\bullet}_\bullet(A,A) \otimes HH^{\Sigma^g_\bullet}_\bullet(A, A)
\end{multline*} is transfered to the cup product $ \HHS{0}(A,A)\otimes   \HHS{g}(A,A)\stackrel{\cup}\to \HHS{g}(A,A)$ by the duality isomorphism $\Theta$. 
\end{lem}
\begin{proof}
The compatibility with the differential follows from an argument similar to the proofof Lemma~\ref{L:cupbimodule}. As for Lemma~\ref{L:Aduality}, the result now follows from the definition of $\Theta$ since  $\Delta$ is the dual (through the duality isomorphism $\Xi$) of the multiplication $A\otimes A \stackrel{\cdot} \to A$.
\end{proof}

\begin{lem}\label{L:pinch0model}
 The map $(\pinch{0,g})_*:CH^{ \Sigma^g_k}_\bullet(A,A) \to CH_\bullet^{\Sigma^0_k}(A,A)\otimes_{A} CH^{ \Sigma^g_k}_\bullet(A,A)$ is a cochain model of  $\Map(\Sigma^0\vee \Sigma^g,M) \to \Map(\Sigma^g,M)$.
\end{lem}

\begin{proof}
Consider the following commutative diagram
\begin{equation}\xymatrix{ & M^{S^2\vee \Sigma^g} \ar[ld]_{\pinch{0,g}}\ar[dd]|\hole \ar[rr]^{p} && M^{I^2}\times_M M^{I^2} \ar[dd] \ar[ld]^{c}\\
  M^{\Sigma^g} \ar[rr]^{\tilde p}\ar[dd] && M^{I^2} \ar[dd] &\\
  &  \Map({\bigvee_{i=1}^{2g} S^1},M) \ar@{=}[ld] \ar[rr]|\hole && M^{S^1}\times_M M^{S^1} \ar[ld]^{\tilde c}\\  \Map({\bigvee_{i=1}^{2g} S^1},M) \ar[rr] && M^{S^1} &}\label{eq:cube}\end{equation}
where the left vertical arrows are induced by the inclusion into $\Sigma^g$ of the boundary $\partial \Sigma^g\cong \bigvee_{i=1}^{2g} S^1$ of the $4g$-gon defining $\Sigma^g$, the map $\tilde p$ by projection of a square onto $\Sigma^g$ identifying the boundary of the square with the boundary of the $4g$-gon as in diagram~\ref{E:sigma0h}, $p$ is the product of $\tilde p$ and the collapsing of $\partial I$ to a point,   $\tilde c$ is   the composition  of loops, $c$ is induced  by pinching a square (in the middle of each edge), and the right vertical arrows are induced by the inclusion of $S^1\cong \partial I^2$ into $I^2$. 

Note that the top face, front and back face of the cube are pullback diagrams. The idea is to  find a semi-free model for $c$, which, by pull back along a model for $\Map({\bigvee_{i=1}^{2g} S^1},M)\to M^{S^1}$, gives rise to a model for $\pinch{0,g}$ that coincides with $CH^{\Sigma^{g}_\bullet}_{\bullet}(A,A) \stackrel{(\pinch{0,g})_*} \longrightarrow   CH^{(\Sigma^{0}\vee \Sigma^{g})_\bullet}_{\bullet}(A,A) $.

Recall from Example~\ref{E:wedge} that the point $pt$ has a simplicial model $pt_\bullet$ which is the constant simplicial set $pt_k=pt$. Then $CH^{pt_k}_{\bullet}(A,A) \cong A$ with constant simplicial structure. Using the invariance of the Hochschild chain complex under quasi-isomorphisms of simplicial sets $X_\bullet \to Y_\bullet$~\cite{P}, it follows that $CH^{X_\bullet}_\bullet(A,A)$ is an $A$-semi-free model for $M^{I^2}$ for any pointed simplicial set $X_\bullet$ whose realization is $I^2$. The simplicial set $\Sigma^g_\bullet$
is, by definition~\eqref{E:sigmag}, defined as a quotient of a simplicial set model for  $I^2$  that we denote by $(I^2_g)_\bullet$. 
That is $\Ig$ is obtained by gluing $g^2$ squares, where the off diagonal squares are subdivided into triangles. The boundary $\partial \Ig$ is a simplicial set realizing the circle $S^1$. 
By Proposition~\ref{qi},  $CH^{\partial \Ig}_\bullet(A,A)$ and  $CH^{\partial I^2_\bullet \vee \partial \Ig}_\bullet(A,A)$ are  $CH^{\partial \Ig}_\bullet(A,A)$-semi free model of $M^{S^1}$ and $M^{S^1\vee S^1}$, respectively, and the inclusion of pointed simplicial sets $\partial \Ig\hookrightarrow \Ig $, $ \partial I^2_\bullet \vee \partial \Ig\hookrightarrow  I^2_\bullet \vee  \Ig$ induce semi-free models for the right vertical maps by functoriality of Hochschild chains.

Similarly to $\Sigma^g_k$, there is a decomposition $CH^{(I^2_g)_k}_\bullet(A,A)\cong A^{\otimes (k^2+2k+3)}\otimes \tilde{B}_k^g(A)$ where $A^{\otimes (k^2+2k+3)}$ are the tensor factors corresponding to the top left square in $\Sigma^g_\bullet$ (without the bottom and right open edges) and $\tilde B_k^g(A)$ is the tensor power of other factors. We write $(b_s)$ for an element in $\tilde B_k^g(A)$. Clearly, this decomposition restricts to $CH^{\partial(I^2_g)_k}_\bullet(A,A)\cong A^{\otimes (2k+3)}\otimes \tilde{\tilde B}_k^g(A)$. 
Let $\rho_{c}: CH^{ (I^2_g)_k}_\bullet(A,A) \to CH^{ I^2_k \vee  (I^2_g)_k}_\bullet(A,A)$ be the map given, for $(a_{ij})\in A^{\otimes (k^2+2k+3)}$, $(b_s)\in B^g_k(A)$, by
\begin{multline}\label{eq:rhoc}
 \rho_{c}((a_{ij})\otimes (b_s)) = \sum_{p=0}^k \left((a_{ij})_{i,j\leq p}\otimes (1)\right)\otimes \left( \big((1) \otimes (a_{i,j})_{ i \text{ or } j>p}\big) \otimes (b_s)\right).
\end{multline}where $(1)$ stands for the tensor products $1\otimes 1 \otimes \cdots$. We also define a linear map $\rho_{\tilde c}: CH^{ \partial (I^2_g)_k}_\bullet(A,A) \to CH^{\partial I^2_k \vee \partial (I^2_g)_k}_\bullet(A,A)$ by the same formula, but restricted to the tensors lying in $CH^{\partial(I^2_g)_k}_\bullet(A,A)\cong A^{\otimes (2k+3)}\otimes \tilde{\tilde B}_k^g(A)$. Note that this formula is indeed very close to formula~\eqref{eq:dual0cup} and can be described by a similar diagram. 

Since $ CH^{pt_\bullet}_\bullet(A,A)\to CH^{ \partial (I^2_g)_\bullet}_\bullet(A,A)$ is an $A$-semi-free quasi-isomorphism, the commutative diagram 
$$\xymatrix{CH^{pt_k}_\bullet(A,A) \ar[d]_{\text{qis}} \ar@{=}[r] & CH^{pt_k}_\bullet(A,A) \ar[d]^{\text{qis}} \\
CH^{ (I^2_g)_k}_\bullet(A,A) \ar[r]^{\rho_c} & CH^{ I^2_k \vee  (I^2_g)_k}_\bullet(A,A)} $$  implies that 
   $\rho_c$ is a cochain model for $M^{I^2}\times_M M^{I^2} \to M^{I^2}$. Note also that there are simplicial morphisms $\partial (I^2_g)_\bullet \to S^1_\bullet$, and $\partial I^2_\bullet \to S^1_\bullet$ obtained by collapsing all edges but the top left one to the basepoint. Recall that $CH^{S^1_k}_\bullet(A,A)\cong A\otimes A^{k}$. A straightforward computation shows that
the following square  $$\xymatrix{
CH^{ \partial (I^2_g)_k}_\bullet(A,A) \ar[d] \ar[r]^{\rho_{\tilde c}} & CH^{ \partial I^2_k \vee  \partial (I^2_g)_k}_\bullet(A,A) \ar[d]\\
CH^{S^1_\bullet}_\bullet(A,A)  \ar[r]^{\Phi\qquad\qquad} & CH^{S^1_\bullet}_\bullet(A,A)\otimes_A CH^{S^1_\bullet}_\bullet(A,A)} $$ is commutative, where $\Phi$ is given by
$$\Phi(a_0\otimes \cdots \otimes a_k)=\sum_{i=0}^p (a_0\otimes a_1\otimes \cdots a_i) \otimes_A (1\otimes a_{i+1}\otimes \cdots \otimes a_k).$$ 
Thus, by~\cite[Lemma 2]{FT}, $\rho_{\tilde c}$ is a cochain model for $M^{S^1\vee S^1}\to M^{S^1}$. 
Hence, $\rho_c$ and $\rho_{\tilde c}$ give a $CH^{\partial (I^2_g)_\bullet}_{\bullet}(A,A)$-semi free cochain model for the right face. 
 It follows 
that \begin{multline}  CH^{  (I^2_g)_\bullet}_\bullet(A,A)\otimes_{CH^{\partial (I^2_g)_\bullet}_{\bullet}(A,A)} CH^{\bigvee_{i=1}^{2g} S^1_\bullet}_\bullet(A,A)\\ \stackrel{\rho_c \otimes id}{\longrightarrow} CH^{ I^2_\bullet \vee (I^2_g)_\bullet}_\bullet(A,A)\otimes_{CH^{\partial (I^2_g)_\bullet}_{\bullet}(A,A)} CH^{\bigvee_{i=1}^{2g} S^1_\bullet}_\bullet(A,A) \end{multline}
is a cochain model for $ M^{S^2\vee\Sigma^g}\to M^{\Sigma^g}$. Note that, by Corollary~\ref{C:pushoutalgebra},  there are isomorphisms of chain complexes $$CH^{  (I^2_g)_\bullet}_\bullet(A,A)\otimes_{CH^{\partial (I^2_g)_\bullet}_{\bullet}(A,A)} CH^{\bigvee_{i=1}^{2g} S^1_\bullet}_\bullet(A,A) \cong CH^{\Sigma^g_\bullet}_\bullet(A,A), $$
$$ CH^{ I^2_\bullet \vee (I^2_g)_\bullet}_\bullet(A,A)\otimes_{CH^{\partial (I^2_g)_\bullet}_{\bullet}(A,A)} CH^{\bigvee_{i=1}^{2g} S^1_\bullet}_\bullet(A,A) \cong CH^{\Sigma^0_k}(A,A)\otimes_{A} CH^{ \Sigma^g_k}_\bullet(A,A).$$
 Under these isomorphisms, it is straightforward to check that $\rho_c$ transfers to $(\pinch{0,g})_*$.
\end{proof}

\subsubsection{End of the proof of Theorem~\ref{T:surface=cup}}
For $g,h>0$, Proposition~\ref{P:dualsurface=cup} and Lemma~\ref{L:Aduality}, Lemma~\ref{L:surfaceduality} imply the result. If either $g=0$ or $h=0$, the result follows 
from Lemma~\ref{L:pinch0model}, Lemma~\ref{L:surfaceduality}, Lemma~\ref{L:umkehrdual}, and Lemma~\ref{L:pinch0cup}. Lemma~\ref{L:umkehrdual} and Lemma~\ref{L:pinch0model} have obvious analogues for the cases $g\neq0, h=0$ and $g=h=0$, which can be proved similarly. $\hfill \Box$

\section{Surface Hochschild (co)homology of symmetric algebras}\label{S:surfacesymmetric}
In this section we compute the surface product of symmetric algebras and use it as a tool for explicit computations.

\subsection{Reduction to Hochschild complexes over a square and a wedge of circles}\label{S:derivedtensor}

To any topological space $X$ one can associate a Hochschild chain complex $CH^{S_\bullet(X)}_\bullet(A,M)$ (see Definition~\ref{Singular}).
According to~\cite[Theorem 2.4]{P}, if $f:X_\bullet\to Y_\bullet$ is a map of (pointed) simplicial sets inducing an isomorphism in homology, then the induced map $f_*: HH^{X_\bullet}_{\bullet}(A,M)\to HH^{Y_\bullet}_{\bullet}(A,M)$ is a quasi-isomorphism. In particular, the adjunction map $\eta: X_\bullet \to S_\bullet(|X_\bullet|)$ (Definition~\ref{Chen-map}) induces, for any space $X$ and any simplicial model $X_\bullet$ of $X$ (that is $|X_\bullet|\cong X$) a natural quasi-isomorphism $\eta: CH^{X_\bullet}_\bullet(A,M)\to CH^{S_\bullet(X)}_\bullet(A,M)$. It follows from this that, to any space $X$,  and any differential graded  commutative algebra  $(A,d)$ and  $A$-module $(M,d)$, one can associate a natural object $CH^X_\bullet(A,M):=CH^{S_\bullet(X)}_\bullet(A,M)$ in the derived category of chain complexes which is functorial in $X$, $A$ and $M$. Furthermore,  Proposition~\ref{P:shuffleinvariance} implies that, if $M=A$ equipped with its canonical $A$-module structure by multiplication, then $CH^X_\bullet(A,M)$ is a well-defined object in the homotopy category of differential graded commutative algebras (over a field of characteristic zero). If $X_\bullet$ is a simplicial model for $X$, then $CH^{X_\bullet}_\bullet(A,M)$ is naturally isomorphic to $CH^X_\bullet(A,M)$ in the derived category of chain complexes, respectively in the homotopy category of differential graded commutative algebras if $M=A$. For details on the rational homotopy theory for commutative differential graded algebras and their module, see~\cite{FHT, Q, S1}.

We denote by $S\mathop{\otimes}\limits^{\mathbb{L}}_R T$ the derived tensor product of differential graded modules $S$ and $T$ over a differential graded algebra $R$.

\begin{lem}\label{L:derivedtensor}
Let $(A,d)$ be a differential graded commutative algebra and $(M,d)$ an $A$-module. 
\begin{itemize}
 \item[i)]There is a natural isomorphism
\begin{equation} \label{E:derivedtensor} CH^{\Sigma^g}(A,M) \cong CH^{\bigvee_{i=1}^{2g}S^1}(A,M)\mathop{\otimes}\limits^{\mathbb{L}}_{CH^{S^1}(A,A)} CH^{I^2}(A,A)\end{equation} where the module structures are induced by the inclusion $S^1\cong \partial I^2\hookrightarrow I^2$ and the map $S^1\to \Sigma^g$ given by the boundary of the model for $\Sigma^g$.  \item[ii)] If furthermore $M=A$ with its canonical $A$-module structure, then the isomorphism~\eqref{E:derivedtensor} is an isomorphism of differential graded commutative algebras.
\end{itemize}
\end{lem}
\begin{proof}
Note that $\Sigma^g \cong \bigvee_{i=1}^{2g} S^1 \bigcup_{S^1} I^2$ where the maps $S^1\to I^2$ and $S^1\to \bigvee_{i=1}^{2g} S^1$ are given as in  Lemma~\ref{L:derivedtensor}.i).
Consider the standard simplicial model $S^1_\bullet$ for $S^1$ and the induced model for $\bigvee_{i=1}^{2g} S^1_\bullet$ (see Definition~\ref{D:wedge}). We consider a model $(I^2_g)_\bullet$ for $I^2$ obtained by taking the simplicial model for $\Sigma^g_\bullet$ (see Section~\ref{S:surfacemodel}) without identifying the boundary edges, \emph{i.e.} $(I^2_g)_\bullet$ consists of $g^2$-squares glued together along edges or vertices with the standard simplicial model $I^2_\bullet$ for the $g$-diagonal squares and off diagonal squares subdivided into two triangles (with model $\Triangle_\bullet$) identified along an edge. Then $\partial (I^2_g)_\bullet$ is a simplicial model for $S^1$ and moreover one has an isomorphism of simplicial sets $\Sigma^g_\bullet \cong \bigvee_{i=1}^{2g} S^1_\bullet \cup_{\partial (I^2_g)_\bullet} (I^2_g)_\bullet$.  Thus, by Corollary~\ref{C:pushoutalgebra}, there are natural quasi-isomorphisms 
\begin{eqnarray*}CH^{\Sigma^g_\bullet}(A,M)& \cong &  CH^{\bigvee_{i=1}^{2g}S^1_\bullet}(A,M)\mathop{\otimes}\limits^{}_{CH^{\partial (I^2_g)_\bullet}(A,A)} CH^{(I^2_g)_\bullet}(A,A) \\
&\cong & CH^{\bigvee_{i=1}^{2g}S^1_\bullet}(A,M)\mathop{\otimes}\limits^{\mathbb{L}}_{CH^{\partial (I^2_g)_\bullet}(A,A)} CH^{(I^2_g)_\bullet}(A,A)
\end{eqnarray*}where the last line follows because $CH^{(I^2_g)_\bullet}(A,A)$ is a free differential graded module over $CH^{\partial (I^2_g)_\bullet}(A,A)$. Furthermore, if $M=A$, the above isomorphism  is an isomorphism of differential graded commutative algebras. Since an homology isomorphism of simplicial sets induces a quasi-isomorphisms of  algebras, the result follows. 
 \end{proof}
Note that, since $I^2$ is contractible, given any commutative differential graded algebra $(A,d)$, we have a sequence of quasi-isomorphisms  \begin{equation}\label{eq:HKRI2}(A,d)\stackrel{\sim}\hookrightarrow CH^{pt_\bullet}_\bullet(A,A) \stackrel{\sim}{\to} CH^{I^2}_\bullet(A,A).\end{equation}

\bigskip

We now want to give a model for computing $HH^{\Sigma^g}_{\bullet}(A,M)$ when $(A,d)=(S(V),d)$ is a  free commutative differential graded algebra\footnote{Note that any differential graded commutative algebra is quasi-isomorphic to a   free commutative differential graded algebra.}. In view of the isomorphism~\eqref{E:derivedtensor}, we first compute $HH^{\bigvee_{i=1}^{2g}S^1}_{\bullet}(S(V),M)$. 

\subsection{HKR type Theorem for wedges of circles}
\begin{notation}\label{N:wedge}
 We denote $a_1,b_1,\dots,a_g, b_g$ the fundamental loops in $\bigvee_{i=1}^{2g}S^1$ (one for each of the $2g$-circles) whose homology classes are the generators of $H_1(\bigvee_{i=1}^{2g}S^1,\mathbb{Z})$ (and vanish on all but one circle).
\end{notation}
Let $S^1_\bullet$ be the standard simplicial model for the circle (see Example~\ref{S^1}). Then a model for the wedges of $2g$-circles is given by $\bigvee_{i=1}^{2g}S^1_\bullet$ which, in simplicial degree $n$, is the finite pointed set $\bigvee_{i=1}^{2g}S^1_n \cong [2gn]=\{0,1,2,\dots, 2gn\}$. In particular, $CH_\bullet^{\bigvee_{i=1}^{2g}S^1_\bullet}(A,M)\cong M\otimes \big(A^{\otimes n}\big)^{\otimes 2g}$. We write $m\otimes (x_i^j\otimes y_i^j)_{i=1\dots n}^{j=1\dots g}$ for an homogeneous   tensor  $m\otimes x_1^1\otimes y_1^1\otimes x_2^2\otimes \cdots \otimes y_n^{g}$ in $M\otimes \big(A^{\otimes n}\big)^{\otimes 2g}$.

The homology $H_\bullet(\bigvee_{i=1}^{2g}S^1)$ can be identified canonically with $k\oplus (\bigoplus_{i=1}^{g} k [a_i] \oplus k [b_i])$, where $[a_i], [b_j]$  (of homological degree 1) are the fundamental classes of the circle factors $a_i,b_j$  in the wedge $\bigvee_{i=1}^{2g}S^1$. 

\smallskip

The linear maps $V\ni v \mapsto [a_i] v $ and $V\ni v \mapsto [b_i] v$ uniquely extend to (degree -1) derivations $s_a^i:S(V) \to S(H_\bullet(\bigvee_{i=1}^{2g}S^1)\otimes V)$ and $s_b^i: S(V) \to S(H_\bullet(\bigvee_{i=1}^{2g}S^1)\otimes V)$ (with $S(V)$-module structure given by multiplication on the factor $k\otimes S(V) \cong S(H_0(\bigvee_{i=1}^{2g}S^1)\otimes V)$). We also extend  $s_a^i$ and $s_b^i$ as derivations of $S(H_\bullet(\bigvee_{i=1}^{2g}S^1)\otimes V)$ by setting $s_a^i(H_1(\bigvee_{i=1}^{2g}S^1)\otimes V)=0$ and $s_b^i(H_1(\bigvee_{i=1}^{2g}S^1)\otimes V)=0$. It follows, since $s_a^i$, $s_b^j$ are degree -1 derivations, that $(s_a^i)^2=0$ and $(s_b^j)^2=0$.
Similarly a differential $d$ on $S(V)$ naturally extends to a differential $d^\vee$ on $S(H_\bullet(\bigvee_{i=1}^{2g}S^1)\otimes V)$ by the formula $d^\vee(v)=d(v)$,  $d^\vee([a_i] v)= -s_a^i(d(v))$ and $d^\vee([b_j]v)= -s_b^j(d(v))$. 

\smallskip

Let $\piwedge:CH^{\bigvee_{i=1}^{2g}S^1_\bullet}_\bullet(S(V),M)\to M\otimes_{S(V)} S(H_\bullet(\bigvee_{i=1}^{2g}S^1)\otimes V)$  be the map, which for $m\otimes (x_i^j\otimes y_i^j)\in CH_\bullet^{\bigvee_{i=1}^{2g}S^1_n}(S(V),M)$, is given by
\begin{multline*}
 \piwedge\big(m\otimes (x_i^j\otimes y_i^j)\big) = \hspace{-0.4cm}
\sum_{p_1+q_1+\cdots+p_g+q_g=n}\hspace{-0.4cm}\pm m   \cdot \prod_{i=1}^{g} \frac{1}{p_i !} x_{1}^i\cdots\\ \cdots x_{p_1+\cdots+p_{i-1}}^i\cdot \big(s_a^i(x^i_{p_1+\cdots +p_{i-1}+1})\cdots \cdots s_a^i(x^i_{p_1+\cdots +p_{i}})\big) \cdot x_{p_1+\cdots +p_i+1}^i\cdots\\\cdots  x_{n}^i \cdot\prod_{i=1}^{g} \frac{1}{q_i !} y_{1}^i\cdots y_{q_1+\cdots +q_{i-1}}^i  \big(s_b^i(b^i_{q_1+\cdots +q_{i-1}+1})\cdots s_b^i(b^i_{q_1+\cdots +q_i})\big)\cdot y_{q_1+\cdots +q_{i}+1}^i \cdots y_{n}^i.  
\end{multline*}
\begin{remk}\label{R:piwedge}
  Iterating the Alexander-Whitney diagonal yields a quasi-isomorphism
$CH_\bullet^{\bigvee_{i=1}^{2g}S^1_\bullet}(S(V),M)\to M \mathop{\otimes}_{S(V)} \big(CH_\bullet^{S^1_\bullet}(S(V),S(V))\big)^{\otimes_{S(V)} 2g}$ where the right hand side is a tensor-product of chain complexes. Then it is easy to check that $\piwedge$ is the composition
\begin{multline*}CH_\bullet^{\bigvee_{i=1}^{2g}S^1_\bullet}(S(V),M)\to M \mathop{\otimes}_{S(V)} \big(CH_\bullet^{S^1_\bullet}(S(V),S(V))\big)^{\otimes_{S(V)} 2g} \\ \stackrel{id\otimes (\pi)^{\otimes 2g}}\longrightarrow M\mathop{\otimes}_{S(V)} S(V\oplus V[1])^{\otimes_{S(V)} 2g}\cong M\mathop{\otimes}_{S(V)} S(H_\bullet(\vee_{i=1}^{2g}S^1)\otimes V)\end{multline*} where $\pi: CH_\bullet^{S^1_\bullet}(S(V),S(V)) \to \Omega^\bullet_{S(V)}\cong S(V\oplus V[1])$ is the usual Hochschild-Kostant-Rosenberg map $x_0\otimes \cdots \otimes x_n \mapsto \cfrac{1}{n!} \; x_0 d(x_1)\cdots d(x_n)$ (here $ \Omega^\bullet_{S(V)}$ is the module of K\"ahler differentials), see~\cite{L}.
In particular there is a canonical isomorphism of differential graded algebras $\left(S\big(H_\bullet(\vee_{i=1}^{2g}S^1)\otimes V\big),d^{\vee}\right) \cong  \Omega^\bullet_{S(V)}\mathop{\otimes}\limits_{S(V)}\cdots\mathop{\otimes}\limits_{S(V)}\Omega^\bullet_{S(V)}$ (here $ \Omega^\bullet_{S(V)}$ is equipped with its usual internal differential induced by the one on $S(V)$). 
\end{remk}

There is also a morphism of  graded algebras\footnote{(but not of \emph{differential} graded algebras)} 
\begin{eqnarray}\label{eq:epsilonwedge}\varepsilon^{\bigvee_{i=1}^{2g}S^1}: S\left(H_\bullet\left(\bigvee_{i=1}^{2g}S^1\right)\otimes V\right) \to CH_\bullet^{\bigvee_{i=1}^{2g}S^1_\bullet}(S(V),S(V))\end{eqnarray} (the algebra structure on $CH_\bullet^{\bigvee_{i=1}^{2g}S^1_\bullet}(S(V),S(V))$ is given by the shuffle product) which maps an element $m\otimes ([a_i] v)$ to $m\otimes (\delta_{i,j}(v)\otimes 1)^{j=1\dots g} \in CH^{\bigvee_{i=1}^{2g}S^1_1}_\bullet(S(V),M)$  (where $\delta_{i,j}(v)=1$ if $i\neq j$ and $\delta_{i,i}(v)=v$) and  maps an element $m\otimes ([b_i] v)$ to $m\otimes (1\otimes \delta_{i,j}(v))^{j=1\dots g} \in CH^{\bigvee_{i=1}^{2g}S^1_1}_\bullet(S(V),M)$. In other words, $\varepsilon^{\bigvee_{i=1}^{2g}S^1}$ sends an element $[a_i]  v$ to the elements $1\otimes \cdots 1\otimes v\otimes 1 \cdots \otimes 1\in CH_1^{\bigvee_{i=1}^{2g}S^1_\bullet}(S(V),S(V))$ where $v$ is the tensor indexed by the circle in the wedge $\bigvee_{i=1}^{2g}S^1$ with fundamental class $[a_i]$ (and similarly for $[b_i] v$). Clearly, $\varepsilon^{\bigvee_{i=1}^{2g}S^1}$ is a morphism of $(S(V),d)$-algebras. Hence, it induces a morphism of $(S(H_\bullet(\bigvee_{i=1}^{2g}S^1)\otimes V),d^\vee)$-modules:
$$\varepsilon^{\bigvee_{i=1}^{2g}S^1}:(M\mathop{\otimes}_{S(V)} S(H_\bullet(\vee_{i=1}^{2g}S^1)\otimes V),d^\vee) \to CH^{\bigvee_{i=1}^{2g}S^1_\bullet}(S(V),M). $$
 
\begin{remk}
 The map $\varepsilon^{\bigvee_{i=1}^{2g}S^1}$ is the composition
\begin{multline*}M\mathop{\otimes}_{S(V)} S(H_\bullet(\vee_{i=1}^{2g}S^1)\otimes V) \cong M\mathop{\otimes}_{S(V)} S(V+V[1])^{\otimes_{S(V)} 2g} \\ \stackrel{id \otimes (\varepsilon)^{2g}} \longrightarrow  M \mathop{\otimes}_{S(V)} \big(CH_\bullet^{S^1_\bullet}(S(V),S(V))\big)^{\otimes_{S(V)} 2g}\to CH_\bullet^{\bigvee_{i=1}^{2g}S^1_\bullet}(S(V),M)\end{multline*} where the last map is the iterated Eilenberg-Zilber map and $\epsilon: S(V+V[1])\cong \Omega^\bullet(S(V)) \to CH_\bullet^{S^1_\bullet}(S(V),S(V))$ is the classical inverse of the Hochschild-Kostant-Rosenberg morphism, namely, the unique algebra morphism defined by $v [1]\mapsto 1\otimes v\in S(V)^{\otimes 2} = CH_\bullet^{S^1_1}(S(V),S(V))$, see~\cite[Section 3]{L}. 
\end{remk}
\begin{lem}\label{L:HKRwedge} Let $V$ be a graded vector space.
 \begin{enumerate} 
\item  The maps $$\piwedge:CH_\bullet^{\bigvee_{i=1}^{2g}S^1_\bullet}(S(V),M)\leftrightarrows M\otimes_{S(V)} S(H_\bullet(\bigvee_{i=1}^{2g}S^1)\otimes V):\varepsilon^{\bigvee_{i=1}^{2g}S^1}$$ are quasi-isomorphisms (of algebras if $M=S(V)$).
\item  Let $(S(V),d)$ be a differential graded free commutative algebra.  The map $\piwedge:CH_\bullet^{\bigvee_{i=1}^{2g}S^1_\bullet}(S(V),M)\to M\otimes_{S(V)} \big(S(H_\bullet(\bigvee_{i=1}^{2g}S^1)\otimes V), d\big)$ is a quasi-isomorphism of differential graded algebras.
\item $\piwedge\circ \varepsilon^{\bigvee_{i=1}^{2g}S^1}= id$.
\end{enumerate}
\end{lem}
\begin{proof}
By the same argument as in the proof of Lemma~\ref{L:derivedtensor},there are isomorphisms of $S(V)$-modules
\begin{eqnarray*}CH_\bullet^{\bigvee_{i=1}^{2g}S^1_\bullet}(S(V),M) &\cong & M \mathop{\otimes}\limits_{S(V)}  CH_\bullet^{S^1_\bullet}(S(V),S(V)) \mathop{\otimes}\limits_{S(V)}
 \cdots\mathop{\otimes}\limits_{S(V)} CH_\bullet^{S^1_\bullet}(S(V),S(V))\\&\cong & M \mathop{\otimes}\limits_{S(V)}^{\mathbb{L}}  CH_\bullet^{S^1_\bullet}(S(V),S(V)) \mathop{\otimes}\limits_{S(V)}^{\mathbb{L}}
 \cdots\mathop{\otimes}\limits_{S(V)}^{\mathbb{L}} CH_\bullet^{S^1_\bullet}(S(V),S(V)).
\end{eqnarray*}
Thus, according to Remark~\ref{R:piwedge}, the map $\piwedge$ is identified with
\begin{multline*}
 M \mathop{\otimes}\limits_{S(V)}^{\mathbb{L}}  CH_\bullet^{S^1_\bullet}(S(V),S(V)) \mathop{\otimes}\limits_{S(V)}^{\mathbb{L}}
 \cdots\mathop{\otimes}\limits_{S(V)}^{\mathbb{L}} CH_\bullet^{S^1_\bullet}(S(V),S(V)) \\
\stackrel{id\mathop{\otimes}\limits^\mathbb{L} {\pi^{\mathop{\otimes}\limits_{S(V)}^{\mathbb{L}} \!{2g}} } }\longrightarrow M\mathop{\otimes}\limits_{S(V)} S(V+V[1])^{\mathop{\otimes}\limits_{S(V)} 2g}\cong \left(M\mathop{\otimes}\limits_{S(V)} S(H_\bullet(\vee_{i=1}^{2g}S^1)\otimes V),d^\vee\right),
\end{multline*}
 where $\pi: CH_\bullet^{S^1_\bullet}(S(V),S(V)) \to \Omega^\bullet(S(V))\cong S(V\oplus V[1])$ is the usual Hochschild-Kostant-Rosenberg map $x_0\otimes \cdots x_n \mapsto 1/n! \; x_0 d(x_1)\cdots d(x_n)$. Since $\pi$ is a quasi-isomorphisms of algebras, (2) and the  first part of (1) follows.  Since $\piwedge$ and $\varepsilon^{\bigvee_{i=1}^{2g}S^1}$ are maps of algebras, it is sufficient to prove (3) for elements of the form $[a_i] v$, $[b_j]v$ for which the result holds trivially. 
We now prove the last part of the claim (1). By construction, $ \varepsilon^{\bigvee_{i=1}^{2g}S^1}$ is a morphism of algebras. Again, it is sufficient to check that $\varepsilon^{\bigvee_{i=1}^{2g}S^1}$ takes value in cocycles for elements of the form $[a_i] v$, $[b_j]v$, for which the result is straightforward. Thus $\varepsilon^{\bigvee_{i=1}^{2g}S^1}$ is a chain map and (3) and the first part of (1) imply that it this a quasi-isomorphism. 
\end{proof}

\begin{remk}
There is an obvious generalisation of Lemma~\ref{L:HKRwedge} for arbitrary wedge $\bigvee_{i=1}^k S^1$. For instance, there is a natural quasi-isomorphisms
$$\pi^{\bigvee _{i=1}^k S^1}: CH_\bullet^{\bigvee_{i=1}^{2g}S^1_\bullet}(S(V),S(V))\rightarrow \big( S(H_\bullet(\bigvee_{i=1}^{k}S^1)\otimes V), d^\vee\big)$$ of differential graded algebras. All statements and proofs for arbitrary wedges of circles are similar to those of even wedges $\bigvee_{i=1}^2g S^1$ and left to the reader.
\end{remk}

Let $pinch:S^1\to S^1\vee S^1$ be the standard ($k$-times iterated) pinching map. By Example~\ref{E:edgewise}, the edegewise subdivision $sd_2(S^1_\bullet)$ is the simplicial set $sd_2(S^1_n)=[2n+1]$ and it satisfies $$sd_2(S^1_\bullet)= \big(I_\bullet\mathop{\cup}\limits_{pt_\bullet} I_\bullet\big)\mathop{\cup}\limits_{pt_\bullet \coprod pt_\bullet} pt_\bullet$$ where the wedge $I_\bullet\cup_{pt_\bullet} I_\bullet$ is with respect to the maps $t$ and $s$ respectively (see Example~\ref{E:wedge}). Identifying the two $0$-simplices of $sd_2(S^1_\bullet)$ yields a simplicial map $\widetilde{pinch}_\bullet: sd_2(S^1_\bullet)\to S^1_\bullet \vee S^1_\bullet$. Explicitly, for any $n$, one has $\widetilde{pinch}_n(a(n+1)+i)=a(n)+i$ if $1\leq i\leq n$, $a=0,1$ and $\widetilde{pinch}_n(a(n+1))=0$. 
\begin{lem}\label{L:pinchS^1}
 The following diagram is commutative
$$\xymatrix{CH^{S^1}_{\bullet}(A,M) \ar[rr]^{pinch_*}  & & CH^{S^1\vee S^1}_{\bullet}(A,M) \\ 
CH^{S^1_\bullet}_{\bullet}(A,M) \ar[u]^{\eta} \ar[r]^{\hspace{-0.1cm}\mathcal{D}_\bullet(2)} & CH^{sd_2(S^1_\bullet)}_{\bullet}(A,M) \ar[r]^{\widetilde{pinch}_*} & CH^{S^1_\bullet\vee S^1_\bullet}_{\bullet}(A,M) \ar[u]^{\eta} }.$$
\end{lem}
\begin{proof}
 Note that $|sd_2(S^1_\bullet|)\cong (I\vee I)_{/\sim}$ where $\sim$ identifies the two (non glued) boundary points $(0,1)$ and $(1,1)$ of $I\vee I$. Then $|\widetilde{pinch}_\bullet|:|sd_2(S^1_\bullet|)\to S^1\vee S^1$ is the map identifying all boundary points of each intervall $I$ in $(I\vee I)_{/\sim}$. Thus $pinch:S^1\to S^1\vee S^1$ is the composition $$S^1\cong|S^1_\bullet|\stackrel{D^{-1}} \longrightarrow  |sd_2(S^1_\bullet)|\stackrel{|\widetilde{pinch}_\bullet|} \longrightarrow  |S^1_\bullet\vee S^1_\bullet|\cong S^1\vee S^1.$$ Now the result follows by naturality of $\eta$ and the fact that $\mathcal{D}_\bullet(2)$ realizes $D^{-1}$, see~\cite[Proposition 3.4]{McC}.
\end{proof}

Let $c_1,\dots, c_k$ be fundamental loops in $S^1$, i.e. $[c_i]=[S^1]\in H_1(S^1)$, and $f:S^1\to \bigvee_{i=1}^k S^1$ be the map obtained by glueing the paths $c_1,c_2,\dots, c_k$ in this order. The map $f$ induces a map $f_*: CH^{S^1}_{\bullet}(A,M) \to  CH^{\bigvee_{i=1}^k S^1}_{\bullet}(A,M)$ in the derived category of chain complexes. We identify $S(V\oplus V[1])\cong S(H_\bullet(S^1)\otimes V)$.
\begin{lem}\label{L:S^1towedge}
 Let $(S(V),d)$ be a differential graded free commutative algebra and $M$ an $(S(V),d)$-module. The following diagram is commutative in the derived category (respectively in the homotopy category of CDG algebras if $M=S(V)$) $$\xymatrix{CH^{S^1}_{\bullet}(A,M) \ar[r]^{f_*} \ar[d]_{\pi}^{\sim} & \ar[d]^{\pi^{\bigvee^{k}_{i=1}S^1}}_{\sim} CH^{\bigvee_{i=1}^k S^1}_{\bullet}(A,M) \\ \big(S(V\oplus V[1],d) \big) \ar[r]^{\hspace{-1cm}\tilde{f}}& \big(M\mathop{\otimes}\limits_{S(V)} S(H_\bullet(\bigvee_{i=1}^k)\otimes V),d^\vee\big)} $$ where $\tilde{f}$ is the unique map of $S(V)$-algebras given, for any $v\in V$, by $$[S^1]v\mapsto c_1v+c_2v\cdots +c_kv.$$
\end{lem}
\begin{proof} 
By functoriality and homotopy invariance of $\eta$ and of the Hochschild chain complex with respect to simplicial sets, there is a natural commutative diagram $$\xymatrix{CH^{S^1}_{\bullet}(A,M) \ar[rr]^{f_*} \ar[d]^{\sim} & &\ar[d]_{\sim} CH^{\bigvee_{i=1}^k S^1}_{\bullet}(A,M) \\ CH^{S^1}_{\bullet}(A,M) \ar[r]^{\hspace{-0.2cm}pinch^k}  &  CH^{\bigvee_{i=1}^k S^1}_{\bullet}(A,M)\ar[r]^{\bigvee\limits_{i=1}^k id}  & CH^{\bigvee_{i=1}^k S^1}_{\bullet}(A,M)} $$ where $pinch^k:S^1\to \bigvee_{i=1}^k S^1$ is the $(k-1)$-times iterated pinching map $$S^1\stackrel{pinch} \to S^1\vee S^1\stackrel{id\vee pinch}\to S^1\vee S^1\vee S^1 \dots \stackrel{id\vee pinch}\to \bigvee_{i=1}^k S^1.$$ Since $CH^{S^1_\bullet}_{\bullet}(A,M)\cong M\otimes_{A} CH^{S^1_\bullet}_{\bullet}(A,A)$, it is sufficient to prove the result for $M=S(V)$. 

Note that there is a natural isomorphism (in the derived category) of differential graded algebras $$CH^{\bigvee_{i=1}^k S^1_\bullet}_{\bullet}(S(V),S(V))\cong CH^{S^1_\bullet}_{\bullet}(S(V),S(V))\otimes_{S(V)} \cdots \otimes_{S(V)} CH^{S^1_\bullet}_{\bullet}(S(V),S(V))$$ by Corollary~\ref{C:pushoutalgebra}. Hence, by homotopy assocaitivity of $pinch$,  it is sufficient to prove the result for $k=2$. 

By Lemma~\ref{L:pinchS^1}, the result follows once we proved that  the following diagram.
\begin{equation}\label{E:S^1towedge}\xymatrix{CH^{S^1_\bullet}_{\bullet}(S(V),S(V)) \ar[d]^{\pi} \ar[r]^{\hspace{-0.1cm}\mathcal{D}_\bullet(2)} & CH^{sd_2(S^1_\bullet)}_{\bullet}(S(V),S(V))  \ar[r]^{\widetilde{pinch}_*}& CH^{S^1_\bullet\vee S^1_\bullet}_{\bullet}(S(V),S(V)) \ar[d]^{\pi^{\bigvee_{i=1}^2 S^1}} \\ \big(S(V\oplus V[1],d) \big) \ar[rr]^{\hspace{-1cm}\tilde{f}}&& \big( S(H_\bullet(\bigvee_{i=1}^k)\otimes V),d^\vee\big)
} \end{equation} is commutative (up to homotopy). By Lemma~\ref{L:HKRwedge}, the vertical maps in diagram~\eqref{E:S^1towedge} are quasi-isomorphisms of algebras. Note that  $f_*:CH^{S^1}_{\bullet}(S(V),S(V))\to CH^{S^1\vee S^1}_\bullet(S(V),S(V))$ is also an algebra morphism, and that $S(V\oplus V[1])$ is free. Hence it is sufficient to prove that diagram~\eqref{E:S^1towedge} is commutative in simplicial degrees $0$ and $1$, since the generators (as an algebra) of $S(V\oplus V[1])$ lies in the subspace $\pi(CH^{S^1_{\leqslant 1}}_{\bullet}(S(V),S(V)))$. In simplicial degree $0$, all the maps in diagram~\eqref{E:S^1towedge} are the identity map. Recall (see formula~\eqref{E:D(2)}) that $$\mathcal{D}_1(2):CH^{S^1_{ 1}}_{\bullet}(A,A)\cong A^{\otimes 2} \to :CH^{sd_2(S^1_{ 1})}_{\bullet}(A,A)\cong A^{\otimes 4}$$ is given by the formula $\mathcal{D}_1(2)(x\otimes y)=\sum\limits_{(\sigma,\delta)\in \mathcal{S}(2,1)} (-1)^{\sigma}S^1_\bullet(\varepsilon_{(\sigma,\delta)})(x\otimes y)$. By definition of $\mathcal{S}(2,1)$, $\sigma$ is the identity and $\delta\in Hom_{\Delta}([0], [1])$. From identity~\eqref{E:epsilonsigmadelta} defining $\varepsilon_{(\sigma,\delta)} $, we get that  $\mathcal{D}_1(2)(x\otimes y)=x\otimes y\otimes 1\otimes 1 + x\otimes 1\otimes 1\otimes y$. Now, the commutativity of diagram~\eqref{E:S^1towedge} easily follows.
\end{proof}

\subsection{HKR type Theorem for surfaces}

For the sphere $S^2\cong \Sigma^0$, there is also a Hochschild-Kostant-Rosenberg type theorem. More precisely, given a differential graded free   commutative algebra $(S(V),d)$  and an $(S(V),d)$-module $M$, there is a natural  isomorphism 
\begin{equation}
 \label{eq:HKRsphere} \pi^{S^2}:HH^{S^2}_\bullet(S(V), M)  \stackrel{\sim}\longrightarrow H_\bullet \big(M\mathop{\otimes}\limits_{S(V)} S(H_\bullet(S^2)\otimes V), d^{S^2}\big). 
\end{equation}
Here, the graded commutative algebra $S(H_\bullet(S^2)\otimes V)$ is equipped with the differential $d^{S^2}$ which is defined as the unique degree 1 derivation satisfying 
$ d^{S^2}(v)=v$  and  $d^{S^2}(\sigma v)= s_\sigma d(v)$ where $\sigma =[S^2]\in H_2(S^2)$ is the fundamental class of $S^2$ and $s_\sigma$ is the unique degre -2 derivation defined by $s_\sigma(v)= \sigma v$ and $s_\sigma(\sigma w)=0$. Note that $S(H_\bullet(S^2)\otimes V)\cong S(V\oplus V[2])$.

\smallskip

Furthermore, if $M\cong S(V)$, $\pi^{S^2}$ is an isomorphism of algebras. See~\cite{P, G} for details. 

\smallskip

For positive genus surfaces, we have the following Hochschild-Kostant-Rosenberg type result. We write $\sigma=[\Sigma^g]\in H_2(\Sigma^g)$ for the fundamental class of $\Sigma^g$ and $[a_1],[b_1],\dots,[a_g],[b_g]$ for the generators of $H_1(\Sigma^g)$. The degree -1 derivations $s_a^i$ and $s_b^j$ on $ S(H_\bullet(\Sigma^g)\otimes V)$ are defined by $s_a^i(H_{\geq 1}(\Sigma^g)\otimes V)=0=s_b^j(H_{\geq 1}(\Sigma^g)\otimes V)$ and $s_a^i(v)=[a_i]v$, $s_b^j(v)=[b_j]v$ for any $v\in V$ and $i,j=1\dots g$. Similarly the degree -2 derivation $s_\sigma$ is defined by $s_\sigma(v)=\sigma v$ and $s_\sigma(H_{\geq 1}(\Sigma^g)\otimes V)=0$. The differential $ d^{\Sigma^g} $ is the unique degree 1 derivation defined by \begin{eqnarray} \label{eq:dsurface1}
     d^{\Sigma^g}([a_i]v) = -s_a^i(d(v)),  &\qquad & d^{\Sigma^g}([b_j]v)=   -s_b^j(d(v))  \\  \label{eq:dsurface2} d^{\Sigma^g}(v)= d(v), & \qquad &   d^{\Sigma^g}(\sigma v)= s_\sigma(d(v))+\sum_{i=1}^g s_a^i(s_b^i (d(v))).\end{eqnarray}
\begin{remk}
Note that the differential $d^{\Sigma^g}$ is based on the coalgebra structure of $H_\bullet(\Sigma^g)$. That is, if $x\in H_\bullet(\Sigma^g)$, then, for any $v\in V$, the differential is given by $d^{\Sigma^g}(x v)=\sum (-1)^{|x_{(1)}|+|x_{(2)}|}\, s_{x_{(1)}}\big(s_{x_{(2)}}(d(v))\big)$ where the coproduct is given by $\Delta(x)=\sum x_{(1)}\otimes x_{(2)}$, and $s_{1}=id$. 
\end{remk}

\begin{remk}
When $X$ is a space with Sullivan model $(S(V),d)$,  the cochain algebra $\big(S(H_\bullet(\Sigma^g)\otimes V),d^{\Sigma^g}\big)$ coincides with the Haefliger model~\cite{Hae} of the sections of the trivial bundle $\Sigma^g\to \Sigma^g\times X$ where one takes $H^\bullet(\Sigma^g)$ as a cochain model for $\Sigma^g$ (which is possible since surfaces are formal spaces). This model also has been  carefully described by Brown and Szczarba~\cite{BS}. Of course, the same remark holds for (wedges) of circles in place of surfaces. It would be interesting to find a proof of this algebras quasi-isomorphism using the universal property of the Haefliger model~\cite{Hae} and the functor homology techniques  introduced in~\cite{P}.
\end{remk}

\begin{thm}\label{T:HKR}
 Let $(S(V),d)$ be a differential free  graded commutative algebra and $M$ a differential  graded  $(S(V),d)$-module. \begin{enumerate}
\item There is a natural isomorphism 
$$\eps{g}: H_\bullet\big( M\otimes_{S(V)} S(H_\bullet(\Sigma^g)\otimes V), d^{\Sigma^g}\big) \stackrel{\sim}\longrightarrow HH^{\Sigma^g_\bullet}_\bullet(S(V),M),$$ which is an isomorphism of algebras if $M=S(V)$.
 \item The following diagram is commutative
$$\xymatrix{H_\bullet \big(M\mathop{\otimes}\limits_{S(V)} S(H_\bullet(\mathop{\bigvee}\limits_{i=1}^{2g}S^1)\otimes V)\big) \ar[d]_{\varepsilon^{\bigvee_{i=1}^{2g}S^1}} \ar[r]^{p} & H_\bullet \big(M\mathop{\otimes}\limits_{S(V)} S(H_\bullet(\Sigma^g)\otimes V)\big) \ar[d]_{\varepsilon^{\Sigma^g}}\ar[r]^{q}&  H_\bullet \big(M\mathop{\otimes}\limits_{S(V)} S(H_\bullet(S^2)\otimes V)\big) \ar[d]_{\varepsilon^{S^2}}\\ HH^{\bigvee_{i=1}^{2g}S^1}_\bullet(S(V),M)  \ar[r]^{(\bigvee_{i=1}^{2d} S^1 \hookrightarrow \Sigma^g)_\bullet}  & HH^{\Sigma^g_\bullet}_\bullet(S(V),M) \ar[r]^{(\Sigma^g \twoheadrightarrow S^2)_\bullet}  & HH^{S^2}_\bullet(S(V),M)} $$
where the horizontal maps $p$ and $q$ are the algebra homomorphisms, respectively induced by the homology maps $H_\bullet(\mathop{\bigvee}\limits_{i=1}^{2g}S^1)\otimes V\to H_\bullet(\Sigma^g)$, and $H_\bullet(\Sigma^g)\to H_\bullet(S^2)$, obtained by the obvious inclusion and surjection of spaces.
\end{enumerate}
\end{thm}
To prove Theorem~\ref{T:HKR}, we want to use the computation in Lemma~\ref{L:HKRwedge} and apply Lemma~\ref{L:derivedtensor}. Hence, we first need a semi-free model of $CH^{I^2}_\bullet(S(V),S(V))$ as a $CH^{S^1}_\bullet(S(V),S(V))$-module.
\begin{proof}
\begin{description}
\item[(a)] Note that $H_\bullet(S^1)\cong k[\xi]$ (with $|\xi|=-1$) and that the standard Hochschild-Kostant-Rosenberg theorem yields a natural isomorphism $$\varepsilon^{S^1}:H_\bullet \big(M\mathop{\otimes}\limits_{S(V)} S(H_\bullet(S^1)\otimes V)\big) \stackrel{\sim}\longrightarrow HH^{S^1}_\bullet(S(V),M).$$
\item[(b)]Since $I^2$ is contractible, for any (DG commutative) algebra $A$, there are natural isomorphisms $HH^{I^2}_\bullet(A,A) \cong HH^{\rm pt}_\bullet(A,A)\cong H_\bullet(A).$ Further, the canonical map $$CH^{(I^2_g)_\bullet}_\bullet(A,A) \to CH^{(I^2_g)_0}_\bullet(A,A)\cong A^{\otimes g^2} \to A,$$ where $(I^2_g)_\bullet$ is the simplicial model for the square described in the proof of Lemma~\ref{L:derivedtensor} and the right map $A^{\otimes g^2} \to A$ is induced by  the  map of pointed sets $(I^2_g)_0\to \{ 0\}$, is a quasi-isomorphism of algebras. 
\item[(c)] The algebra morphism (coming from {\bf (a)} and {\bf (b)}) $$H_\bullet(S(k[\xi]\otimes V))\cong HH^{S^1}_\bullet (S(V),S(V)) \to HH^{I^2}_\bullet(S(V),S(V)) \cong H_\bullet(S(V))$$ is induced by the unique (differential graded) commutative algebra morphism $S(k[\xi]\otimes V)\to S(V)$ satisfying $\xi \otimes v\mapsto 0$, $1\otimes v\mapsto v$ for any $v\in V$. This follows since the image of $\xi\otimes v$ in $CH^{S^1_\bullet}(S(V),S(V))$ lies in simplicial degree $1$.  
\end{description}

\smallskip

By Lemma~\ref{L:derivedtensor}, Lemma~\ref{L:HKRwedge}, and {\bf (a)}, {\bf (b)}, and {\bf (c)} above, there is a natural isomorphism 
\begin{equation}\label{eqT:HKR}
  S(H_\bullet(\mathop{\bigvee}\limits_{i=1}^{2g}S^1)\otimes V) \mathop{\otimes}\limits_{S(k[\xi]\otimes V)} K^{I^2}(S(V)) \stackrel{\sim}\longrightarrow HH^{\Sigma^g}_\bullet(S(V),M)
\end{equation}
for any $(S(k[\xi]\otimes V),d)$-semifree resolution $K^{I^2}(S(V))$ of $(S(V),d)$. Further, if $K^{I^2}(S(V))$ is also a resolution as an algebra and $M\cong S(V)$, then the isomorphism~\eqref{eqT:HKR} is an isomorphism of algebras. We now construct an explicit resolution  $K^{I^2}(S(V))$.

\medskip

We first recall a  $S(V\oplus V)$-semifree resolution of $(S(V),d)$, which we denote by $K^I(S(V))$. (Note that $S(V\oplus V)\hookrightarrow CH^{\partial I_\bullet}_\bullet(S(V),S(V))$ and $S(V) \hookrightarrow CH^{ I_\bullet}_\bullet(S(V),S(V))$ are quasi-isomorphisms.) We identify $S(V\oplus V[1] \oplus V)\cong S(k[x_0]\oplus k[\xi] \oplus k[x_1]) \otimes V$ where $[x_0],[x_1]$ are of degree $0$ and $[\xi]$ is of  degree -1. Let $s$ be the unique degree -1 derivation of $S(V\oplus V[1] \oplus V)$ given by $s([x_i]v)= [\xi] v$ and $s([\xi] v)=0$. Then, by~\cite[Section 15.(c), Example 1]{FHT},  $$K^I(S(V)):=(S(V\oplus V[1] \oplus V), d^I)$$ is an $S(V\oplus V)$-semifree resolution of $(S(V),d)$ where $d^I$ is the unique degree 1 derivation defined by $d^I([x_i] v)=[x_i] d(v)$ and \begin{equation}     d^I([\xi] v)=  [x_0]v -[x_1]v -\sum_{i=1}^\infty \cfrac{(sd^I)^n}{n!} ([x_0]v).\end{equation} Hence there is an isomorphism $K^I(S(V))=(S(V\oplus V[1] \oplus V),d^I)\cong CH^{I_\bullet}(S(V),S(V))$ in the homotopy category of commutative differential graded algebras and a commutative diagram 
$$\xymatrix{ CH^{pt_\bullet}(S(V),S(V)) \ar[r]^{s_*} \ar[d]^{\sim} & CH^{I_\bullet}(S(V),S(V)) \ar[d]^{\sim}& CH^{pt_\bullet}(S(V),S(V)) \ar[l]_{t_*} \ar[d]^{\sim}\\ (S(V),d) \ar[r]^{\hspace{-0.2cm}v\mapsto [x_0]v} & (K^I(S(V),d^I) & (S(V),d) \ar[l]_{\hspace{0.3cm}v\mapsto [x_1]v}} $$ where $s,t$ are the two inclusions $pt_\bullet \to I_\bullet$ defined in Example~\ref{E:wedge}.
Note that, by Corollary~\ref{C:product}, for any differential graded algebra $(A,d)$, there is a natural isomorphism 
\begin{equation}\label{E:isoI^2} CH^{I^2_\bullet}_\bullet(A,A)\cong CH^{I_\bullet}_\bullet(CH^{I_\bullet}_\bullet(A,A), CH^{I_\bullet}_\bullet(A,A)) \end{equation}
 where $CH^{I_\bullet}_\bullet(A,A)$ is equipped with the Hochschild total differential and the algebra structure given by the shuffle product. Note that since $I_\bullet$, $I_\bullet^2$ are contractible, one can simply notice that the canonical inclusion $A\hookrightarrow CH^{I_\bullet}_\bullet(CH^{I_\bullet}_\bullet(A,A), CH^{I_\bullet}_\bullet(A,A))$  is a quasi-isomorphism instead of using  Corollary~\ref{C:product}. Thus, there is an isomorphism \begin{equation}\label{E:K^I(K^I)} K^I(K^I(S(V)))\cong  CH^{I^2_\bullet}_\bullet(S(V), S(V))\end{equation} in the homotopy category of commutative differential graded algebras. We denote $K^{I^2}(S(V)):=K^I(K^I(S(V)))$ and write $d^{I^2}$ for its differential. By construction $$K^{I^2}(S(V))\cong S(V^{\oplus 4} \oplus V[1]^{\oplus 4} \oplus V[2])\cong S\left(\big(\bigoplus_{i,j=1,2} k[x_{ij}] \oplus \bigoplus_{i=0}^1 (k[\xi_i] \oplus k[\xi'_i]) \oplus k[\sigma]\big) \otimes V\right)$$ where $|x_{ij}|=0$, $|\xi_i|=|\xi'_j|=-1$ and $|\sigma|=-2$. One may think of $x_{ij}$ as points, $\xi_i$ as a path from $x_{i0}$ to $x_{i1}$, $\xi'_j$ as a path from $x_{0j}$ to $x_{1j}$ as in the following picture:
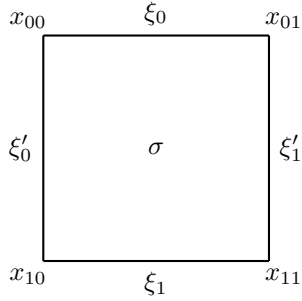
\begin{equation*}
\begin{pspicture}(0,0)(4,4)
 \psline(0.5,.5)(0.5,3.5) \psline(.5,.5)(3.5,.5) \psline(3.5,.5)(3.5,3.5) \psline(.5,3.5)(3.5,3.5)
 \rput(2,3.8){$\xi_0$} \rput(2,.2){$\xi_1$} \rput(0.2,2){$\xi'_0$}  
\rput(3.8,2){$\xi'_1$}  \rput(2,2){$\sigma$} \rput(0.3,3.7){$x_{00} $} \rput(3.7,3.7){$x_{01} $} \rput(0.3,0.3){$x_{10} $} \rput(3.7,0.3){$x_{11} $}
\end{pspicture} \end{equation*}
In particular the subalgebra $K^{\xi_0}(V):=S((k[x_{00}]\oplus k[\xi_0] \oplus k[x_{01}])\otimes V)\subset K^{I^2}(S(V))$ is a differential graded subalgebra which is canonically isomorphic to $K^I(S(V))$. The same holds for the 3 other subalgebras: $K^{\xi_1}(V):=S((k[x_{10}]\oplus k[\xi_1] \oplus k[x_{11}])\otimes V)$, $K^{\xi'_0}(V):=S((k[x_{00}]\oplus k[\xi'_0] \oplus k[x_{10}])\otimes V)$ and $K^{\xi'_1}(V):=S((k[x_{01}]\oplus k[\xi'_1] \oplus k[x_{11}])\otimes V)$.
The differential $d^{I^2}$ is the degree 1 derivation defined by $d^{I^2}([x_{ij}]v)=[x_{ij}]d(v)$, $d^{I^2}([\xi_{i}]v)=d^I([x_{i0}] v)$, $d^{I^2}([\xi'_{j}]v)=d^I([x_{0j}] v)$ and 
$$d^{I^2}([\sigma]v)=[\xi'_0] v - [\xi'_1] v -\sum_{n=1}^\infty \cfrac{(\tilde{s}d^{I^2})^n}{n!} ([\xi_0]v) $$ where $\tilde{s}$ is the degree -1 derivation defined by $\tilde{s}([x_{ij}] v)=[\xi_j] v$, $\tilde{s}([\xi_{i}] v)=[\sigma]v=\tilde{s}([\xi'_{j}] v)$ and $\tilde{s}([\sigma] v)=0$.

Since the boundary of $I^2_\bullet$ is the wedge $\big(I_{\bullet} \,{}_{t}\!\cup_{s} I_{\bullet}\big) \cup_{pt_\bullet \coprod pt_\bullet} \big(I_{\bullet} \,{}_{t}\!\cup_{s} I_{\bullet}\big)$, the natural isomorphism~\eqref{E:isoI^2} induces a natural isomorphism \begin{eqnarray}\label{E:boundarymodel} CH^{\partial I^2_\bullet}(S(V),S(V)) &\cong& K^{\partial I^2}(S(V)),\end{eqnarray} where $K^{\partial I^2}(S(V))$ denotes the differential graded commutative algebra $$K^{\partial I^2}(S(V)):=\big(K^{\xi_0}(V) \hspace{-0.1cm}\mathop{\otimes}\limits_{S([x_{01}]V)}\hspace{-0.1cm} K^{\xi'_1}(V)\big)\hspace{-0.2cm}\mathop{\otimes}_{S([x_{00}]V\oplus [x_{11}] V)} \hspace{-0.2cm} \big(K^{\xi'_0}(V)\hspace{-0.1cm} \mathop{\otimes}_{S([x_{10}]V)} \hspace{-0.1cm}K^{\xi_1}(V)\big).$$  

\medskip

We now need to identify the induced map $\partial I^2 \to \bigvee_{i=1}^{2g} S^1$. We first consider the genus 1 case. We still identify $\partial I^2$ with $(I\vee I)\cup_{\{*\}\coprod \{*\}} (I\vee I)$ where the endpoints of the two (lenght 2) intervals are identified. For a surface of  genus $g=1$, given as a quotient of $(I^2_1)$ by the path $a_1 b_1 a_1^{-1} b_1^{-1}$, the boundary map $\partial I^2 \to S^1\vee S^1$ factors as the composition
$$\partial I^2\cong (I\vee I)\cup_{\{*\}\coprod \{*\}} (I\vee I) \stackrel{\bigvee_{i=1}^4 col}\longrightarrow (S^1\vee S^1) \vee (S^1\vee S^1) \stackrel{(a_1\vee b_1)\vee (a_1\vee b_1)}\longrightarrow S^1\vee S^1$$ where the first map $ \bigvee_{i=1}^4 col$ collapses each interval to a circle.

 The map $s_*\otimes t_*: S(V)\otimes S(V)\cong CH^{pt_\bullet\coprod pt_\bullet}_{\bullet}(S(V),S(V))  \to CH^{I_\bullet}_\bullet(A,A)$ induces a quasi-isomorphism of algebras $$CH^{S^1_\bullet}_\bullet (S(V),S(V))\cong CH^{I_\bullet}_\bullet(S(V),S(V))\mathop{\otimes}_{S(V)\otimes S(V)} S(V) \cong K^I(S(V))\mathop{\otimes}_{S(V)\otimes S(V)} S(V)$$ since $CH^{I_\bullet}_\bullet(S(V),S(V))$ is a semi-free $S(V)\otimes S(V)$-algebra.
 Thus the algebra quasi-isomorphism~\eqref{E:boundarymodel} transfers the algebra homomorphism $CH^{S^1_\bullet}_\bullet (S(V),S(V))\cong CH^{\partial I^2_\bullet}_\bullet (S(V),S(V)) \to CH^{S^1\vee S^1}_\bullet(S(V),S(V))$ to the algebra homomorphism 
$$ \varphi_1: K^{\partial I^2}(S(V)) \longrightarrow S(H_\bullet(S^1\vee S^1)\otimes V)$$
defined by $\varphi_1 ([x_{ij}]v)=v$, $\varphi_1([\xi_i]v)=[a_1]v$ and  $\varphi_1([\xi'_j]v)=[b_1]v$. 
Now, since $s$, $\tilde{s}$ are degree -1 derivations, and furthermore since $s_a^i$, $s_b^j$ are degree -1 derivations with square 0, a straightforward computation gives an algebra isomorphism $$S(H_\bullet(\mathop{\bigvee}\limits_{i=1}^{2}S^1)\otimes V)\big) \mathop{\otimes}\limits_{S(k[\xi]\otimes V)} K^{I^2}(S(V)) \cong (S(H_\bullet(\Sigma^1)\otimes V),d^{\Sigma^1}). $$

\smallskip

For a surface of genus $g>1$, our model $\Sigma^g$ (see Section~\ref{S:surfacemodel}) is also obtained as a quotient of a square. Now the  the boundary map $\partial I^2 \to \bigvee_{i=1}^{2g} S^1$ factors as the composition
$$\partial I^2\cong (I\vee I)\cup_{\{*\}\coprod \{*\}} (I\vee I) \stackrel{\bigvee_{i=1}^4 col}\longrightarrow (S^1\vee S^1) \vee (S^1\vee S^1) \stackrel{(f_1\vee f_2)\vee (f_3\vee f_4)}\longrightarrow S^1\vee S^1$$ where the first map $ \bigvee_{i=1}^4 col$ still collapses each interval to a circle and the maps $f_i:S^1\to \bigvee_{i=1}^g S^1$ are the loops  given by
$$f_1= a_1b_1\dots a_{g/2}b_{g/2}, \quad  f_2=a_{g/2+1}b_{g/2+1}\dots a_g b_g, \quad f_3=b_1a_1\dots b_{g/2}a_{g/2}$$  and   $f_4=b_{g/2+1}a_{g/2+1}\dots b_g a_g$ when $g$ is even, and by $$f_1= a_1b_1\dots a_{m}, \quad  f_2=b_{m}a_{m+1}\dots a_g b_g, \quad f_3=b_1a_1\dots b_{m}$$  and   $f_4=a_{m}b_{m+1}\dots b_g a_g$ when $g=2m-1$ is odd. Now by the argument for $g=1$ and Lemma~\ref{L:S^1towedge}, it follows that 
the algebra quasi-isomorphism~\eqref{E:boundarymodel} transfers the algebra homomorphism $CH^{S^1_\bullet}_\bullet (S(V),S(V))\cong CH^{\partial I^2_\bullet}_\bullet (S(V),S(V)) \to CH^{\bigvee_{i=1}^{2g}S^1}_\bullet(S(V),S(V))$ to the algebra homomorphism 
$$ \varphi_g: K^{\partial I^2}(S(V)) \longrightarrow S(H_\bullet(\bigvee_{i=1}^{2g}S^1)\otimes V)$$
defined by $\varphi_g ([x_{ij}]v)=v$, and, for $g$ even by \begin{align*}
 \varphi_g([\xi'_0]v)=[a_1]v+[b_1]v+\cdots+[b_{g/2}]v, \quad & \quad \varphi_g([\xi_1]v)=[a_{g/2+1}]v+\cdots [a_g]v+[b_g]v,\\ 
\varphi_g([\xi_0]v)= [b_1]v +[a_1]v+\cdots [a_{g/2}] v, \quad & \quad \varphi_g([\xi'_1]v)= [b_{g/2+1}]v+\cdots +[b_g] v+[a_g] v,\end{align*}
 and, for $g=2m-1$ odd by  
 \begin{align*}
 \varphi_g([\xi'_0]v)=[a_1]v+[b_1]v+\cdots+[a_{m}]v, \quad & \quad \varphi_g([\xi_1]v)=[b_{m}]v+\cdots [a_g]v+[b_g]v,\\ 
\varphi_g([\xi_0]v)= [b_1]v +[a_1]v+\cdots [b_{m}] v, \quad & \quad \varphi_g([\xi'_1]v)= [a_{m}]v+\cdots +[b_g] v+[a_g] v.\end{align*} Again,  a straightforward computation gives an algebra isomorphism $$S(H_\bullet(\mathop{\bigvee}\limits_{i=1}^{2g}S^1)\otimes V)\big) \mathop{\otimes}\limits_{S(k[\xi]\otimes V)} K^{I^2}(S(V)) \cong (S(H_\bullet(\Sigma^g)\otimes V),d^{\Sigma^1}). $$
This proves claim (1) of Theorem~\ref{T:HKR}. The commutativity of the left square in claim (2) follows from the isomorphism~\eqref{eqT:HKR} and the computation above. Note that \begin{eqnarray*} CH^{S^2_\bullet}_\bullet (A,M) & \cong & CH^{pt_\bullet}_\bullet(A,M)\mathop{\otimes}_{CH^{\partial I^2_\bullet}_\bullet (A,A)} CH^{I^2_\bullet}_\bullet (A,A) \\ &\cong &  CH^{pt_\bullet}_\bullet(A,M)\mathop{\otimes}^{\mathbb{L}}_{CH^{\partial I^2_\bullet}_\bullet (A,A)} CH^{I^2_\bullet}_\bullet (A,A).\end{eqnarray*} Hence there  is a commutative diagram 
$$\xymatrix{CH^{\Sigma^g_\bullet}_\bullet(A,M) \ar[d]^{\sim} \ar[r]^{\Sigma^g_\bullet \to S^2_\bullet}  & CH^{S^2_\bullet}_\bullet(A,M) \ar[d]^{\sim}\\ CH^{\bigvee_{i=1}^{2g}S^1_\bullet}_\bullet(A,M) \mathop{\otimes}\limits_{CH^{\partial I^2_\bullet}_\bullet (A,A)} CH^{I^2_\bullet}_\bullet (A,A)\ar[r]^{\quad p_*\otimes id} &  CH^{pt_\bullet}_\bullet(A,M)\mathop{\otimes}\limits_{CH^{\partial I^2_\bullet}_\bullet (A,A)} CH^{I^2_\bullet}_\bullet (A,A) } $$ where $p: \bigvee_{i=1}^{2g} S^1_\bullet \to pt_\bullet$ is the canonical map. Let $\tilde{p}: K^I(S(V))\to S(V)$ be the algebra map defined by $\tilde{p}([x_i] v)=v$ and $\tilde{p}([\xi] v)=0$. This is a map of differential graded algebras and, furthermore, since the composition $S(V)\hookrightarrow K^I(S(V)) \stackrel{\tilde{p}}\to S(V)$ is the identity,   $\tilde{p}$ is a quasi-isomorphism and the following diagram is commutative
$$\xymatrix{CH^{I_\bullet}_\bullet(S(V),S(V) \ar[r]^{p_*} & CH^{pt_\bullet}_\bullet (S(V),S(V)) \\
K^I(S(V))\ar[u]^{\sim} \ar[r]^{\tilde{p}} & S(V) \ar[u]^{\sim} } $$ in the homotopy category of differential graded commutative algebras. Since $p_*:CH^{\bigvee_{i=1}^{2g}S^1_\bullet}_\bullet(S(V),S(V))\to H^{pt_\bullet}_\bullet (S(V),S(V))$ is the composition 
\begin{multline*}
 CH^{\bigvee_{i=1}^{2g}S^1_\bullet}_\bullet(S(V),S(V)) \cong CH^{S^1_\bullet}_\bullet(S(V),S(V)) \mathop{\otimes}_{S(V)} \cdots \mathop{\otimes}_{S(V)}  CH^{S^1_\bullet}_\bullet(S(V),S(V)) \\
\cong \left(CH^{I_\bullet}_\bullet(S(V),S(V))\mathop{\otimes}_{S(V\oplus V)} S(V)\right) \mathop{\otimes}_{S(V)} \cdots \mathop{\otimes}_{S(V)}  \left(CH^{I_\bullet}_\bullet (S(V),S(V))\mathop{\otimes}_{S(V\oplus V)} S(V)\right)\\
 \stackrel{{p}_*\mathop{\otimes}_{S(V\oplus V)} id}\longrightarrow  S(V)
\end{multline*}
it follows that  the isomorphism of Theorem~\ref{T:HKR} claim (1) transfers $(\Sigma^g_\bullet \to S^2_\bullet)_* : CH^{\Sigma^g_\bullet}_\bullet(S(V),S(V)) \to CH^{S^2_\bullet}_\bullet(S(V),S(V)) $ to the map \begin{multline*} (\tilde{p}\otimes id)^{\otimes 2g}\otimes id:\left((K^I(S(V))\mathop{\otimes}_{S(V\oplus V)} S(V)\right)^{\otimes 2g} \mathop{\otimes}_{K^{\partial I}(S(V))} K^{I^2}(S(V))\\ \to S(V) \mathop{\otimes}_{K^{\partial I}(S(V))} K^{I^2}(S(V))\cong S(H_\bullet(S^2)\otimes V),\end{multline*}
 which proves the commutativity of the right square in claim (2).
\end{proof}
The importance of Theorem~\ref{T:HKR} folllows from the fact that any differential graded commutative algebra is quasi-isomorphic, as an algebra, to a graded symmetric one. 

\begin{cor}\label{C:HKR}
There is a natural isomorphism $${\eps{g}}^*: HH_{\Sigma^g_\bullet}^\bullet(S(V),M)\stackrel{\sim}\longrightarrow  H_\bullet\left( Hom_{S(V)}\big( S(H_\bullet(\Sigma^g)\otimes V),S(V)\big ), {d^{\Sigma^g}}^*\right).$$ 
\end{cor}
\begin{proof} Let $X_\bullet$  be any pointed simplicial set and let $A$ be a differential graded commutative algebra. There is a map of differential graded algebras $A=CH^{pt_0}_\bullet(A,A) \to CH^{pt_\bullet}_\bullet(A,A)\to CH^{X_\bullet}_\bullet(A,A) $ where the last map is the unique pointed map. It follows that  the Hochschild complex $(CH^{X_\bullet}_\bullet(A,A),D)$ is a chain complex of semi-free $A$-modules. More explicitly,  the $A$-module structure is given by multiplication on the tensor  $A$ corresponding to the base point in $CH^{X_\bullet}_\bullet(A,A)\cong A \otimes A^{\otimes x_\bullet}$.
Now let $M$ be an $A$-module. Then there is an isomorphism of cochain complexes
$$(CH_{X_\bullet}^\bullet(A,M),D)\cong \big( Hom_A(CH^{X_\bullet}_\bullet(A,A),M), D^* \big), $$
where the differential $D^*$ is the dual of the differential $D$ on the Hochschild chain complex $CH^{X_\bullet}_\bullet(A,A)$. Since $\big(S(H_\bullet(\Sigma^g)\otimes V),d^{\Sigma^g}\big)$ is also semi-free, the result follows  from the first statement in Theorem~\ref{T:HKR}.
\end{proof}

\subsection{The surface product for Lie groups}\label{S:Liegroups}
In this section, we apply Theorem~\ref{T:HKR} and Lemma~\ref{L:Aduality} to compute the Hochschild surface product for odd spheres and Lie groups. The idea is that in both cases, the commutative differential graded algebra of the forms $\Omega^\bullet {M}$ (where $M=S^{2n+1}$ or $M=G$ is a Lie group) is  quasi-isomorphic as a  differential graded algebra to a symmetric algebra $S(V)$ with zero differential. 

\medskip

Let $A=(S(V),0)$ be a free graded commutative algebra (with zero differential).  Then the identities~\eqref{eq:dsurface1} and~\eqref{eq:dsurface2} immediately implies that the differential $d^{\Sigma^g}=0$ for any genus $g$. Similarly, the differentials $d^\vee$ and $d^{S^2}$ vanish, too. Hence, for any $S(V)$-module $M$, by Theorem~\ref{T:HKR}, there is a commutative diagram (natural in $M$ and $V$),
\begin{eqnarray}\label{eq:diagramd=0}\xymatrix{\\ M\mathop{\otimes}\limits_{S(V)} S(H_\bullet(\mathop{\bigvee}\limits_{i=1}^{2g}S^1)\otimes V) \ar[d]_{\varepsilon^{\bigvee_{i=1}^{2g}S^1}} \ar[r]^{p} & M\mathop{\otimes}\limits_{S(V)} S(H_\bullet(\Sigma^g)\otimes V) \ar[d]_{\varepsilon^{\Sigma^g}}\ar[r]^{q}&  M\mathop{\otimes}\limits_{S(V)} S(H_\bullet(S^2)\otimes V) \ar[d]_{\varepsilon^{S^2}}\\ HH^{\bigvee_{i=1}^{2g}S^1}_\bullet(S(V),M)  \ar[r]^{(\bigvee_{i=1}^{2d} S^1 \hookrightarrow \Sigma^g)_\bullet}  & HH^{\Sigma^g_\bullet}_\bullet(S(V),M) \ar[r]^{(\Sigma^g \twoheadrightarrow S^2)_\bullet}  & HH^{S^2}_\bullet(S(V),M)}\end{eqnarray} with the vertical arrows being isomorphisms (of algebras if $M=S(V)$). 
Note that $S(H_\bullet(\Sigma^g)\otimes V)$ splits of as a tensor product $S(H_\bullet(\Sigma^g)\otimes V)\cong S(H_\bullet(\mathop{\bigvee}\limits_{i=1}^{2g}S^1)\otimes V) \otimes_{S(V))} S(H_\bullet(S^2)\otimes V) $. By Lemma~\ref{L:HKRwedge} and formula~\eqref{eq:epsilonwedge}, we already have an explicit morphism of algebras at the chain level for the restriction $S(H_\bullet(\mathop{\bigvee}\limits_{i=1}^{2g}S^1)\otimes V)\to CH_\bullet^{\Sigma^g_\bullet}(S(V),S(V))$ of  $\varepsilon^{\Sigma^g}$  to  $S(H_\bullet(\mathop{\bigvee}\limits_{i=1}^{2g}S^1)\otimes V)$. 
It is easy to check that the formula 
\begin{eqnarray}
\label{eq:epsilonsphere} \varepsilon^{S^2}(\sigma v)=\left(
\begin{matrix}
1&&\\ 
& \otimes 1& \otimes\, v \\
& \otimes 1& \otimes 1 
\end{matrix}
\right) -\left(
\begin{matrix}
1&&\\ 
& \otimes 1& \otimes 1 \\
& \otimes \, v & \otimes 1 
\end{matrix}
\right) 
\end{eqnarray} (where we use the notation $\left(
\begin{matrix}
a_{(0,0)}&&\\ 
& \otimes a_{(1,1)}& \otimes a_{(1,2)} \\
& \otimes a_{(2,1)}& \otimes a_{(2,2)} 
\end{matrix}
\right)$ for an homogeneous element in $CH^{S^2_2}_\bullet(S(V),S(V))$ as in Example~\ref{E:2sphere}) defines a cocycle in $CH^{S^2_2}_\bullet(S(V),S(V))$ and induces a  quasi-isomorphism of algebras $S(H_\bullet(S^2)\otimes V)\to CH^{S^2}_\bullet(S(V),S(V))$ (see the proof of Theorem~\ref{T:HKR} and~\cite{G}). For any  $v\in V$, choose any cocycle $\varepsilon^{\Sigma^g}(\sigma v) \in CH_\bullet^{\Sigma^g_2}(S(V),S(V))$ such that $\varepsilon^{\Sigma^g}(\sigma v)$ is maped to $\varepsilon^{S^2}(\sigma v)$ by the map $(\Sigma^g \twoheadrightarrow S^2)_\bullet$ in diagram~\eqref{eq:diagramd=0}.
We have thus defined an \emph{explicit} quasi-isomorphism of algebras 
 $\varepsilon^{\Sigma^g}:S(H_\bullet(\Sigma^g)\otimes V) \to CH^{\Sigma^g_\bullet}_\bullet(S(V),S(V))$ which, by abuse of notation, could be rewritten as the tensor product $\varepsilon^{\Sigma^g}= \varepsilon^{\mathop{\bigvee}\limits_{i=1}^{2g}S^1}\otimes_{S(V)}\varepsilon^{S^2}$ through the isomorphism $S(H_\bullet(\Sigma^g)\otimes V)\cong S(H_\bullet({\bigvee}_{i=1}^{2g}S^1)\otimes V) \otimes_{S(V))} S(H_\bullet(S^2)\otimes V) $.
\begin{remk} There is a standard choice for $\varepsilon^{\Sigma^g}(\sigma v)$, given as follows.  Recall, that $\Sigma^g_\bullet$ is obtained by gluing $g$ standard models for the square $I^2_\bullet$, and $g(g-1)$ models for triangles  $T_\bullet$. In particular any element in $I_2^2$ or $T_2$ is a sum of tensors which can be written in the form 
$\left(
\begin{matrix}
\mbox{boundary terms}&&\\ 
& \otimes a_{(1,1)}& \otimes a_{(1,2)} \\
& \otimes a_{(2,1)}& \otimes a_{(2,2)} 
\end{matrix}
\right)$
where the boundary terms are tensor powers of elements lying in the boudary $(\partial I^2)_\bullet$ $(\partial T)_2$. Thus, for $v\in V$, and any square or triangle $C_\bullet\subset \Sigma^g_\bullet$ we can define the element $\varepsilon_{C_\bullet} (v)= \left(
\begin{matrix}
1s&&\\ 
& \otimes 1& \otimes \, v \\
& \otimes 1& \otimes 1 
\end{matrix}
\right) -\left(
\begin{matrix}
1s&&\\ 
& \otimes 1& \otimes 1 \\
& \otimes \, v& \otimes 1 
\end{matrix}
\right)\in CH_\bullet^{\Sigma^g_2}(S(V),S(V))$   where the $1s$ in the top left corner means that any tensor in the boundary of $C_2$ or in $\Sigma^g_2 -C_2$ is 1.

\smallskip

 It is straightforward to check that, in the normalized chain complex, a possible choice for $\varepsilon^{\Sigma^g}(\sigma v)\in CH^{\Sigma^g_2}_{\bullet}(S(V),S(V))$ is given by $$\varepsilon^{\Sigma^g}(\sigma v)=\cfrac{1}{g^2}\sum_{C_\bullet \subset \Sigma^g_\bullet} \varepsilon_{C_\bullet}(v) $$ where the sum is over all triangles and squares in the simplicial model for $\Sigma^g$.
\end{remk}

\medskip

We now define a multiplication on  $\bigoplus_{g\geq 0}  Hom_{S(V)}\big( S(H_\bullet(\Sigma^g)\otimes V),S(V)\big )$.  

Consider the natural  surjective map $H_\bullet(\Sigma^g ) \bigoplus H_\bullet(\Sigma^h)\to H_\bullet(\Sigma^g \vee \Sigma^h)$ (whose kernel is isomorphic to $H_0(pt)$), tensor both sides with $V$ over the ground field and apply the free graded commutative algebra functor $S$.  Dualizing as $S(V)$-modules\footnote{When $X$ is an arcwise connected space, we endow  $S(H_\bullet(X)\otimes V)$ with its  natural  $S(V)\cong S(H_0(X)\otimes V)$-module structure.},  the multiplication $\mu:S(V)\otimes S(V)\to S(V)$ induces a linear map

\begin{multline*}Hom_{S(V)}\big( S(H_\bullet(\Sigma^g)\otimes V),S(V)\big ) \otimes Hom_{S(V)}\big( S(H_\bullet(\Sigma^h)\otimes V),S(V)\big )\\ \stackrel{\mu_*}\longrightarrow Hom_{S(V)}\big( S(H_\bullet(\Sigma^g\vee \Sigma^h)\otimes V),S(V)\big )\end{multline*}
for any $g,h\geqslant 0$. Furthermore, the pinching map $\pinch{g,h}$ yields a linear map $H_\bullet(\Sigma^{g+h}) \to H_\bullet(\Sigma^g\vee \Sigma^h)$  and thus an $S(V)$-algebra map 
\begin{equation}\label{eq:pinchd=0}
p_{g,h}: S(H_\bullet(\Sigma^{g+h})\otimes V)\to S(H_\bullet(\Sigma^g\vee \Sigma^h)\otimes V).\end{equation}

\begin{defn}\label{D:cupd=0} The multiplication $\cup$ on $\bigoplus_{g\geq 0}  Hom_{S(V)}\big( S(H_\bullet(\Sigma^g)\otimes V),S(V)\big )$ is induced by the composition
\begin{multline*}
Hom_{S(V)}\big( S(H_\bullet(\Sigma^g)\otimes V),S(V)\big ) \otimes Hom_{S(V)}\big( S(H_\bullet(\Sigma^h)\otimes V),S(V)\big )\\ \stackrel{\mu_*}\to Hom_{S(V)}\big( S(H_\bullet(\Sigma^g\vee \Sigma^h)\otimes V),S(V)\big )\\
\stackrel{p_{g,h}^*} \to Hom_{S(V)}\big( S(H_\bullet(\Sigma^{g+h})\otimes V),S(V)\big ).
\end{multline*}
\end{defn}
It is immediate to check that $\cup$ makes $\bigoplus_{g\geq 0}  Hom_{S(V)}\big( S(H_\bullet(\Sigma^g)\otimes V),S(V)\big )$ into an associative unital algebra.

\begin{remk} Let $(S(V),d)$ be a free graded commutative algebra with non-zero differential. It is easy to check that Definition~\ref{D:cupd=0} indeed yields a differential  graded unital algebra structure for $\bigoplus_{g\geq 0} \left( Hom_{S(V)}\big( S(H_\bullet(\Sigma^g)\otimes V),S(V)\big ), {d^{\Sigma^g}}^*\right)$. That is, the differential ${d^{\Sigma^g}}^*$ is a derivation for the multiplication $\cup$. 
\end{remk}

\medskip

The following theorem expresses that the surface product for $S(V)$ (with zero differential) corresponds to the multiplication $\cup$ in Definition~\ref{D:cupd=0}. 
\begin{thm}\label{T:HKRd=0} Let $S(V)$ be a  free graded commutative algebra  (with no differential).  
\begin{enumerate}
\item There is an isomorphism (natural in $V$),
$$\varepsilon^*= \bigoplus_{g\geq 0} {\varepsilon^{\Sigma^g}}^*:\bigoplus_{g\geq 0} HH_{\Sigma^g}^\bullet(S(V),S(V)) \stackrel{\sim}\to \bigoplus_{g\geq 0} Hom_{S(V)}\left(S(H_\bullet(\Sigma^g)\otimes V),S(V) \right)    $$ 
\item The following diagram is commutative
$$\xymatrix{ \left(\mathop{\bigoplus}\limits_{g\geq 0}HH_{\Sigma^g}^\bullet(S(V),S(V))\right)^{\otimes 2} \ar[r]^{\cup} \ar[d]^{\sim}_{({\varepsilon}^*)^{\otimes 2}}&  \mathop{\bigoplus}\limits_{g\geq 0}HH_{\Sigma^g}^\bullet(S(V),S(V))\ar[d]^{\sim}_{{\varepsilon}^*}
\\
\left(\mathop{\bigoplus}\limits_{g\geq 0}Hom_{S(V)}\left(S(H_\bullet(\Sigma^g)\otimes V),S(V) \right)\right)^{\otimes 2} \ar[r]^{\,\,\,\cup} & \mathop{\bigoplus}\limits_{g\geq 0}Hom_{S(V)}\left(S(H_\bullet(\Sigma^g)\otimes V),S(V) \right).}
$$\end{enumerate}
\end{thm}
\begin{proof}
Since $d^{\Sigma^g}=0$, the first statement follows from Corollary~\ref{C:HKR} and the definition of $\varepsilon^{\Sigma^g}$.

\medskip

By Corollary~\ref{C:pushoutalgebra}, there is a quasi-isomorphism of differential graded algebras $$CH_{\bullet}^{\Sigma^g_\bullet}(S(V),S(V)) \mathop{\otimes}\limits_{S(V)} CH_{\bullet}^{\Sigma^h_\bullet}(S(V),S(V)) \cong CH_{\bullet}^{(\Sigma^g\vee \Sigma^h)_\bullet}(S(V),S(V)).$$
Since the map $\varepsilon^{\Sigma^g}$ and $\varepsilon^{\Sigma^h}$ coincide on $S(V)$, the map \begin{multline*} 
\varepsilon^{\Sigma^g\vee \Sigma^h}: S(H_\bullet(\Sigma^g\vee \Sigma^h)\otimes V)\cong S(H_\bullet(\Sigma^g)\otimes V)\mathop{\otimes}\limits_{S(V)} S(H_\bullet(\Sigma^g)\otimes V) \\
\stackrel{\varepsilon^{\Sigma^g}\mathop{\otimes}\limits_{S(V)}\varepsilon^{\Sigma^h} } \longrightarrow CH_{\bullet}^{(\Sigma^g\vee \Sigma^h)_\bullet}(S(V),S(V))
\end{multline*} is well defined and an algebra quasi-isomorphism by Theorem~\ref{T:HKR}. Similarly, Corollary~\ref{C:pushoutalgebra} yields a natural quasi-isomorphism of algebras $$CH^{\Sigma^g_\bullet}_\bullet(S(V),S(V)) \otimes CH^{\Sigma^h_\bullet}_\bullet(S(V),S(V))\stackrel{\simeq}{\longrightarrow}   CH_{\bullet}^{(\Sigma^g\coprod \Sigma^h)_\bullet}(S(V),S(V))$$ and thus  $$\varepsilon^{\Sigma^g \coprod \Sigma^h}= \varepsilon^{\Sigma^g}\mathop{\otimes} \varepsilon^{\Sigma^h}:S(H_\bullet(\Sigma^g\coprod \Sigma^h)\otimes V)\to CH_{\bullet}^{(\Sigma^g\coprod \Sigma^h)_\bullet}(S(V),S(V))$$ is a quasi-isomorphism of algebras. Note that, for any algebra $A$ and module $M$ there is a natural isomorphism of cosimplicial modules
$$CH_\bullet^{(\Sigma^g \coprod \Sigma^h)_n}(A,M) \cong CH_\bullet^{\Sigma^g_n}(A,M) \times CH_\bullet^{\Sigma^h_n}(A,A) $$ with the diagonal simplicial structure on the right hand side, {\it cf.} Definition~\ref{cross}. Further the pointed  maps $j^g:pt_\bullet \to \Sigma^g_\bullet$ and $j^h:pt_\bullet \to \Sigma^h_\bullet$ yields a simplicial map $pt_\bullet\coprod pt_\bullet \stackrel{j^g\coprod j^h}\longrightarrow \Sigma^g_\bullet\coprod \Sigma^h_\bullet $ which in turn gives a structure of $A\otimes A \stackrel{\sim} \hookrightarrow CH^{pt_\bullet \coprod pt_\bullet }_\bullet(A,A)$-module to $ CH_{\bullet}^{(\Sigma^g\coprod \Sigma^h)_\bullet}(A,A)$. There is an isomorphism of simplicial modules $$CH_{\Sigma^g_\bullet}^\bullet(A,A) \times CH_{\Sigma^h_\bullet}^\bullet(A,A)\cong Hom_{A\otimes A}\big(CH_{\bullet}^{(\Sigma^g\coprod \Sigma^h)_\bullet}(A,A), A\otimes A\big)$$ under which the map $\vee: CH_{\Sigma^g_\bullet}^\bullet(A,A) \times CH_{\Sigma^h_\bullet}^\bullet(A,A)\to CH^\bullet_{(\Sigma^g \vee \Sigma^h)_\bullet}(A,A)$ from Definition~ \ref{D:cupg} identifies with the composition
\begin{multline*} Hom_{A\otimes A}\big(CH_{\bullet}^{(\Sigma^g\coprod \Sigma^h)_\bullet}(A,A), A\otimes A\big) \stackrel{ \mu_*} \to Hom_{A\otimes A}\big(CH_{\bullet}^{(\Sigma^g\coprod \Sigma^h)_\bullet}(A,A), A\big) \\ \cong  Hom_{A}\big(CH^{\bullet}_{(\Sigma^g\vee \Sigma^h)_\bullet}(A,A), A\big)=CH^{\bullet}_{(\Sigma^g\vee \Sigma^h)_\bullet}(A,A)\end{multline*}

It is now straightforward to check that the following diagram is commutative:
\begin{eqnarray} \label{eq:Diagmu}
\xymatrix{\\ CH_{\Sigma^g_\bullet}^\bullet(S(V),S(V)) \times CH_{\Sigma^h_\bullet}^\bullet(S(V),S(V))     \ar[r]^{\qquad \quad \vee}  \ar[d]^{\big(\varepsilon^{\Sigma^g\coprod \Sigma^h}\big)^*} & CH_\bullet^{\Sigma^g \vee \Sigma^h}(S(V),S(V)) \ar[d]^{\big(\varepsilon^{\Sigma^g\vee \Sigma^h}\big)^*} \\
Hom_{S(V)^{\otimes 2}}\big(S(H_\bullet(\Sigma^g\coprod \Sigma^h)\otimes V), S(V)\otimes S(V)\big) \ar[r]^{\qquad \mu_*}&
Hom_{S(V)}\big(S(H_\bullet(\Sigma^g\vee \Sigma^h)\otimes V), S(V)\big)  .} 
\end{eqnarray}
Let $x$ be any element in $S(H_\bullet(\Sigma^g)\otimes V)$ and $y$  be any element in $S(H_\bullet(\Sigma^h)\otimes V)$. Then, by definition $\varepsilon^{\Sigma^g\coprod \Sigma^h}(x\cdot y)=sh(\varepsilon^{\Sigma^g}(x), \varepsilon^{\Sigma^h}(y))$ where $sh$ is the shuffle product (see Section~\ref{section-product}). Since the Alexander-Whitney map is inverse to the shuffle product on normalized chains, it follows that the following diagram of $S(V)\otimes S(V)$-linear maps
\begin{eqnarray}
\label{eq:DiagAW}
\xymatrix{\\ CH_\bullet^{\Sigma^g}(S(V),S(V)) \otimes CH_\bullet^{\Sigma^h}(S(V),S(V))   & CH_\bullet^{\Sigma^g \coprod \Sigma^h}(S(V),S(V)) \ar[l]_{\qquad \qquad  AW}    \\
S(H_\bullet(\Sigma^g)\otimes V) \otimes S(H_\bullet(\Sigma^h) \otimes V)  \ar[u]^{\varepsilon^{\Sigma^g}\otimes {\varepsilon^{\Sigma^h}}}\ar[r]^{\qquad \sim} & S(H_\bullet(\Sigma^g\coprod \Sigma^h)\otimes V) \ar[u]^{{\varepsilon^{\Sigma^g\coprod \Sigma^h}}}.}  
\end{eqnarray}
is commutative on normalized chains.
Diagrams~\eqref{eq:Diagmu} and~\eqref{eq:DiagAW} imply that the following diagram is commutative,
\begin{eqnarray} 
\label{eq:DiagAWmu}
\xymatrix{ \\ \left(\mathop{\bigoplus}\limits_{g\geq 0}HH^\bullet_{\Sigma^g}(S(V),S(V))\right)^{\otimes 2}  \ar[r]^{\qquad \quad \vee \circ AW} \ar[d]_{{\varepsilon}^*\otimes {\varepsilon}^*} & \mathop{\bigoplus}\limits_{g,h\geq 0} HH^\bullet_{\Sigma^g \vee \Sigma^h}(S(V),S(V))   \ar[d]_{\mathop{\bigoplus}\limits_{g,h\geq 0}{\varepsilon^{\Sigma^g\vee \Sigma^h}}^*}  \\
\left(\mathop{\bigoplus}\limits_{g\geq 0} Hom_{S(V)}(S(H_\bullet(\Sigma^g)\otimes V),S(V))\right)^{\otimes 2}   \ar[r]^{\mu^*} & \mathop{\bigoplus}\limits_{g,h\geq 0} Hom_{S(V)}(S(H_\bullet(\Sigma^g\vee \Sigma^h)\otimes V), S(V)) }  
\end{eqnarray}

\smallskip

Statement (2) in Theorem~\ref{T:HKRd=0} now follows from the commutativity of diagram~\eqref{eq:DiagAWmu} and of the following diagram
 \begin{eqnarray}
\label{eq:DiagPinch}
\xymatrix{ S(H_\bullet(\Sigma^{g+h})\otimes V) \ar[r]^{p_{g,h}} \ar[d]^{\varepsilon^{\Sigma^{g+h}}} & S(H_\bullet(\Sigma^g\vee \Sigma^h)\otimes V)  \ar[d]^{\varepsilon^{\Sigma^g\vee \Sigma^h}}\\ 
CH^{\Sigma^{g+h}_\bullet}_\bullet(S(V),S(V))\ar[r]^{{\pinch{g,h}}_*} &   CH^{(\Sigma^g\vee \Sigma^h)_\bullet}_\bullet(S(V),S(V))}  
\end{eqnarray}
where $p_{g,h}$ is the map~\eqref{eq:pinchd=0} from Definition~\ref{D:cupd=0}. Since $S(H_\bullet(\Sigma^{g+h})\otimes V) $ is a free graded commutative algebra, and all the maps involved in Diagram~\eqref{eq:DiagPinch} are maps of algebras, it is enough to check the commutativity of Diagram~\eqref{eq:DiagPinch} on the generators. 
This is obvious for the generators lying in $H_{\bullet \leqslant 1}(\Sigma^{g+h})\otimes V$ since they are of simplicial degree 1.  
As for the generators lying in $H_{2}(\Sigma^{g+h})\otimes V$, by functoriality and the definition of $\varepsilon^{\Sigma^g}(\sigma v)$, it is sufficient to prove that the following diagram is commutative
$$ 
\xymatrix{ S(H_\bullet(S^2)\otimes V) \ar[rr]^{p} \ar[d]^{\varepsilon^{S^2}} && S(H_\bullet(S^2\vee S^2)\otimes V)  \ar[d]^{\varepsilon^{S^2\vee S^2}}\\ 
CH^{S^2_\bullet}_\bullet(S(V),S(V)) \ar[r]^{\mathcal{D}_\bullet(2)} & CH^{sd_2(S^2_\bullet)}_\bullet(S(V),S(V))\ar[r]^{{\pinch{0,0}}_*} &   CH^{(S^2\vee S^2)_\bullet}_\bullet(S(V),S(V))  }  
$$ (on the normalized chains)
Here, $p$ is the algebra map defined on the generators by the pinching map  on homology $$H_\bullet(S^2)\otimes V \stackrel{{\pinch{0,0}}_*\otimes id_V}\longrightarrow  H_\bullet(S^2\vee S^2)\otimes V.$$ Now the result follows from a straightforward computation.
\end{proof}

If $(S(V),d)$ is a free model of a differential graded algebra $(A,d_A)$, then by Proposition~\ref{P:cupinvariance}, there exists an algebra isomorphism $$\bigoplus_{g\geq 0} HH_{\Sigma^g}^\bullet(A,A) \cong \bigoplus_{g\geq 0} HH_{\Sigma^g}^\bullet (S(V),S(V)).$$ Further, $CH_{\Sigma^g_\bullet}^\bullet (S(V),S(V))$ is  a  filtred differential graded algebra with respect to the filtration induced by the internal degree. This yields a spectral sequence of algebras whose $(E_0^{\bullet,\bullet},d_0)$ page is the Hochschild cochain complex  $CH_{\Sigma^g_\bullet}^\bullet (S(V),S(V))$ of $S(V)$ \emph{equipped with the zero differential}. From Theorem~\ref{T:HKRd=0}, we easily deduce 
\begin{cor} Let $(S(V),d)$ be any free model of $(A,d_A)$. The $E^1$-term of the above spectral sequence is  $$E_1^{p,q}=\bigoplus_{g\geq 0} Hom_{S(V)}\left(S_{(q)}(H_\bullet(\Sigma^g)\otimes V),S(V) \right)$$ where the right hand side is equipped with the multiplication of Definition~\ref{D:cupd=0}, and $S_{(q)}(H_\bullet(\Sigma^g)\otimes V)$ consists of those polynomial of total external (homological) degree $q$ (that is the total degree coming from $H_\bullet(\Sigma^g)$ is q); in particular $q\geq 0$.
\end{cor}

\begin{ex}[Odd spheres] \label{E:oddspheres} Since spheres are formal, there is a (chain of) quasi-isomorphism of differential graded commutative algebras between the forms $\Omega^\bullet{S^{2n+1}}$ and $S(x)\cong H^\bullet(S^{2n+1})$, where $x$ is of degree $|x|=2n+1$. Applying Proposition~\ref{P:cupinvariance} and Theorem~\ref{T:HKRd=0}  we get that
$$HH_{\Sigma^g}^\bullet(\Omega^\bullet {S^{2n+1}},\Omega^\bullet {S^{2n+1}})\cong S(x^g,\alpha_1^g,\dots \alpha^g_g, \beta^g_1,\dots \beta^g_g,\omega^g)$$ where $|x^g|=2n+1, |\alpha_1^g|=\dots =|\alpha_g^g|=|\beta_1^g|=\dots=|\beta_g^g|=-2n$, and $|\omega^g|=1-2n$. The cup-product is given, for any polynomial $P=P(x^g,\alpha_i^g,\beta^g_j,\omega^g)\in S(x^g,\alpha_1^g,\dots \alpha^g_g, \beta^g_1,\dots \beta^g_g,\omega^g)$, by the formul\ae:
\begin{eqnarray*}
P(x^g,\alpha_i^g,\beta^g_j,\omega^g) \cup x^h &=& P(x^{g+h},\alpha_i^{g+h},\beta^{g+h}_j,\omega^{g+h})x^{g+h}, \\
P(x^g,\alpha_i^g,\beta^g_j,\omega^g) \cup \omega^h &=& P(x^{g+h},\alpha_i^{g+h},\beta^{g+h}_j,\omega^{g+h})\omega^{g+h}, \\
P(x^g,\alpha_i^g,\beta^g_j,\omega^g) \cup \alpha_i^h &=& P(x^{g+h},\alpha_i^{g+h},\beta^{g+h}_j,\omega^{g+h})\alpha_{g+i}^{g+h},\\
P(x^g,\alpha_i^g,\beta^g_j,\omega^g) \cup \beta_j^h &=& P(x^{g+h},\alpha_i^{g+h},\beta^{g+h}_j,\omega^{g+h})\beta_{g+j}^{g+h},
\end{eqnarray*}
where the products on the right hand side are taken in the free graded commutative algebra $S(x^{g+h},\alpha_1^{g+h},\dots \alpha^{g+h}_g, \beta^{g+h}_1,\dots \beta^{g+h}_g,\omega^{g+h})$. Note that the center of $\bigoplus_{g\geq 0} HH_{\Sigma^g}^\bullet(\Omega^\bullet {S^{2n+1}},\Omega^\bullet {S^{2n+1}})$ is exactly $HH_{\Sigma^0}^\bullet(\Omega^\bullet {S^{2n+1}},\Omega^\bullet{S^{2n+1}})\cong S(x^0,\omega^0)$. 

\medskip

By Theorem~\ref{T:surface=cup}, if $n\geq 1$, then\footnote{where, by convention, degrees are intended to be of cohomological type} $\mathbb{H}_\bullet(\Map(\Sigma^g,M))\cong HH^\bullet_{\Sigma^g}(\Omega^\bullet S^{2n+1},\Omega^\bullet S^{2n+1})$ and the surface product agrees with the cup product.
\end{ex}

\begin{ex}[Lie groups]\label{E:LieGroups}
It is well-known that if $G$ is a Lie group, then $G$ is rationally homotopy equivalent to a product $S^{2d_1+1}\times \cdots \times S^{2d_e+1} $ of odd spheres where $e$ is the exponent of the group. Thus, by Proposition~\ref{P:cupinvariance},  Theorem~\ref{T:HKRd=0}, and Example~\ref{E:oddspheres}, we find that
$$HH_{\Sigma^g}^\bullet(\Omega^\bullet G,\Omega^\bullet G)\cong  S(x_k^g,\alpha_{k,i}^g, \beta^g_{k,j},\omega_k^g),$$
where $k=1,\dots, e$, and $i,j=1, \dots, g$, and the degrees of the generators are given by $|x^g_k|=2d_k +1$, $|\alpha^g_{k,i}|=|\beta^g_{k,j}|=-2d_k$, and $|\omega^g_k|=1-2d_k$. The formulae for the cup product are similar to those in Example~\ref{E:oddspheres} (except for the additional subscript $k$).

\medskip

If $G$ is simply connected, then it is automatically 2-connected, and, by Theorem~\ref{T:surface=cup}, the surface product agrees with the cup product through the  isomorphism $\mathbb{H}_\bullet(\Map(\Sigma^g,G))\cong HH^\bullet_{\Sigma^g}(\Omega^\bullet G,\Omega^\bullet G)$.
\end{ex}

\end{document}